\documentclass[11pt,reqno]{amsart}

\usepackage[english]{babel}
\usepackage{mathtools,amssymb,amsthm,latexsym,amsmath,amstext}
\usepackage{xcolor}
\usepackage{titletoc}
\usepackage{enumitem}
\usepackage{mathrsfs}
\usepackage{commath} 
\usepackage{listings}
\usepackage{caption}
\usepackage{hhline}
\usepackage{tikz}
\usepackage{tikz-cd}
\usepackage{bbold} 
\usepackage{todonotes}
\usepackage{varioref} 
\usepackage{nameref}
\usepackage{amsthm,xpatch}
\usepackage{lipsum}
\usepackage{eucal} 
\usepackage{lipsum}
\usepackage{mathtools}
\usepackage{dsfont}
\usepackage{hyperref}

\newcommand \M{\mathfrak{M}}
\newcommand \Image{\mathrm{Im}}
\newcommand \Aff{\mathrm{Aff}}
\newcommand \kR{\mathfrak{R}}

\newcommand \Ker{\mathrm{Ker}}
\newcommand \X{\mathcal{X}}

\newcommand \kX{\mathfrak{X}}
\newcommand \N{\mathbb N}
\newcommand \Z{\mathbb Z}
\newcommand \Q{\mathbb Q}
\newcommand \R{\mathbb R}

\newcommand \Class{\mathcal C}

\DeclareMathOperator \Or {O}
\DeclareMathOperator \RCD {\mathrm{RCD}}
\DeclareMathOperator \GH {d_{\mathrm{GH}}}
\DeclareMathOperator \D {\mathfrak{D}}
\DeclareMathOperator \CD {\mathrm{CD}}
\DeclareMathOperator \Iso {Iso}

\DeclareMathOperator \Diam {Diam}
\DeclareMathOperator \Spt {Spt}

\DeclareMathOperator \Dis {Dis}
\DeclareMathOperator \id {id}
\DeclareMathOperator \dist {d}
\DeclareMathOperator \length {\mathcal{L}}
\DeclareMathOperator \m {\mathfrak{m}}

\DeclareMathOperator \Dpeq {\D_{\mathrm{p}}^{\mathrm{eq}}
}

\newtheorem{thmintro}{Theorem}

\newtheorem{corintro}{Corollary}

\theoremstyle{definition}
\newtheorem*{questionintro}{Question}
\theoremstyle{definition}
\theoremstyle{definition}

\theoremstyle{definition}

\theoremstyle{definition}

\theoremstyle{definition}

\theoremstyle{remark}
\newtheorem{remAppendix}{Remark}

\makeatletter
\newcommand{\proofpart}[2]{%
 \par
 \addvspace{\medskipamount}%
 \noindent \fbox{Part #1: #2}\\ \par\nobreak
 \addvspace{\smallskipamount}%
 \@afterheading
}
\makeatother

\theoremstyle{definition}
\newtheorem{Notation}{Notation}[section]
\theoremstyle{definition}

\theoremstyle{definition}

\theoremstyle{definition}
\newtheorem{Theorem}{Theorem}[section]
\theoremstyle{definition}
\newtheorem{cor}{Corollary}[section]
\theoremstyle{definition}
\newtheorem{Lemma}{Lemma}[section]
\theoremstyle{definition}
\newtheorem{Proposition}{Proposition}[section]
\theoremstyle{definition}
\newtheorem{Definition-Proposition}{Definition-Proposition}[section]
\theoremstyle{definition}
\newtheorem{Definition}{Definition}[section]
\theoremstyle{remark}
\newtheorem{Remark}{Remark}[section]

\begin{document}

\title{Moduli spaces of compact {RCD}({0},{N})-structures}

\author{Andrea Mondino}
\address{University of Oxford, Mathematical Institute}
\email{Andrea.Mondino@maths.ox.ac.uk}
\author{Dimitri Navarro}
\address{University of Oxford, Mathematical Institute}
\email{navarro@maths.ox.ac.uk}

\date{\today}

\maketitle

\begin{abstract}
The goal of the paper is to set the foundations and prove some topological results about moduli spaces of non-smooth metric measure structures with non-negative Ricci curvature in a synthetic sense (via optimal transport) on a compact topological space; more precisely, we study moduli spaces of $\RCD(0,N)$-structures. First, we relate the convergence of $\RCD(0,N)$-structures on a space to the associated lifts' equivariant convergence on the universal cover. Then we construct the Albanese and soul maps, which reflect how structures on the universal cover split, and we prove their continuity. Finally, we construct examples of moduli spaces of $\RCD(0,N)$-structures that have non-trivial rational homotopy groups.
\end{abstract}

\tableofcontents

\newpage

\section{Introduction}

One of Riemannian geometry's most fundamental problems is studying metrics satisfying a particular curvature constraint on a fixed smooth manifold. Three thoroughly studied types of curvature are the sectional, the Ricci, and the scalar curvature; common curvature constraints are lower (resp. upper) bounds on the corresponding curvature. Also, when a smooth manifold admits a metric of the desired type, it is interesting to describe such metrics' space. A way to tackle this problem is to study the topological properties of the associated moduli space, i.e., the quotient of the space of all metrics satisfying the curvature condition by isometry equivalence. In the last decades, moduli spaces of metrics with positive scalar curvature (resp. negative sectional curvature) have been studied intensively (see \cite{Tuschmann-Wraith_15} for a comprehensive introduction). Yet, there are not as many results on moduli spaces of non-negatively Ricci curved metrics. In 2017, Tuschmann and Wiemeler published the first result on these moduli spaces' homotopy groups (see Theorem 1.1 in \cite{Tuschmann-Wiemeler_17}).\\
Until recently, most of the results on moduli spaces in metric geometry focused on smooth metrics. Nevertheless, Belegradek \cite{Belegradek_17} recently tackled the case of non-negatively curved Alexandrov spaces, studying the moduli space of non-negatively curved length metrics on the $2$-sphere. Whereas Alexandrov introduced curvature lower bound in the setting of length spaces (generalizing sectional curvature lower bound), $\RCD$ spaces generalize lower bounds on the Ricci curvature to the setting of metric measure spaces. Indeed, roughly, $\RCD(0,N)$ spaces should be thought as possibly non smooth spaces with dimension bounded above by $N\in [1,\infty)$ and non-negative Ricci curvature, in a synthetic sense. The first goal of the present paper is to set the foundations for studying moduli spaces of $\RCD(0,N)$-structures. The question we will be studying in this paper is the following:

\begin{questionintro}Let $N\in[1,\infty)$, and let $X$ be a compact topological space that admits an $\RCD(0,N)$-structure. What can be said about the topology of the moduli space of $\RCD(0,N)$-structures on $X$?
\end{questionintro}

We will start by recalling in Section \ref{sec1.1} what is an $\RCD(0,N)$ space. Then, in Section \ref{sec1.2}, we will introduce $\RCD(0,N)$-structures on a fixed topological space, together with the associated moduli space. In Section \ref{sec1.3}, we will introduce the notions of lift and push-forward of an $\RCD(0,N)$-structure. Afterwards, in Section \ref{sec1.4}, we will present concisely the Albanese map and the soul map associated to a compact topological space that admits an $\RCD(0,N)$-structure. Finally, in Section \ref{sec1.5}, we will present our main results.

\subsection{RCD(0,N)-spaces}\label{sec1.1}

The story of $\RCD$ spaces has its roots in Gromov's precompactness Theorem (see Corollary 11.1.13 of \cite{Petersen_06}). The result states that sequences of compact Riemannian manifolds with a lower bound on the Ricci curvature, and an upper bound on both the dimension and the diameter are precompact in the Gromov--Hausdorff topology (GH topology for short). Since then, there has been much work to understand the properties of limits of such sequences, called Ricci limit spaces. In the early '00s, Cheeger and Colding published in \cite{Cheeger-Colding_I_97}, \cite{Cheeger-Colding_II_00} and \cite{Cheeger-Colding_III_00} an extensive study of the aforementioned spaces. One important observation (already noticed by Fukaya) is that, to retain good stability properties at the limit space, it is fundamental to keep track of the Riemannian measures' behaviour associated to the approximating sequence. Since then, it is common to use the measured Gromov--Hausdorff topology (mGH topology for short), endowing Riemannian manifolds with their normalized volume measure.\\
A related (but slightly different) approach is to introduce a new definition of Ricci curvature lower bounds and dimension upper bound, at the more general level of possibly non-smooth metric measure spaces, which generalizes the classical notions and is stable when passing to the limit in the mGH topology. The definition of $\RCD$ spaces is an example of such a definition. Therefore, any result proven for $\RCD$ spaces (with the tools of metric measure theory) would hold a fortiori for Ricci limit spaces.\\

Throughout the paper, we will use the following definition of metric measure spaces.

\begin{Definition}Let $(X,\dist,\m)$ be a triple where $(X,\dist)$ is a metric space and $\m$ is a measure on $X$. We say that $(X,\dist,\m)$ is a metric measure space (m.m.s. for short) when $(X,\dist)$ is a complete separable metric space and $\m$ is a non-negative boundedly finite Radon measure on $(X,\dist)$.
\end{Definition}

$\CD$-spaces were introduced independently by Lott and Villani in \cite{Lott-Villani_09}, and by Sturm in \cite{Sturm_I_06} and \cite{Sturm_II_06}. For simplicity, we only give the definition of $\CD(0,N)$-spaces (following Definition 1.3 in \cite{Sturm_II_06}). An extensive study of $\CD$-spaces is given in \cite[Chapters 29 and 30]{Villani_09}.

Denote by $\mathcal{P}_2(X,\dist,\m)$ the space of probability measures that are absolutely continuous w.r.t. $\m$ and with finite variance, and let $W_2$  be the quadratic Kantorocih-Wasserstein transportation distance. Let also $S_N(\mu\mid\m)$ be the Renyi entropy of $\mu\in \mathcal{P}_2(X,\dist,\m)$ with respect to $\m$, i.e. $S_N(\rho\m\mid \m)=-N \int_X \rho^{1-\frac{1}{N}} \m$.

\begin{Definition}\label{CD(0,N)}Given $N\in[1,\infty)$, a \emph{$\CD(0,N)$-space} is a m.m.s. $(X,\dist,\m)$ such that $\m(X)>0$ and satisfying the following property: for every pair $\mu_0,\, \mu_1 \in \mathcal{P}_2(X,\dist,\m)$, there exists a $W_2$-geodesic $(\mu_t)_{t\in [0,1]}\subset  \mathcal{P}_2(X,\dist,\m)$ from $\mu_0$ to $\mu_1$ such that the function $[0,1]\ni t \mapsto S_N(\mu_t \mid \m)$ is convex.
\end{Definition}

The class of $\CD$ spaces includes (some) non-Riemannian Finsler structures which, from the work of Cheeger-Colding \cite{Cheeger-Colding_I_97, Cheeger-Colding_II_00, Cheeger-Colding_III_00}, cannot appear as Ricci limits.  In order to single out the ``Riemannian'' $\CD$ structures, it is convenient to add the assumption that the Sobolev space $W^{1,2}$ is a Hilbert space  \cite{Ambrosio-Gigli-Savare_14, Gigli15} (see also \cite{AmbrosioGigliMondinoRajala15}): indeed   for a Finsler manifold, $W^{1,2}$ is a Banach space and $W^{1,2}$ is a Hilbert space if and only if the Finsler structure is actually Riemannian.  Such a  condition is known as ``infinitesimal Hilbertianity''.

\begin{Definition}Given $N\in[1,\infty)$, a m.m.s. $(X,\dist,\m)$ is an $\RCD(0,N)$ space if it is an infinitesimally Hilbertian $\CD(0,N)$-space.
\end{Definition}

\noindent
An equivalent way to characterise the $\RCD(0,N)$ condition is via the validity of the  Bochner inequality  \cite{AmbrosioGigilSavareAnnProb,ErbarKuwadaSturm15, AmbrosioMondinoSavare15}  (see also \cite{CavallettiMilman16} for the globalization of general Ricci lower bounds $\RCD(K,N)$).

\subsection{Moduli spaces of $\RCD(0,N)$-structures on compact topological spaces}\label{sec1.2}

In this section, we will start by defining $\RCD(0,N)$-structures on a fixed topological space; then, we will introduce the associated moduli space, together with its topology.

\begin{Definition}Given a topological space $X$ and $N\in[1,\infty)$, an \emph{$\RCD(0,N)$-structure on $X$} is an $\RCD(0,N)$ space $(X,\dist,\m)$ such that $\dist$ metrizes the topology of $X$ and $\Spt(\m)=X$. 
\end{Definition}

It is common to identify m.m.s. that are isomorphic. There are two distinct notions of isomorphisms for m.m.s. (see the discussion in Chapter 27, Section ``adding the measure'' in \cite{Villani_09}). However, both notions coincide when restricted to the set of $\RCD(0,N)$-structures on a topological space (since we imposed measures to have full support in that case). We will adopt the following definition of isomorphism between metric measure spaces.

\begin{Definition}\label{def:mms}Two m.m.s. $(X_1,\dist_1,\m_1)$ and $(X_2,\dist_2,\m_2)$ are \emph{isomorphic} when there is a bijective isometry $\phi\colon(X_1,\dist_1)\rightarrow(X_2,\dist_2)$ such that $\phi_*\m_1=\m_2$.
\end{Definition}



Thanks to \cite[Theorem 2.3]{Abraham-Delma-Hoscheit_13} and \cite[Remark 3.29]{Gigli-Mondino-Savare_15}, the following result holds.

\begin{Theorem}The compact Gromov--Hausdorff--Prokhorov distance $\dist_{\mathrm{GHP}}^{\mathrm{c}}$ (see \cite[Section 2.2]{Abraham-Delma-Hoscheit_13}) is a complete separable metric on the set $\kX$ of isomorphism classes of compact metric measure spaces. Moreover, $\dist_{\mathrm{GHP}}^{\mathrm{c}}$ metrizes the mGH topology (see for instance \cite[Definition 27.30]{Villani_09}).
\end{Theorem}

We are now in position to introduce the moduli space $\M_{0,N}(X)$ of $\RCD(0,N)$-structures on a compact topological space $X$; this will be the main object of study in the paper.

\begin{Notation}Let $N\in[1,\infty)$, and let $X$ be a compact topological space that admits an $\RCD(0,N)$-structure. We introduce the following spaces:
\begin{itemize}
\item[(i)] $\mathfrak{RCD}(0,N)\subset\kX$ is the set of isomorphism classes of compact $\RCD(0,N)$ spaces with full support, endowed with the mGH-topology (seen as a subspace of $\kX$),
\item[(ii)] $\mathfrak{R}_{0,N}(X)$ is the set of all $\RCD(0,N)$-structures on $X$,
\item[(iii)] $\M_{0,N}(X)$ is the quotient of $\mathfrak{R}_{0,N}(X)$ by isomorphisms, endowed with the mGH-topology (seen as a subspace of $\mathfrak{RCD}(0,N)$).\end{itemize}
We call $\M_{0,N}(X)$ the \emph{moduli space of $\RCD(0,N)$-structures on $X$}.
\end{Notation}

\subsection{Lift and push-forward}\label{sec1.3}

To our aims, a fundamental result is the existence of a universal cover for an $\RCD$ space (see \cite[Theorem 1.1]{Mondino-Wei_19}).

\begin{Theorem}\label{th3.2}Let $N\in[1,\infty)$, and let $X$ be a compact topological space that admits an $\RCD(0,N)$-structure. Then $X$ admits a universal cover $p\colon \tilde{X}\to X$. We denote $\overline{\pi}_1(X)$ the associated group of deck transformation, also called the \emph{revised fundamental group of $X$}.
\end{Theorem}

\begin{Remark}In Theorem \ref{th3.2}, the universal cover must be understood in a general sense (as defined in \cite[Chap. 2, Sect. 5] {Spanier_81}). To be precise, it is not known whether it is simply connected or not. In particular, the revised fundamental group $\overline{\pi}_1(X)$  is a quotient of the fundamental group $\pi_1(X)$ of $X$, and may a priori not be isomorphic to $\pi_1(X)$.
\end{Remark}

As in the case of a Riemannian manifold, it is possible to lift an $\RCD(0,N)$-structure on a compact topological space to its universal cover, and, conversely, to push-forward an equivariant $\RCD(0,N)$-structure on the universal cover back to the base space. Indeed, let us fix a real number $N\in[1,\infty)$, a compact topological space $X$ that admits an $\RCD(0,N)$-structure, and denote $p\colon \tilde{X}\to X$ the universal cover of $X$. We will see in Section \ref{section4} that:
\begin{itemize}
\item given an $\RCD(0,N)$-structure $(X,\dist,\m)$ on $X$, there exists a unique $\RCD(0,N)$-structure $(\tilde{X},\tilde{\dist},\tilde{\m})$ on $\tilde{X}$ (called the \emph{lift of $(X,\dist,\m)$}) such that $p\colon(\tilde{X},\tilde{\dist},\tilde{\m})\to (X,\dist,\m)$ is a local isomorphism, and $\overline{\pi}_1(X)$ acts by isomorphism on $(\tilde{X},\tilde{\dist},\tilde{\m})$ (see Corollary \ref{corlift});
\item given an $\RCD(0,N)$-structure $(\tilde{X},\tilde{\dist},\tilde{\m})$ on $\tilde{X}$ such that $\overline{\pi}_1(X)$ acts by isomorphism on $(\tilde{X},\tilde{\dist},\tilde{\m})$, there exists a unique $\RCD(0,N)$-structure $(X,\dist,\m)$ on $X$ (called the \emph{push-forward of $(\tilde{X},\tilde{\dist},\tilde{\m})$}) such that $p\colon(\tilde{X},\tilde{\dist},\tilde{\m})\to (X,\dist,\m)$ is a local isomorphism (see Proposition \ref{propX}).
\end{itemize}
Theorem \ref{thmA}, which will be introduced in Section \ref{sec1.5}, relates the convergence of $\RCD(0,N)$-structures on $X$ to the convergence of the associated lifts.

\subsection{Soul map and Albanese map}\label{sec1.4}

In this section, we fix a real number $N\in[1,\infty)$, a compact topological space $X$ that admits an $\RCD(0,N)$-structure, and denote $p\colon \tilde{X}\to X$ the universal cover of $X$. 

A special property enjoyed by $\RCD(0,N)$ spaces is the existence of splittings. More precisely, given an $\RCD(0,N)$-structure $(X,\dist,\m)$ on $X$, and denoting $(\tilde{X},\tilde{\dist},\tilde{\m})$ the associated lift (see Section \ref{sec1.3}), there exists (thanks to Theorem 1.3 in \cite{Mondino-Wei_19}) an isomorphism:
$$
\phi\colon(\tilde{X},\tilde{\dist},\tilde{\m})\to(\overline{X},\overline{\dist},\overline{\m})\times\R^k,
$$
where $k\in\N\cap[0,N]$ is called the \emph{degree of $\phi$}, $\R^k$ is endowed with the Euclidean distance and the Lebesgue measure, and $(\overline{X},\overline{\dist},\overline{\m})$ is a compact $\RCD(0,N-k)$-space such that $\overline{\pi}_1(\overline{X})=0$; the space $\overline{X}$ is called the \emph{soul of $\phi$} (see Theorem \ref{th3.2} for the definition of $\overline{\pi}_1$). Such a map $\phi$ is called a \emph{splitting of $(\tilde{X},\tilde{\dist},\tilde{\m})$}, and induces an isomorphism: 
$$
\phi_*\colon\Iso(\tilde{X},\tilde{\dist},\tilde{\m})\to\Iso(\overline{X},\overline{\dist},\overline{\m})\times\Iso(\R^k).
$$
Moreover, the revised fundamental group $\overline{\pi}_1(X)$ acts by isomorphism onto $(\tilde{X},\tilde{\dist},\tilde{\m})$. Therefore, applying $\phi_*$ and projecting onto $\Iso(\R^k)$, we get a group homomorphism $\rho_{\R}^{\phi}\colon\overline{\pi}_1(X)\to\Iso(\R^k)$. In Section \ref{sec5}, we will prove the following properties:
\begin{itemize}
\item the degree $k$ does not depend either on the chosen splitting $\phi$, or the chosen $\RCD(0,N)$-structure $(X,\dist,\m)$ on $X$ (see Corollary \ref{cor5.1});
\item the image $\Gamma_{\phi}\coloneqq\Image(\rho_{\R}^{\phi})$ is a crystallographic subgroup of $\Iso(\R^k)$; moreover, up to conjugating with an affine transformation, $\Gamma_{\phi}$ does not depend either on the chosen splitting $\phi$, or the chosen $\RCD(0,N)$-structure $(X,\dist,\m)$ on $X$ (see Proposition \ref{prop5.1}).
\end{itemize}
Therefore, it is possible to introduce $k(X)\coloneqq k$ (called the \emph{splitting degree of $X$}), and the set $\Gamma(X)$ of crystallographic subgroups of $\Iso(\R^k)$ that are conjugated to $\Gamma_{\phi}$ by an affine transformation (called the \emph{crystallographic class of $X$}), both being topological invariants of $X$.\\

Thanks to the above discussion, to any $\RCD(0,N)$-structure $(X,\dist,\m)$ on $X$ with lift $(\tilde{X},\tilde{\dist},\tilde{\m})$, and to any splitting $\phi$ of $(\tilde{X},\tilde{\dist},\tilde{\m})$, we can associate:
\begin{itemize}
\item the soul $(\overline{X},\overline{\dist},\overline{\m})$ of $\phi$ which is a compact $\RCD(0,N-k(X))$-space; 
\item the compact $k(X)$-dimensional flat orbifold $(\R^{k(X)}/\Gamma_{\phi},\dist_{{\Gamma}_\phi})$ with orbifold fundamental group equal to $\Gamma_{\phi}\in\Gamma(X)$ (see the discussion preceding Definition \ref{def7.2}).
\end{itemize}
We will denote $A(X)$ the set of compact flat $k(X)$-dimensional orbifolds whose orbifold fundamental group belong to $\Gamma(X)$ (called the \emph{Albanese class of $X$}).

\begin{Remark}Let $(M^N,g)$ be a compact $N$-dimensional Riemannian manifold with non-negative Ricci curvature, and such that $\pi_1(M)=\Z^k$. Observe that, in that case, $(M,\dist_{g},\m_g)$ is an $\RCD(0,N)$ space, where $\dist_g$ is the geodesic distance and $\m_g$ is the Riemannian measure. It is possible to show that the orbifold obtained following the discussion above is nothing but the usual Albanese variety of $(M,g)$ (up to isometry).
\end{Remark}

In Section \ref{sec7.2}, we will see that, up to isomorphism, the orbifold (resp. the soul) does not depend either on the choice of the splitting map $\phi$, or on the isomorphism class of $(X,\dist,\m)$ (see Lemma \ref{lem:albwelldef}). Therefore, we will be able to introduce:
\begin{itemize}
\item the Albanese map $\mathcal{A}\colon\M_{0,N}(X)\to\mathscr{M}_{\mathrm{flat}}(A(X))$ associated to $X$, where $\mathscr{M}_{\mathrm{flat}}(A(X))$ is the quotient of $A(X)$ by isometry equivalence (endowed with the GH topology);
\item the soul map $\mathcal{S}\colon\M_{0,N}(X)\to\mathfrak{RCD}(0,N-k(X))$ associated to $X$;
\end{itemize}
such that for any $\RCD(0,N)$-structure $(X,\dist,\m)$ on $X$ with lift $(\tilde{X},\tilde{\dist},\tilde{\m})$ and any splitting $\phi$ of $(\tilde{X},\tilde{\dist},\tilde{\m})$ with soul $(\overline{X},\overline{\dist},\overline{\m})$, we have $\mathcal{A}([X,\dist,\m])=[\R^{k(X)}/\Gamma_{\phi},\dist_{\Gamma_{\phi}}]$ and $\mathcal{S}([X,\dist,\m])=[\overline{X},\overline{\dist},\overline{\m}]$ (where brackets denote the equivalence class in the appropriate moduli space).\\

Theorem \ref{thmB}, which will be presented in Section \ref{sec1.5}, states the continuity of the Albanese and soul maps.

\subsection{Main results}\label{sec1.5}

Our first result relates the convergence of $\RCD(0,N)$-structures on a compact topological space to the convergence of the associated lifts.

\begin{thmintro}\label{thmA}Let $N\in[1,\infty)$, let $X$ be a compact topological space that admits an $\RCD(0,N)$-structure, and denote $p\colon\tilde{X}\to X$ the universal cover of $X$. Assume that for every $n\in\N\cup\{\infty\}$:
\begin{itemize}
\item $\X_n=(X,\dist_n,\m_n,*_n)$ is a pointed $\RCD(0,N)$-structure on $X$;
\item $\tilde{\X}_n=(\tilde{X},\tilde{\dist}_n,\tilde{\m}_n,\tilde{*}_n)$ is the associated pointed lift, where $\tilde{*}_n$ is any point in $p^{-1}(*_n)$.
\end{itemize}
Then $\{\X_n\}_{n\in\N}$ converges to $\X_{\infty}$ in the pmGH topology if and only if $\{\tilde{\X}_n\}_{n\in\N}$ converges to $\tilde{\X}_{\infty}$ in the equivariant pmGH topology (introduced in Section \ref{subsection6.3}).
\end{thmintro}

\begin{Remark}\label{rem:pointed or not}Note that, since $X$ is compact, it is also possible to formulate Theorem \ref{thmA} as follows (forgetting about the reference points in the base space):\\
Assume that for every $n\in\N\cup\{\infty\}$, $\X_n^*=(X,\dist_n,\m_n)$ is an $\RCD(0,N)$-structure on $X$ with lift $\tilde{\X}_n^*=(\tilde{X},\tilde{\dist}_n,\tilde{\m}_n)$. Then $\{\X_n^*\}_{n\in\N}$ converges to $\X_{\infty}^*$ in the mGH topology if and only if for every $n\in\N\cup\{\infty\}$, there exist $\tilde{*}_n\in\tilde{X}$ such that $\{(\tilde{\X}_n^*,\tilde{*}_n)\}_{n\in\N}$ converges to $(\tilde{\X}_{\infty}^*,\tilde{*}_{\infty})$ in the equivariant pmGH topology.
\end{Remark}

As we will observe at the end of Section \ref{section7.1},  Theorem \ref{thmA} implies the following corollary, which is particularly useful when computing the homeomorphism type of specific examples of moduli spaces (see for example the case of $\R\mathbb{P}^2$ in \cite{Navarro}).

\begin{corintro}\label{corA}Let $N\in[1,\infty)$, let $X$ be a compact topological space that admits an $\RCD(0,N)$-structure, and denote $p\colon\tilde{X}\to X$ the universal cover of $X$.\\
Then the lift map:
$$
p^*\colon\M^{\mathrm{p}}_{0,N}(X)\to \M_{0,N}^{\mathrm{p,eq}}(\tilde{X}),
$$
is a homeomorphism (introduced in Section \ref{section7.1}), where $\M_{0,N}^{\mathrm{p,eq}}(\tilde{X})$ and $ \M^{\mathrm{p}}_{0,N}(X)$ are respectively the moduli space of equivariant pointed $\RCD(0,N)$-structures on $\tilde{X}$ and the moduli space of pointed $\RCD(0,N)$-structures on ${X}$ (introduced in Section \ref{section6}).
\end{corintro}

Observe that it is more straightforward to obtain Corollary \ref{corA} by using Theorem \ref{thmA} than its equivalent version given in Remark \ref{rem:pointed or not}.

Our next result states the continuity of the Albanese map and the soul map associated to a compact topological space that admits an $\RCD(0,N)$-structure (with $N\in[1,\infty)$). On the first hand, this result is essential when computing the homeomorphism type of specific examples of moduli spaces (see for example the case of the M\"obius band $\mathbb{M}^2$ and the finite cylinder $\mathbb{S}^1\times[0,1]$ in \cite{Navarro}). On the other hand, the continuity of the Albanese map will be crucial in the proof of Theorem \ref{thmC}.

\begin{thmintro}\label{thmB}Let $N\in[1,\infty)$, and let $X$ be a compact topological space that admits an $\RCD(0,N)$-structure. Then, the \emph{Albanese map} $\mathcal{A}\colon\M_{0,N}(X)\to\mathscr{M}_{\mathrm{flat}}(A(X))$ and the \emph{soul map} $\mathcal{S}\colon\M_{0,N}(X)\to\mathfrak{RCD}(0,N-k(X))$ are continuous, where $\mathscr{M}_{\mathrm{flat}}(A(X))$ and $\mathfrak{RCD}(0,N-k(X))$ are respectively endowed with the GH and mGH topology.
\end{thmintro}

Let us recall that if $X$ is a compact topological space that admits an $\RCD(0,2)$-structure, then the moduli space $\M_{0,2}(X)$ is contractible (see Theorem 1.1 in \cite{Navarro}). Theorem \ref{thmC} should be put in contrast with that result since it shows that the topology of moduli spaces of $\RCD(0,N)$-structures is not always as trivial. Moreover, Theorem \ref{thmC} can also be seen as a non-smooth analogue of Theorem 1.1 in \cite{Tuschmann-Wiemeler_17}.

\begin{thmintro}\label{thmC}Let $N\in[1,\infty)$ and let $X$ be a compact topological space that admits an $\RCD(0,N)$-structure such that $\overline{\pi}_1(X)=0$ (see Theorem \ref{th3.2} for the definition of $\overline{\pi}_1(X)$). In addition, let $Y$ be either $\mathbb{S}^1\times\mathbb{K}^2$ (where $\mathbb{K}^2$ is the Klein bottle) or a torus of dimension $k\geq 4$ such that $k\neq8,9,10$. Then, the moduli space $\M_{0,N+\dim(Y)}(X\times Y)$ has non-trivial higher rational homotopy groups.
\end{thmintro}

Thanks to Theorem \ref{thmC}, we immediately obtain the following corollary, which can be seen as a non-smooth analogue of Corollary 1.2 in \cite{Tuschmann-Wiemeler_17}.

\begin{corintro}\label{corB}For every $N\geq3$ (resp. $N\geq4$ / $N\geq5$) there exists a compact topological space $X$ such that $\M_{0,N}(X)$ is not simply connected (resp. has non-trivial third rational homotopy group / non-trivial fifth rational homotopy group).
\end{corintro}

In Section \ref{preli}, we will introduce in full details the main objects and constructions of the paper. In Section \ref{proof of main}, we will prove the main results.

\subsection*{Acknowledgments} Both the authors are supported by the European Research Council (ERC), under the European Union Horizon 2020 research and innovation programme, via the ERC Starting Grant CURVATURE, grant agreement No. 802689.
\\The authors wish to thank G\'erard Besson for stimulating conversations on the topics of the paper. 


\section{Preliminaries}\label{preli}

Throughout this section:
\begin{itemize}
\item $N\in[1,\infty)$ is a real number,
\item $X$ is a compact topological space that admits an $\RCD(0,N)$-structure,
\item $p\colon \tilde{X}\to X$ denotes the universal cover of $X$ (whose existence is given by Theorem \ref{th3.2}).
\end{itemize}

In Section \ref{section4}, we will introduce the notions of lift and push forward of an $\RCD(0,N)$-structure on $X$.\\
In Section \ref{sec5}, we will present splittings and use them to construct topological invariants associated to $X$ (splitting degree, crystallographic class and Albanese class).\\
In Section \ref{section6}, we will define the moduli space of pointed $\RCD(0,N)$-structures on $X$; then, we will introduce the moduli space of equivariant pointed $\RCD(0,N)$-structures on the universal cover $\tilde{X}$.\\
In Section \ref{sec7}, we will define the lift and push-forward map (which are important to get Corollary \ref{corA}), and the Albanese and soul maps.

\subsection{Covering space theory of {RCD}({0},{N})-spaces}\label{section4}





In this section, we will start by introducing $\delta$-covers associated to an $\RCD(0,N)$-structure on $X$. Then, we will explain how to lift an $\RCD(0,N)$-structure on $X$ to the associated $\delta$-cover. Afterwards, we will explain how the universal cover of $X$ is related to $\delta$-covers. Subsequently, we will explain how to lift an $\RCD(0,N)$-structure on $X$ to its universal cover $\tilde{X}$, and, conversely, how to push-forward an equivariant $\RCD(0,N)$-structure on $\tilde{X}$ onto $X$. Finally, we will introduce the Dirichlet domain associated to an $\RCD(0,N)$-structure on $X$.\\

Before introducing $\delta$-covers, we recall the following result (Chapter 2, Sections 4 and 5 of \cite{Spanier_81}) that associates a regular covering $p^{\mathcal{U}}\colon X^{\mathcal{U}}\to X$ to any open cover $\mathcal{U}$ of $X$.

\begin{Proposition}\label{opencover}Given an open cover $\mathcal{U}$ of $X$, there exists a unique regular covering $p^{\mathcal{U}}\colon X^{\mathcal{U}}\to X$ (up to equivalence) such that:
$$
\forall y \in X^{\mathcal{U}}, p^{\mathcal{U}}_*\pi_1(X^{\mathcal{U}},y)=\pi_1(\mathcal{U},p^{\mathcal{U}}(y)),
$$
where $\pi_1(\mathcal{U},p^{\mathcal{U}}(y))$ is composed of homotopy classes of loops of the form $\omega^{-1}*\alpha*\omega$, where $\alpha$ is a loop contained in some $U\in\mathcal{U}$ and $\omega$ is a path from $p^{\mathcal{U}}(y)$ to $\alpha(0)$. Moreover, every connected open set $U\in\mathcal{U}$ is evenly covered bu $p_{\mathcal{U}}$.
\end{Proposition}

The notion of $\delta$-cover was introduced first by Sormani and Wei to prove the existence of a universal cover for Ricci limit spaces (see Theorem 1.1 in \cite{Wei-Sormani_01}). Later, it has also been used by Mondino and Wei in \cite{Mondino-Wei_19} to prove Theorem \ref{th3.2}. These covering spaces will be very important in the proof of Theorem \ref{thmA}.

\begin{Definition}Given $\delta>0$ and $(X,\dist,\m)$ an $\RCD(0,N)$-structure on $X$, the \emph{$\delta$-cover associated to $(X,\dist,\m)$} is the regular covering $p_{\dist}^{\delta}\colon X^{\delta}_{\dist}\to X$ associated to the open covering $\mathcal{U}(\delta,\dist)$ consisting of balls of radius $\delta$ for the distance $\dist$ (see Proposition \ref{opencover}). We write $G(\delta,\dist)$ the associated group of deck transformations.
\end{Definition}

In the following result, we introduce the lift of an $\RCD(0,N)$-structure on $X$ to a $\delta$-cover.

\begin{Proposition}\label{prop3.4}Given $\delta>0$ and $(X,\dist,\m)$ an $\RCD(0,N)$-structure on $X$, there exists a unique $\RCD(0,N)$-structure $(X_{\dist}^{\delta},\dist_{\delta},\m_{\delta})$ on $X_{\dist}^{\delta}$ such that $p_{\dist}^{\delta}\colon(X_{\dist}^{\delta},\dist_{\delta},\m_{\delta})\to(X,\dist,\m)$ is a local isomorphism. Moreover, we have the following properties:
\begin{itemize}
\item[(i)] for every $\tilde{x},\tilde{y}\in X_{\dist}^{\delta}$, we have $\dist_{\delta}(\tilde{x},\tilde{y})=\inf\{\length_{\dist}(p_{\dist}^{\delta}\circ\tilde{\gamma})\}$, where the infimum is taken over all continuous path $\tilde{\gamma}\colon[0,1]\to X_{\dist}^{\delta}$ from $\tilde{x}$ to $\tilde{y}$ and $\length_{\dist}$ is the length structure induced by $\dist$,
\item[(ii)] for every Borel subset $\tilde{E}\subset X_{\dist}^{\delta}$ such that ${p_{\dist}^{\delta}}_{\lvert \tilde{E}}$ is an isometry, we have $\m_{\delta}(\tilde{E})=\m(E)$,
\item[(iii)] the group of deck transformations $G(\delta,\dist)$ is a subgroup of $\Iso_{\mathrm{m.m.s.}}(X_{\dist}^{\delta},\dist_{\delta},\m_{\delta})$,
\item[(iv)] for every $\tilde{x}\in X_{\dist}^{\delta}$ and every $r\leq\delta$, the restriction of $p_{\dist}^{\delta}$ to $B_{\dist_{\delta}}(\tilde{x},r)$ is a homeomorphism onto $B_{\dist}(p_{\dist}^{\delta}(\tilde{x}),r)$,
\item[(v)] for every $\tilde{x}\in X_{\dist}^{\delta}$ and every $r\leq\delta/2$, the restriction of $p_{\dist}^{\delta}$ to $(B_{\dist_{\delta}}(\tilde{x},r),\dist_{\delta},\m_{\delta})$ is an isomorphism onto $(B_{\dist}(x,r),\dist,\m)$.
\end{itemize}
\end{Proposition}

\begin{proof}First of all, there is obviously at most one $\RCD(0,N)$-structure on $X_{\dist}^{\delta}$ such that $p_{\dist}^{\delta}$ is a local isomorphism.\\
Then, thanks to Lemma 2.18 of \cite{Mondino-Wei_19}, $(X_{\dist}^{\delta},\dist_{\delta},\m_{\delta})$ is an $\RCD(0,N)$ space (where $\dist_{\delta}$ and $\m_{\delta}$ are defined as in (i) and (ii)). Moreover, it is readily checked that $\dist_{\delta}$ metrizes the topology of $X_{\dist}^{\delta}$, that $\Spt(\m_{\delta})=X_{\dist}^{\delta}$, and that $G(\delta,\dist)$ acts by isomorphism. Therefore, $(X_{\dist}^{\delta},\dist_{\delta},\m_{\delta})$ is an $\RCD(0,N)$-structure on $X_{\dist}^{\delta}$ satisfying point (i) to (iii).\\
Finally, thanks to Proposition 15 of \cite{Reviron_08}, and by definition of $\m_{\delta}$, point (iv) and (v) are satisfied.
\end{proof}

Now, we put the universal cover of $X$ in relation with $\delta$-covers (see Theorem 2.7 in \cite{Mondino-Wei_19} for a proof).

\begin{Theorem}\label{3.3}Let $(X,\dist,\m)$ be an $\RCD(0,N)$-structure on $X$, and let $\delta(X,\dist)$ be the supremum of all $\delta>0$ such that every ball of radius $\delta$ in $(X,\dist)$ is evenly covered by $p$. Then $\delta(X,\dist)>0$, and for every $\delta<\delta(X,\dist)$, $p$ and $p^{\delta}_{\dist}$ are equivalent, and every equivalence map is an isomorphism between $(\tilde{X},\tilde{\dist},\tilde{\m})$ and $(X_{\dist}^{\delta},\dist_{\delta},\m_{\delta})$.
\end{Theorem}

Thanks to Proposition \ref{prop3.4} and Theorem \ref{3.3}, we can introduce the lift of an $\RCD(0,N)$-structure on $X$ to the universal cover $\tilde{X}$.

\begin{cor}\label{corlift}Let $(X,\dist,\m)$ be an $\RCD(0,N)$-structure on $X$. There is a unique $\RCD(0,N)$-structure $(\tilde{X},\tilde{\dist},\tilde{\m})$ on $\tilde{X}$ (called the \emph{lift of $(X,\dist,\m)$}) such that $p\colon(\tilde{X},\tilde{\dist},\tilde{\m})\rightarrow (X,\dist,\m)$ is a local isomorphism. Moreover, the revised fundamental group $\overline{\pi}_1(X)$ acts by isomorphism on $(\tilde{X},\tilde{\dist},\tilde{\m})$.
\end{cor}

The following proposition is a sort of converse to Corollary \ref{corlift}; it introduces the push-forward of an equivariant $\RCD(0,N)$-structure on $\tilde{X}$ (cf. \cite[Lemma 2.18]{Mondino-Wei_19} and \cite[Lemma 2.24]{MMP}).

\begin{Proposition}\label{propX}Let $(\tilde{X},\tilde{\dist},\tilde{\m})$ be an $\RCD(0,N)$-structure on $\tilde{X}$ such that $\overline{\pi}_1(X)$ acts by isomorphisms on $(\tilde{X},\tilde{\dist},\tilde{\m})$. There is a unique $\RCD(0,N)$-structure $(X,\dist,\m)$ on $X$ (called the \emph{push-forward of $(\tilde{X},\tilde{\dist},\tilde{\m})$}) such that $p\colon (\tilde{X},\tilde{\dist},\tilde{\m})\to(X,\dist,\m)$ is a local isomorphism. It satisfies the following properties:
\begin{itemize}
\item[(i)] for every $x,y\in X$, we have $\dist(x,y)=\inf\{\tilde{\dist}(\tilde{x},\tilde{y})\}$, where the infimum is taken over all $\tilde{x}\in p^{-1}(x)$ and $\tilde{y}\in p^{-1}(y)$,
\item[(ii)] for every open set $U\subset X$ that is evenly covered by $p$, we have $\m(U)= \tilde{\m}(\tilde{U})$, where $\tilde{U}$ is any open set in $\tilde{X}$ such that $p\colon\tilde{U}\to U$ is a homeomorphism.
\end{itemize}
\end{Proposition}

\begin{proof}First of all, there is obviously at most one $\RCD(0,N)$-structure on $X$ such that $p$ is a local isomorphism.\\
Then, let us define $\dist$ and $\m$ as in points (i) and (ii). Observe that since $\tilde{X}$ is locally compact, and since $\overline{\pi}_1(X)$ acts by isometries, $\m$ is well defined, and the infimum in the definition of $\dist$ is achieved. It is then readily checked that $\dist$ is a distance on $X$, and that $\m$ defines a measure on $X$ (using the fact that the Borel $\sigma$-algebra of $X$ is generated by evenly covered open sets).\\
Let us now show that $p$ is a local isomorphism. Let $\tilde{x}\in\tilde{X}$ and define $x\coloneqq p(\tilde{x})$. There exists an open neighborhood $\tilde{U}$ of $\tilde{x}$ such that $p\colon \tilde{U}\to U\coloneqq p(\tilde{U})$ is a homeomorphism. Moreover, there exists $r>0$ such that $B_{\tilde{\dist}}(\tilde{x},r)\subset \tilde{U}$. Let us show that, for every $0<r'\leq r$, $p$ is a homeomorphism from $B_{\tilde{\dist}}(\tilde{x},r')$ onto $B_{\dist}(x,r')$. First, notice that $p$ is distance decreasing; in particular, we have $p(B_{\tilde{\dist}}(\tilde{x},r'))\subset B_{\dist}(x,r')$. Now, let $y\in B_{\dist}(x,r')$. Since the infimum in the definition of $\dist$ is achieved, there exists $\tilde{y}\in p^{-1}(y)$ such that $\tilde{\dist}(\tilde{x},\tilde{y})=\dist(x,y)<r'$. Hence, $p(B_{\tilde{\dist}}(\tilde{x},r'))=B_{\dist}(x,r')$. Since $B_{\tilde{\dist}}(\tilde{x},r')$ is a subset of $\tilde{U}$, $p$ is injective on $B_{\tilde{\dist}}(\tilde{x},r')$. Hence, $p$ is a bijective map from $B_{\tilde{\dist}}(\tilde{x},r')$ onto $B_{\dist}(x,r')$. However, $p$ is an open map, hence it is a homeomorphism from $B_{\tilde{\dist}}(\tilde{x},r')$ onto $B_{\dist}(x,r')$.\\
Now, let $\tilde{y},\tilde{z}\in B_{\tilde{\dist}}(\tilde{x},r/3)$. Looking for a contradiction, let us suppose that $\dist(y,z)<\tilde{\dist}(\tilde{y},\tilde{z})$, where $y\coloneqq p(\tilde{y})$ and $z\coloneqq p(\tilde{z})$. In that case, there exists $\tilde{z}'\in p^{-1}(z)$ such that $\tilde{\dist}(\tilde{y},\tilde{z}')=\dist(y,z)<\tilde{\dist}(\tilde{y},\tilde{z})\leq\tilde{\dist}(\tilde{x},\tilde{y})+\tilde{\dist}(\tilde{x},\tilde{z})<2r/3<r$. However, $p$ is a homeomorphism from $B_{\tilde{\dist}}(\tilde{x},r)$ onto $B_{\dist}(x,r)$, so we should have $\tilde{z}'=\tilde{z}$, which is the contradiction we were looking for. Hence, $p$ is an isometry from $B_{\tilde{\dist}}(\tilde{x},r/3)$ onto $B_{\dist}(x,r/3)$. Moreover, by definition of $\m$, this implies that $p$ is an isomorphism of metric measure space from $B_{\tilde{\dist}}(\tilde{x},r/3)$ onto $B_{\dist}(x,r/3)$.\\
Then, it is easy to check that $\dist$ metrizes the topology of $X$ and that $\Spt(\m)=X$. To conclude, we just need to show that $(X,\dist,\m)$ is an $\RCD(0,N)$ space.
To this aim, first of all observe that  $\RCD(0,N)$ is equivalent to $\RCD^{*}(0,N)$  (by the explicit form of the distortion coefficients), which in turn is equivalent to $\CD^{e}(0,N)$ plus infinitesimally Hilbertianity.

 Observe that, since $p$ is a local isometry, it preserves the length of curves; therefore, $(X,\dist)$ is a compact geodesic space. Moreover, since $X$ is compact and $\tilde{\m}$ is boundedly finite, $\m$ is necessarily a finite measure on $X$. Summarising: $(X,\dist,\m)$ is a compact geodesic space endowed with a finite measure, and it is locally isomorphic to an $\RCD(0,N)$ space  in the sense that for every point $x\in X$ there exists a closed metric ball $\bar{B}_r(x)$ centred at $x$ isomorphic to a closed metric ball $\tilde{\bar{B}}_{r}(\tilde{x})$ inside the $\RCD(0,N)$ space $\tilde{X}$. 

Notice that, by  triangle inequality,  if $y,z\in \bar{B}_{r/4}(x)$ then any geodesic joining them is contained in  $\bar{B}_{r}(x)$. Recall also that, given two absolutely continuous probability measures with compact support in an $\RCD$ space, there exists a unique $W_{2}$-geodesic joining them \cite[Theorem 1.1]{GigliRajalaSturm16}. It follows that, given two  absolutely continuous probability measures with compact support contained in  $\bar{B}_{r/4}(x)$ (which in turm is isomorphic to $\tilde{\bar{B}}_{r/4}(x)\subset \tilde{X}$, and $\tilde{X}$ satisfies $\RCD(0,N)$), there exists a unique $W_{2}$-geodesic joining them,  its support is contained in $\bar{B}_r(x)$, and it satisfies the convexity property of the $\CD^{e}(0,N)$ condition. 

 In particular, $(X,\dist, \m)$ satisfies the strong $\CD^{e}_{loc}(0,N)$  condition in the sense of \cite{ErbarKuwadaSturm15} and it is locally infinitesimally Hilbertian. Then, using \cite[Theorem 3.25]{ErbarKuwadaSturm15}, we obtain that  $(X,\dist, \m)$ satisfies $\RCD(0,N)$. 
\end{proof}

\begin{Remark}\label{rem:pushlift}Observe that if $(X,\dist,\m)$ is an $\RCD(0,N)$-structure on $X$, then the push-forward of the lift of $(X,\dist,\m)$ is equal to $(X,\dist,\m)$, thanks to Proposition \ref{propX} and Corollary \ref{corlift}. The same is true in the other direction; if $(\tilde{X},\tilde{\dist},\tilde{\m})$ is an $\RCD(0,N)$-structure on $\tilde{X}$ such that $\overline{\pi}_1(X)$ acts by isomorphisms, then the lift of the push-forward of $(\tilde{X},\tilde{\dist},\tilde{\m})$ is equal to $(\tilde{X},\tilde{\dist},\tilde{\m})$.
\end{Remark}

We conclude this section with the following results that introduces the Dirichlet domain associated to an $\RCD(0,N)$-structure on $X$.

\begin{Proposition}\label{prop3.6}Let $(X,\dist,\m)$ be an $\RCD(0,N)$-structure on $X$ and let $\tilde{x}\in\tilde{X}$. We define the \emph{Dirichlet domain with center $\tilde{x}$ associated to $(X,\dist,\m)$} by:
$$
\mathcal{F}(\tilde{x})\coloneqq\bigcap\limits_{\eta\in\overline{\pi}_1(X)}\phi_{\eta}^{-1}(\R_{\geq0})
$$
where $\phi_{\eta}(\tilde{y})\coloneqq\tilde{\dist}(\tilde{y},\eta\tilde{x})-\tilde{\dist}(\tilde{y},\tilde{x})$, for $\tilde{y}\in\tilde{X}$. The Dirichlet domain satisfies the following two properties:
\begin{itemize}
\item[(i)] for every $\tilde{y}\in\tilde{X}$, there exists $\eta\in\overline{\pi}_1(X)$ such that $\eta\tilde{y}\in\mathcal{F}(\tilde{x})$,
\item[(ii)] for every $\tilde{y}\in\mathcal{F}(\tilde{x})$, we have $\tilde{\dist}(\tilde{x},\tilde{y})=\dist(x,y)$, where $x\coloneqq p(\tilde{x})$ and $y\coloneqq p(\tilde{y})$.
\end{itemize}
In particular, $\mathcal{F}(\tilde{x})\subset B_{\tilde{\dist}}(\tilde{x},D)$, where $D\coloneqq\Diam(X,\dist)$.
\end{Proposition}

\begin{proof}We start with the proof of (i). Let $\tilde{y}\in\tilde{X}$ and define $R\coloneqq\tilde{d}(\tilde{x},\tilde{y})$. Then, $p^{-1}(x)\cap \overline{B}_{\tilde{d}}(\tilde{y},R)$ is a compact, discrete, non empty set; hence, it contains finitely many points. In particular, there exists $\eta\in\overline{\pi}_1(X)$ such that $\eta\tilde{x}\in \overline{B}_{\tilde{d}}(\tilde{y},R)$, and such that:
\begin{equation}\label{eqj}
\forall\tilde{z}\in p^{-1}(x)\cap \overline{B}_{\tilde{d}}(\tilde{y},R),\tilde{\dist}(\tilde{y},\eta\tilde{x})\leq\tilde{\dist}(\tilde{y},\tilde{z}).
\end{equation} 
Now, assume that $\mu\in\overline{\pi}_1(X)$. If $R\leq\tilde{\dist}(\tilde{y},\mu\tilde{x})$, we have $\tilde{\dist}(\tilde{y},\eta\tilde{x})\leq\tilde{\dist}(\tilde{y},\mu\tilde{x})$ since $\tilde{\dist}(\tilde{y},\eta\tilde{x})\leq R$. If $\tilde{\dist}(\tilde{y},\mu\tilde{x})<R$, then, thanks to equation \ref{eqj}, we also have $\tilde{\dist}(\tilde{y},\eta\tilde{x})\leq\tilde{\dist}(\tilde{y},\mu\tilde{x})$. Thus, for every $\mu\in\overline{\pi}_1(X)$, we get $
\tilde{\dist}(\tilde{y},\eta\tilde{x})\leq\tilde{\dist}(\tilde{y},\mu\tilde{x})$. Hence, for every $\mu\in\overline{\pi}_1(X)$, we have $\phi_{\mu}(\eta^{-1} \tilde y)=\tilde{\dist}(\eta^{-1} \tilde{y},\mu\tilde{x})-\tilde{\dist}(\eta^{-1} \tilde{y},\tilde{x})=\tilde{\dist}(\tilde{y},\eta\mu\tilde{x})-\tilde{\dist}(\tilde{y},\eta\tilde{x})\geq 0$. In conclusion, $\eta^{-1}\tilde{y}\in\mathcal{F}(\tilde{x})$.\\
Now we prove (ii). Assume that $\tilde{y}\in\mathcal{F}(\tilde{x})$ and let $\tilde{\gamma}\colon[0,1]\to\tilde{X}$ be minimizing geodesic from $\tilde{x}$ to $\tilde{y}$. Then, we define $y\coloneqq p(\tilde{y})$ and we assume that $\beta\colon[0,1]\to X$ is a minimizing geodesic from $x$ to $y$. Let $\tilde{\beta}$ be the lift of $\beta$ starting at $\tilde{x}$ and let $\eta\in\overline{\pi}_1(X)$ such that $\tilde{\beta}(1)=\eta \tilde{y}$. Looking for a contradiction, let us suppose that $\dist(x,y)<\tilde{\dist}(\tilde{x},\tilde{y})$. Then, observe that $\tilde{\dist}(\tilde{x},\eta\tilde{y})\leq\length(\beta)=\dist(x,y)$; in particular, $\tilde{\dist}(\tilde{x},\eta\tilde{y})=\dist(x,y)$, since $p$ contracts distances. Hence, we have $\phi_{\eta^{-1}}(\tilde{y})=\tilde{\dist}(\tilde{y},\eta^{-1}\tilde{x})-\tilde{\dist}(\tilde{x},\tilde{y})=\tilde{\dist}(\eta\tilde{y},\tilde{x})-\tilde{\dist}(\tilde{x},\tilde{y})=\dist(x,y)-\tilde{\dist}(\tilde{x},\tilde{y})<0$. In particular, $\tilde{y}\notin\mathcal{F}(\tilde{x})$, which is the contradiction we were looking for. Thus, $\dist(x,y)\geq\tilde{\dist}(\tilde{x},\tilde{y})$, and, since $p$ contracts distances, we have $\dist(x,y)=\tilde{\dist}(\tilde{x},\tilde{y})$. This concludes the proof.
\end{proof}

\subsection{Splittings and topological invariants}\label{sec5}


In this section, we will introduce the notion of splitting associated to an $\RCD(0,N)$-structure on $X$. To any splitting $\phi$, we will associate a degree $k$ and a euclidean homomorphism $\rho^{\phi}_{\R}\colon\overline{\pi}_{1}(X)\to \Iso(\R^k)$, and we will investigate the properties of $\Gamma(\phi)=\Image(\rho^{\phi}_{\R})$. We will prove that the degree $k$ and the affine conjugacy class of $\Gamma(\phi)$ do not depend either on the chosen splitting $\phi$, or the chosen $\RCD(0,N)$-structure on $X$. This will lead us to introduce the splitting degree $k(X)$ and the crystallographic class $\Gamma(X)$ of $X$, which are topological invariants of $X$. Finally, we will introduce the Albanese class $A(X)$ of $X$, which consists of orbifolds whose fundamental group belong to $\Gamma(X)$.\\

First of all, let us introduce the definition of splittings.

\begin{Definition}Let $(X,\dist,\m)$ be an $\RCD(0,N)$-structure on $X$, and denote $(\tilde{X},\tilde{\dist},\tilde{\m})$ its lift. A \emph{splitting of $(\tilde{X},\tilde{\dist},\tilde{\m})$} is an isomorphism $\phi\colon (\tilde{X},\tilde{\dist},\tilde{\m})\rightarrow (\overline{X},\overline{\dist},\overline{\m})\times\R^k$, where $\R^k$ is endowed with the Euclidean distance and Lebesgue measure, $k\in\N\cap[0,N]$ is called the \emph{degree of $\phi$}, and $(\overline{X},\overline{\dist},\overline{\m})$ is a compact $\RCD(0,N-k)$-space with trivial revised fundamental group called the \emph{soul of $\phi$}.
\end{Definition}

Thanks to Theorem 1.3 in \cite{Mondino-Wei_19}, we have the following existence result.

\begin{Theorem}\label{4.1}For every $\RCD(0,N)$-structure $(X,\dist,\m)$ on $X$, $(\tilde{X},\tilde{\dist},\tilde{\m})$ admits a splitting. Moreover, for every splitting $\phi$ of $(\tilde{X},\tilde{\dist},\tilde{\m})$, the group of isomorphism of $(\overline{X},\overline{\dist},\overline{\m})\times\R^k$ splits, i.e., we have:
$$
\Iso_{\mathrm{m.m.s.}}((\overline{X},\overline{\dist},\overline{\m})\times\R^k)=\Iso_{\mathrm{m.m.s.}}(\overline{X},\overline{\dist},\overline{\m})\times\Iso(\R^k),
$$
where $(\overline{X},\overline{\dist},\overline{\m})$ is the soul of $\phi$ and $k$ is the degree of $\phi$.
\end{Theorem}

Thanks to Theorem \ref{4.1}, Theorem \ref{3.3}, and Proposition \ref{prop3.4}, we can introduce the following notations.

\begin{Notation}Let $(X,\dist,\m)$ be an $\RCD(0,N)$-structure on $X$ and let $\phi$ be a splitting of $(\tilde{X},\tilde{\dist},\tilde{\m})$ with degree $k$ and soul $(\overline{X},\overline{\dist},\overline{\m})$. We write:
\begin{itemize}
\item[(i)] ${p_{S}^{\phi}}_*$ (resp. ${p_{\R}^{\phi}}_*$) the projection of $\Iso_{\mathrm{m.m.s.}}((\overline{X},\overline{\dist},\overline{\m})\times\R^k)$ onto $\Iso_{\mathrm{m.m.s.}}(\overline{X},\overline{\dist},\overline{\m})$ (resp. $\Iso(\R^k)$),
\item[(ii)] $\iota$ the inclusion of $\overline{\pi}_1(X)$ into $\Iso_{\mathrm{m.m.s.}}(\tilde{X},\tilde{\dist},\tilde{\m})$,
\item[(iii)] $\phi_*$ the isomorphism from $\Iso_{\mathrm{m.m.s.}}(\tilde{X},\tilde{\dist},\tilde{\m})$ onto $\Iso_{\mathrm{m.m.s.}}((\overline{X},\overline{\dist},\overline{\m})\times\R^k)$ defined by $\phi_*(\eta)\coloneqq\phi\circ\eta\circ\phi^{-1}$ for every $\eta\in \Iso_{\mathrm{m.m.s.}}(\tilde{X},\tilde{\dist},\tilde{\m})$.
\end{itemize}
We call $\rho_{S}^{\phi}\coloneqq{p_S^{\phi}}_*\circ\phi_*\circ\iota$ (resp. $\rho_{\R}^{\phi}\coloneqq {p_{\R}^{\phi}}_*\circ\phi_*\circ\iota$) the \emph{soul homomorphism associated to $\phi$} (resp. the \emph{Euclidean homomorphism associated to $\phi$}) and we write $K(\phi)\coloneqq\Ker(\rho_{\R}^{\phi})$ and $\Gamma(\phi)\coloneqq\Image(\rho_{\R}^{\phi})$.
\end{Notation}

The next result shows that the kernel and the image of the Euclidean homomorphism associated to a splitting enjoy particular group structures.

\begin{Proposition}\label{prop4.1}Let $(X,\dist,\m)$ be an $\RCD(0,N)$-structure on $X$ and let $\phi$ be a splitting of $(\tilde{X},\tilde{\dist},\tilde{\m})$ with degree $k$ and soul $(\overline{X},\overline{\dist},\overline{\m})$. Then, $K(\phi)$ is a finite normal subgroup of $\overline{\pi}_1(X)$ and $\Gamma(\phi)$ is a crystallographic subgroup of $\Iso(\R^k)$ (i.e. it acts cocompactly and discretly on $\R^k$).
\end{Proposition}

\begin{proof}First, let us show that $K(\phi)$ is finite. Observe that every element $\eta\in K(\phi)$ satisfies $\eta(\overline{X}\times\{0\})\cap(\overline{X}\times\{0\})\neq\varnothing$.
However, $\overline{X}\times\{0\}$ is a compact subset of $\overline{X}\times\R^k$ and $\overline{\pi}_1(X)$ acts properly on $\overline{X}\times\R^k$; thus, $K(\phi)$ is finite.\\
Now, let us show that $\Gamma(\phi)$ acts cocompactly on $\R^k$. Thanks to the first isomorphism theorem for topological spaces, there is a continuous map $\mu$ such that the following diagram is commutative:
$$
\begin{tikzcd} \overline{X}\times\R^k \arrow[r,"p_{\R^k}"] \arrow[d,"q_1"] & \R^k \arrow[d,"q_2"]\\
(\overline{X}\times\R^k)/\overline{\pi}_1(X) \arrow[r, "\mu"]
& \R^k/\Gamma(\phi)
\end{tikzcd}
$$
where $q_i$ ($i\in\{1,2\}$) are the quotient maps. Moreover, $\mu$ is surjective since $q_2\circ p_{\R^k}$ is surjective. Finally, $X$ is homeomorphic to $(\overline{X}\times\R^k)/\overline{\pi}_1(X)$; in particular, $(\overline{X}\times\R^k)/\overline{\pi}_1(X)$ is compact and $\R^k/\Gamma(\phi)$ is compact, being the image of a compact topological space by a continuous surjective map. In conclusion, $\Gamma(\phi)$ acts cocompactly on $\R^k$.\\
Let us prove that $\Gamma(\phi)$ acts discretely on $\R^k$ (i.e. its orbits are discrete subsets of $\R^k$). First, observe that it is sufficient to prove that $\Gamma(\phi)$ acts properly on $\R^k$. To prove this, let $K$ be a compact subset of $\R^k$ and let us show that there are only finitely many elements $g\in\Gamma(\phi)$ such that $g K\cap K\neq\varnothing$. By definition of $\Gamma(\phi)$, we have:
$$
\{g\in\Gamma(\phi), g K\cap K\neq\varnothing\}=\rho^{\phi}_{\R}(\{\eta\in\overline{\pi}_1(X), \eta(\overline{X}\times K)\cap(\overline{X}\times K)\neq\varnothing\}).
$$
However, $\overline{\pi}_1(X)$ acts properly on $\overline{X}\times{\R^k}$ and $\overline{X}\times K$ is compact; hence:
$$
\{\eta\in\overline{\pi}_1(X), \eta(\overline{X}\times K)\cap(\overline{X}\times K)\neq\varnothing\}
$$
is finite. Thus, $\{g\in\Gamma(\phi), g K\cap K\neq\varnothing\}$ is finite, being the image of a finite set.
\end{proof}

The following corollary of Proposition \ref{prop4.1} defines the splitting degree of $X$ (cf. \cite[Proposition 2.25]{MMP}).

\begin{cor}[Splitting degree $k(X)$]\label{cor5.1}The revised fundamental group $\overline{\pi}_1(X)$ is a finitely generated group which has polynomial growth of order $k(X)\in\N\cap[0,N]$. Moreover, given any $\RCD(0,N)$-structure $(X,\dist,\m)$ on $X$ with lift $(\tilde{X},\tilde{\dist},\tilde{\m})$, the degree of any splitting $\phi$ of $(\tilde{X},\tilde{\dist},\tilde{\m})$ is equal to $k(X)$. We call $k(X)$ the \emph{splitting degree of $X$}.
\end{cor}

\begin{proof}Thanks to Proposition \ref{prop4.1}, $\Gamma(\phi)$ is a crystallographic subgroup of $\Iso(\R^k)$, where $k\in[0,N]\cap\N$ is the degree of $\phi$. We need to prove that $\overline{\pi}_1(X)$ has polynomial growth of order $k$.\\
By Bieberbach's first Theorem (see Theorem 3.1 in \cite{Charlap_86}), $\Gamma(\phi)$ admits a normal subgroup $\Gamma(\phi)\cap\R^k$ such that $\Gamma(\phi)\cap\R^k$ is isomorphic to $\Z^k$ and $\Gamma(\phi)\cap\R^k$ has finite index in $\Gamma(\phi)$. In particular, $\Gamma(\phi)\cap\R^k$ is finitely generated, has polynomial growth of order $k$, and is a normal subgroup of $\Gamma(\phi)$ with finite index; thus, $\Gamma(\phi)$ is also finitely generated and has polynomial growth of order $k$. Now, $\overline{\pi}_1(X)/K(\phi)$ is isomorphic to $\Gamma(\phi)$; hence it is finitely generated with polynomial growth of order $k$. However, $K(\phi)$ is finite and is a normal subgroup of $\overline{\pi}_1(X)$; thus, $\overline{\pi}_1(X)$ is also finitely generated and has polynomial growth of order $k$.
\end{proof}

The revised fundamental group satisfies the following additional group property (which will be crucial in the proof of Theorem \ref{thmA}).

\begin{Proposition}\label{prop4.2}The revised fundamental group $\overline{\pi}_1(X)$ is a Hopfian group, i.e., every surjective group homomorphism from $\overline{\pi}_1(X)$ onto itself is an isomorphism.
\end{Proposition}

\begin{proof}First of all, let us recall some results on Group theory:
\begin{itemize}
\item[(i)] Noetherian groups (every subgroup is finitely generated) are Hopfian groups.
\item[(ii)] If $H$ is a normal subgroup of $G$ such that both $H$ and $G/H$ are Noetherian, then $G$ is Noetherian.
\item[(iii)] Finite groups are Noetherian.
\item[(iv)] Finitely generated abelian groups are Noetherian.
\end{itemize}
Let us fix an $\RCD(0,N)$-structure $(X,\dist,\m)$ on $X$ and let $\phi$ be a splitting of its lift $(\tilde{X},\tilde{\dist},\tilde{\m})$. By Proposition \ref{prop4.1} and Corollary \ref{cor5.1}, $\Gamma(\phi)$ is a crystallographic subgroup of $\Iso(\R^{k(X)})$. Hence, by Bieberbach's 1st Theorem (Theorem 3.1 in \cite{Charlap_86}), $\Gamma(\phi)\cap\R^{k(X)}$ is isomorphic to $\Z^{k(X)}$. In particular, $\Gamma(\phi)\cap\R^{k(X)}$ is Noetherian thanks to (iv). Moreover, $\Gamma(\phi)\cap\R^k$ is normal in $\Gamma(\phi)$ and the quotient $\Gamma(\phi)/\Gamma(\phi)\cap\R^k$ is finite. Thanks to (iii), $\Gamma(\phi)/\Gamma(\phi)\cap\R^k$ is Noetherian, and, using (ii), $\Gamma(\phi)$ is Noetherian. In addition, $K(\phi)$ is finite by Proposition \ref{prop4.1}, so it is Noetherian by (iii). Finally, $\overline{\pi}_1(X)/K(\phi)$ is isomorphic to $\Gamma(\phi)$ so it is Noetherian. In conclusion, thanks to (ii), $\overline{\pi}_1(X)$ is Noetherian; hence, it is Hopfian using (i).
\end{proof}

Given $k\in\N$, two crystallographic subgroups of $\Iso(\R^k)$ are called equivalent if they are conjugated by an affine transformation. The set $\mathrm{Crys}({k})$ of equivalence classes of crystallographic subgroups of $\Iso(\R^k)$ is a finite set thanks to Bieberbach's third Theorem (see Theorem 7.1 in \cite{Charlap_86}). The following result defines the crystallographic class of $X$.

\begin{Proposition}[Crystallographic class $\Gamma(X)$]\label{prop5.1}For $i\in\{1,2\}$, let $(X,\dist_i,\m_i)$ be an $\RCD(0,N)$-structure on $X$, and let $\phi_i$ be a splitting of $(\tilde{X},\tilde{\dist}_i,\tilde{\m}_i)$. Then $\Gamma(\phi_1)$ and $\Gamma(\phi_2)$ are equivalent as crystallographic subgroups of $\Iso(\R^{k(X)})$. We denote by $\Gamma(X)$ the common equivalence class and call it the \emph{crystallographic class of $X$}.
\end{Proposition}

\begin{proof}By Bieberbach's second Theorem (see Theorem 4.1 of \cite{Charlap_86}), two crystallographic subgroups of $\Iso(\R^{k})$ are conjugated by an affine transformation if and only if they are isomorphic (we let $k\coloneqq k(X)$). We need to show that $\Gamma(\phi_1)$ and $\Gamma(\phi_2)$ are isomorphic. Observe that, for $i\in\{1,2\}$, we have the following exact sequence of groups:
$$
\begin{tikzcd} \{1\} \arrow[r] & K(\phi_i) \arrow[r,"\iota"] & \overline{\pi}_1(X)\arrow[r,"\rho^{\phi_i}_{\R}"]& \Gamma(\phi_i)\arrow[r]&\{1\}\end{tikzcd},
$$
where $\iota$ is just the inclusion, $\Gamma(\phi_i)$ is a crystallographic subgroup of $\Iso(\R^k)$, and $K(\phi_i)$ is finite. By Remark 2.5 of \cite{Wilking_00}, $K(\phi_i)=\iota(K(\phi_i))$ is uniquely characterized as the maximal finite normal subgroup of $\overline{\pi}_1(X)$. In particular, we necessarily have $K(\phi_1)=K(\phi_2)$. In conclusion, $\Gamma(\phi_1)\simeq\overline{\pi}_1(X)/K(\phi_1)=\overline{\pi}_1(X)/K(\phi_2)\simeq\Gamma(\phi_2)$; thus, $\Gamma(\phi_1)$ is isomorphic to $\Gamma(\phi_2)$.
\end{proof}

Given $k\in\N$, and $\Gamma$ a crystallographic subgroup of $\Iso(\R^{k})$, the quotient space $\R^{k}/\Gamma$ has the structure of a compact flat orbifold of dimension $k$, whose orbifold metric $\dist_{\Gamma}$ satisfies:
\begin{equation}\label{orbdist}
\dist_{\Gamma}([x],[y])=\inf\{\lvert x'-y'\rvert\},
\end{equation}
where $x,y\in\R^k$, $[x]$ and $[y]$ are their equivalence class in $\R^k/\Gamma$, and the infimum is taken over all $x'\in[x]$ and $y'\in[y]$. Moreover, equivalent crystallographic groups give rise to orbifolds that are affinely equivalent.\\
Conversely, given a compact flat orbifold $(X,\dist)$ of dimension $k$, the orbifold fundamental group $\pi_1^{\mathrm{orb}}(X)$ acts by isometries on the orbifold universal cover, which is $\R^k$ (the action being discrete and cocompact). Hence, one can associate a crystallographic group to $(X,\dist)$. Finally, two affinely equivalent flat orbifolds of dimension $k$ have isomorphic orbifold fundamental groups. Hence, by Bieberbach's second Theorem (see Theorem 4.1 in \cite{Charlap_86}), they give rise to equivalent crystallographic groups (see the introduction of Section 2.1 in \cite{Bettiol-Derdzinski-Piccione_18} for more details and some references).\\
Therefore, there is a one-to-one correspondence between equivalence classes of crystallographic subgroups of $\Iso(\R^k)$ and affine equivalence classes of compact flat orbifolds of dimension $k$. This leads us to the Definition of the Albanese class of $X$.

\begin{Definition}[Albanese class $A(X)$]\label{def7.2}We write $A(X)$ the affine equivalence class of compact flat orbifold determined by $\Gamma(X)$, and call it the \emph{Albanese class of $X$}. More explicitly, $A(X)$ is the set of all flat orbifolds $(\R^{k(X)}/\Gamma,\dist_{\Gamma})$, where $\Gamma\in\Gamma(X)$, and $\dist_{\Gamma}$ is defined in equation \ref{orbdist}.
\end{Definition}

\subsection{Moduli spaces and their topology}\label{section6}

In Section \ref{subsection6.2}, we will introduce the moduli space of pointed $\RCD(0,N)$-structures on $X$. Then, in Section \ref{subsection6.3}, we will introduce the moduli space of equivariant pointed $\RCD(0,N)$-structures on $\tilde{X}$. In particular, (based on the equivariant distance introduced by Fukaya and Yamaguchi in \cite{Fukaya-Yamaguchi-92}), we will introduce the equivariant pmGH topology.

\subsubsection{Moduli space of pointed {RCD}({0},{N})-structures}\label{subsection6.2}

Throughout the paper, we will use the following definition of pointed metric measure spaces.

\begin{Definition}A pointed metric measure space (p.m.m.s. for short) is a $4$-tuple $(X,\dist,\m,\ast)$, where $(X,\dist,\m)$ is a m.m.s. and $\ast\in X$.
\end{Definition}

As for m.m.s., there are two distinct notions of isomorphism between two pointed metric measure spaces. In this paper, we decided to use the following definition (which emphasizes on the whole space's metric structure, not only on the metric structure of the measure's support).

\begin{Definition}\label{def:pmms}Two p.m.m.s. $(X_1,\dist_1,\m_1,\ast_1)$ and $(X_2,\dist_2,\m_2,\ast_2)$ are \emph{isomorphic} when there is an isomorphism of metric measure spaces $\phi\colon(X_1,\dist_1,\m_1)\rightarrow(X_2,\dist_2,\m_2)$ such that $\phi(\ast_1)=*_2$.
\end{Definition}

Thanks to Theorem 2.7 of \cite{Abraham-Delma-Hoscheit_13}, and Remark 3.29 of \cite{Gigli-Mondino-Savare_15}, we have the following result.

\begin{Theorem}\label{th: pmGH topology}The Gromov--Hausdorff--Prokhorov distance $\dist_{\mathrm{GHP}}$ (see Section 2.3 of \cite{Abraham-Delma-Hoscheit_13}) is a complete separable metric on the set $\kX^{\mathrm{p}}$ of isomorphism classes of pointed metric measure spaces that are locally compact and geodesic. Moreover, $\dist_{\mathrm{GHP}}$ metrizes the pointed measured Gromov--Hausdorff topology (introduced in Definition 27.30 of \cite{Villani_09}).
\end{Theorem}

As we will see in Remark \ref{rem: approximtation} below, it is possible to realize the pmGH convergence using maps with small distortion.

\begin{Notation}\label{not: distorsion}Let $f\colon(X,\dist_X)\to(Y,\dist_Y)$ be a map between metric spaces, we denote:
$$
\Dis(f)\coloneqq\sup\{\lvert \dist_Y(f(x),f(y))-\dist_X(x,y)\rvert\},
$$
where the supremum is taken over all couples $x,y\in X$.
\end{Notation}

\begin{Remark}\label{rem: approximtation}Assume that $\{\mathcal{X}_n=(X_n,\dist_n,\m_n,*_n)\}_{n\in\N\cup\{\infty\}}$ is a familly of locally compact geodesic p.m.m.s. such that $\mathcal{X}_n\to\mathcal{X}_{\infty}$ in the pmGH topology. Then, Theorem \ref{th: pmGH topology} implies that there exists a sequence $\{f_n,g_n,\epsilon_n\}$ where $f_n\colon (X_n,\ast_n)\to(X_{\infty},\ast_{\infty})$ and $g_n\colon(X_{\infty},\ast_{\infty})\to(X_n,\ast_n)$ are pointed Borel maps, $\epsilon_n\to0$, and the following properties are satisfied:
\begin{itemize}
\item[(i)]for every $x\in X_n$, $\dist_n(\ast_n,x)\leq\epsilon_n^{-1}\implies\dist_n(g_n(f_n(x)),x)\leq\epsilon_n$, and, for every $x\in X_{\infty}$, $\dist_{\infty}(\ast_{\infty},x)\leq\epsilon_n^{-1}\implies\dist_{\infty}(f_n(g_n(x)),x)\leq\epsilon_n$,
\item[(ii)]$\max\{\Dis({f_n}_{\lvert B_n(\epsilon_n^{-1})}),\Dis({g_n}_{\lvert B_{\infty}(\epsilon_n^{-1})})\}\leq\epsilon_n$ (see Notation \ref{not: distorsion}),
\item[(iii)]$\max\{\dist_{\mathcal{P}}^{\epsilon_n^{-1}}({f_n}*{\m}_n,{\m}_{\infty}),\dist_{\mathcal{P}}^{\epsilon_n^{-1}}({g_n}*{\m}_{\infty},{\m}_n)\}\leq\epsilon_n$.
\end{itemize}
Such a sequence is said to realize the convergence of $\{\mathcal{X}_n\}$ to $\mathcal{X}_{\infty}$ in the pmGH topology.\end{Remark}

We conclude this section by introducing the moduli space of pointed $\RCD(0,N)$-structures on $X$.

\begin{Notation}We introduce the following spaces:
\begin{itemize}
\item[(i)] $\mathfrak{RCD}^{\mathrm{p}}(0,N)\subset\kX$ is the set of isomorphism classes of pointed $\RCD(0,N)$ spaces with full support, endowed with the pmGH-topology (seen as a subspace of $\kX^{\mathrm{p}}$),
\item[(ii)] $\mathfrak{R}^{\mathrm{p}}_{0,N}(X)$ is the set of all pointed $\RCD(0,N)$-structures on $X$,
\item[(iii)] $\M^{\mathrm{p}}_{0,N}(X)$ is the quotient of $\mathfrak{R}^{\mathrm{p}}_{0,N}(X)$ by isomorphisms, endowed with the pmGH-topology (seen as a subspace of $\mathfrak{RCD}^{\mathrm{p}}(0,N)$).
\end{itemize}
 We call $\M_{0,N}^{\mathrm{p}}(X)$ the \emph{moduli space of pointed $\RCD(0,N)$-structures on $X$}.
\end{Notation}

\subsubsection{Moduli space of equivariant pointed {RCD}({0},{N})-structures}\label{subsection6.3}

First of all, we introduce equivariant pointed $\RCD(0,N)$-structures on $\tilde{X}$. Here, in comparison with the definition of equivariant metric triple given by Fukaya and Yamaguchi in \cite{Fukaya-Yamaguchi-92}, both the topological space and the group action are fixed.

\begin{Definition}\label{def:equiv}A pointed $\RCD(0,N)$-structure $(\tilde{X},\tilde{\dist},\tilde{\m},\tilde{\ast})$ on $\tilde{X}$ is called equivariant if $\overline{\pi}_1(X)$ acts by isomorphisms on $(\tilde{X},\tilde{\dist},\tilde{\m})$.
\end{Definition}

The following definition introduces equivariant isomorphisms between equivariant pointed $\RCD(0,N)$-structures on $\tilde{X}$.

\begin{Definition}For $i\in\{1,2\}$, let $\tilde{\X}_i=(\tilde{X},\tilde{\dist}_i,\tilde{\m}_i,\tilde{\ast}_i)$ be an equivariant pointed $\RCD(0,N)$-structure on $\tilde{X}$. We say that $\tilde{\X}_1$ and $\tilde{\X}_2$ are equivariantly isomorphic when there is an isomorphism $\phi$ of $\overline{\pi}_1(X)$ and an isomorphism $f\colon\tilde{\X}_1\to\tilde{\X}_2$ of p.m.m.s. such that $f(\gamma x)=\phi(\gamma)f(x)$, for every $\gamma\in\overline{\pi}_1(X)$, and every $x\in\tilde{X}$.
\end{Definition}

We now introduce the space (and moduli space) of equivariant pointed $\RCD(0,N)$-structures on $\tilde{X}$.

\begin{Notation}We introduce the following spaces:
\begin{itemize}
\item[(i)] $\kR_{0,N}^{\mathrm{p,eq}}(\tilde{X})$ the set of equivariant pointed $\RCD(0,N)$-structures on $\tilde{X}$,
\item[(ii)] $\M_{0,N}^{\mathrm{p,eq}}(\tilde{X})$ the quotient space of $\kR_{0,N}^{\mathrm{p,eq}}(\tilde{X})$ by equivariant ismormophisms.
\end{itemize}
We call $\M_{0,N}^{\mathrm{p,eq}}(\tilde{X})$ the \emph{moduli space of equivariant pointed $\RCD(0,N)$-structures on $\tilde{X}$}.
\end{Notation}

To define a topological structure on $\M_{0,N}^{\mathrm{p,eq}}(\tilde{X})$ , we start by introducing the equivariant pointed distance on $\kR_{0,N}^{\mathrm{p,eq}}(\tilde{X})$.

\begin{Definition}\label{def:equivtopo}Let $\epsilon>0$, and, for $i\in\{1,2\}$, let $\tilde{\X}_i=(\tilde{X},\tilde{\dist}_i,\tilde{\m}_i,\tilde{\ast}_i)$ be an equivariant pointed $\RCD(0,N)$-structure on $\tilde{X}$. An \emph{equivariant pointed $\epsilon$-isometry} between $\tilde{\X}_1$ and $\tilde{\X}_2$ is a triple $(f,g,\phi)$ where $f\colon\tilde{X}\to\tilde{X}$ and $g\colon\tilde{X}\to\tilde{X}$ are Borel maps and $\phi$ is an isomorphism of $\overline{\pi}_1(X)$ such that:
\begin{itemize}
\item[(i)]$f(\tilde{\ast}_1)=\tilde{\ast}_2$ and $g(\tilde{\ast}_2)=\tilde{\ast}_1$,
\item[(ii)]for every $\gamma\in\overline{\pi}_1(X)$ and $x\in\tilde{X}$, $f(\gamma x)=\phi(\gamma)f(x)$ and $g(\gamma x)=\phi^{-1}(\gamma)g(x)$,
\item[(iii)]for every $x,y\in\tilde{X}$, $\tilde{\dist}_1(x,y)\leq\epsilon^{-1}\implies \lvert \tilde{\dist}_2(f(x),f(y))-\tilde{\dist}_1(x,y)\rvert\leq\epsilon$, and $\tilde{\dist}_2(x,y)\leq\epsilon^{-1}\implies \lvert \tilde{\dist}_1(g(x),g(y))-\tilde{\dist}_2(x,y)\rvert\leq\epsilon$,
\item[(iv)]for every $x\in\tilde{X}$, $\tilde{\dist}_1(g\circ f(x),x)\leq \epsilon$ and $\tilde{\dist}_2(f\circ g (x),x)\leq \epsilon$,
\item[(v)]$\max\{\dist_{\mathcal{P}}^{\epsilon^{-1}}(f_*\tilde{\m}_1,\tilde{\m}_2),\dist_{\mathcal{P}}^{\epsilon^{-1}}(g_*\tilde{\m}_2,\tilde{\m}_1)\}\leq\epsilon$.
\end{itemize}
We define $\D_{\mathrm{p}}^{\mathrm{eq}}(\tilde{\X}_1,\tilde{\X}_2)$ the \emph{equivariant pointed distance between $\tilde{\X}_1$ and $\tilde{\X}_2$} as the minimum between $1/24$ and the infimum of all $\epsilon>0$ such that there exists an equivariant pointed $\epsilon$-isometry between $\tilde{\X}_1$ and $\tilde{\X}_2$.
\end{Definition}

The following result shows that we can endow $\M_{0,N}^{\mathrm{p,eq}}(\tilde{X})$ with a metrizable topology.

\begin{Proposition}\label{prop:equivtopo}$\Dpeq$ induces a metrizable uniform structure on $\M_{0,N}^{\mathrm{p,eq}}(\tilde{X})$.
\end{Proposition}

\begin{proof}See \hyperref[Appendix]{Appendix}.
\end{proof}

From now on, we endow $\M_{0,N}^{\mathrm{p,eq}}(\tilde{X})$ with the topology induced by $\Dpeq$, which we call the \emph{equivariant pmGH-topology}.

\subsection{Maps between moduli spaces}\label{sec7}

In Section \ref{section7.1}, we are going to introduce the lift and push-forward maps. As we will explain at the end of that section, a consequence of Theorem \ref{thmA} is that these maps are homeomorphisms and respectively inverse to each other (see Corollary \ref{corA}). Then, in Section \ref{sec7.2}, we will introduce the Albanese map and the soul map associated to $X$.

\subsubsection{Lift and push-forward maps}\label{section7.1}

Thanks to Corollary \ref{corlift}, we can define the lift of a pointed $\RCD(0,N)$-structure.

\begin{Definition}\label{def9.1}Let $(X,\dist,\m,x)$ be a pointed $\RCD(0,N)$-structure on $X$ and let $\tilde{x}\in p^{-1}(x)$. We define $p_{\tilde{x}}^*(X,\dist,\m,x)\coloneqq(\tilde{X},\tilde{\dist},\tilde{\m},\tilde{x})$, where $(\tilde{X},\tilde{\dist},\tilde{\m})$ is the lift of $(X,\dist,\m)$.
\end{Definition}

\begin{Remark}\label{rem:p*}For $i\in\{1,2\}$, let $(X,\dist_i,\m_i,x_i)$ be a pointed $\RCD(0,N)$-structure on $X$, and let $\tilde{x}_i\in p^{-1}(x_i)$. If $(X,\dist_1,\m_1,x_1)$ and $(X,\dist_2,\m_2,x_2)$ are isomorphic, then $p_{\tilde{x}_1}^*(X_1,\dist_1,\m_1,x_1)$ is equivariantly isomorphic to $p_{\tilde{x}_1}^*(X_2,\dist_2,\m_2,x_2)$.
\end{Remark}

Thanks to Remark \ref{rem:p*}, we can define the lift map associated to $X$.

\begin{Definition}[Lift map] The \emph{lift map associated to $X$} is the unique map $p^*\colon\M^{\mathrm{p}}_{0,N}(X)\to \M_{0,N}^{\mathrm{p,eq}}(\tilde{X})$ that satisfies $p^*[X,\dist,\m,x]=[p_{\tilde{x}}^*(X,\dist,\m,x)]$ for every $(X,\dist,\m,x)\in\mathfrak{R}_{0,N}^{\mathrm{p}}(X)$ and $\tilde{x}\in p^{-1}(x)$.
\end{Definition}

Thanks to Proposition \ref{propX}, we can define the push-forward of an equivariant pointed $\RCD(0,N)$-structure.

\begin{Definition}\label{def:push}Let $(\tilde{X},\tilde{\dist},\tilde{\m},\tilde{x})$ be an equivariant pointed $\RCD(0,N)$-structures on $\tilde{X}$. We define $p_*(\tilde{X},\tilde{\dist},\tilde{\m},\tilde{x})$ as the unique pointed $\RCD(0,N)$-structure on $X$ such that $p\colon (\tilde{X},\tilde{\dist},\tilde{\m},\tilde{x})\to p_*(\tilde{X},\tilde{\dist},\tilde{\m},\tilde{x})$ is a pointed local isomorphism.
\end{Definition}

\begin{Remark}\label{rem:p_*}For $i\in\{1,2\}$, let $(\tilde{X},\tilde{\dist}_i,\tilde{\m}_i,\tilde{x}_i)$ be an equivariant pointed $\RCD(0,N)$-structure on $\tilde{X}$. If $(\tilde{X},\tilde{\dist}_1,\tilde{\m}_1,\tilde{x}_1)\sim(\tilde{X},\tilde{\dist}_2,\tilde{\m}_2,\tilde{x}_2)$, then $p_*(\tilde{X},\tilde{\dist}_1,\tilde{\m}_1,\tilde{x}_1)$ is isomorphic to $p_*(\tilde{X},\tilde{\dist}_2,\tilde{\m}_2,\tilde{x}_2)$.
\end{Remark}

Thanks to Remark \ref{rem:p_*}, we can define the push-forward map associated to $X$.

\begin{Definition}[Push-forward map] The \emph{push-forward map associated to $X$} is the unique map \[p_*\colon \M_{0,N}^{\mathrm{p,eq}}(\tilde{X})\to \M^{\mathrm{p}}_{0,N}(X)\] satisfying $p_*[\tilde{X},\tilde{\dist},\tilde{\m},\tilde{x}]=[p_*(\tilde{X},\tilde{\dist},\tilde{\m},\tilde{x})]$ for every $(\tilde{X},\tilde{\dist},\tilde{\m},\tilde{x})\in\kR_{0,N}^{\mathrm{p,eq}}(\tilde{X})$.
\end{Definition}

Thanks to Remark \ref{rem:pushlift}, we have the following proposition.

\begin{Proposition}\label{propliftpush}The lift map $p^*\colon\M^{\mathrm{p}}_{0,N}(X)\to \M_{0,N}^{\mathrm{p,eq}}(\tilde{X})$ and the push-forward map $p_*\colon \M_{0,N}^{\mathrm{p,eq}}(\tilde{X})\to \M^{\mathrm{p}}_{0,N}(X)$ are respectively inverse to each other.
\end{Proposition}

Observe that Corollary \ref{corA} immediately follows from Proposition \ref{propliftpush} and Theorem \ref{thmA} (which we will prove in Section \ref{proofA}).

\subsubsection{Albanese and soul maps}\label{sec7.2}

First of all, we introduce the moduli space of flat metrics on the Albanese class $A(X)$ (introduced in Definition \ref{def7.2}). This moduli space will act as the codomain of the Albanese map.

\begin{Definition}[$\mathscr{M}_{\text{flat}}(A(X))$]The \emph{moduli space of flat metrics on $A(X)$} is the quotient of $A(X)$ by isometry equivalence, endowed with the Gromov--Hausdorff distance $\GH$ (see Definition 7.3.10 in \cite{Burago-Ivanov_01}).
\end{Definition}

The following remark will be helpful in the proof of Theorem \ref{thmC}. It is also interesting on its own as it gives a more explicit way to see the moduli space $\mathscr{M}_{\text{flat}}(A(X))$.

\begin{Remark}\label{rem5.1}Given any element $\Gamma\in\Gamma(X)$, the moduli space of flat metrics on $A(X)$ is isometric to the moduli space of flat metrics on the compact orbifold $\R^k/\Gamma$ (endowed with the Gromov--Hausdorff distance), which we denote $\mathscr{M}_{\mathrm{flat}}(\R^k/\Gamma)$ (see Section 4.2 in \cite{Bettiol-Derdzinski-Piccione_18} for more details on $\mathscr{M}_{\mathrm{flat}}(\R^k/\Gamma)$).
\end{Remark}



The next lemma is fundamental to introduce the Albanese and soul maps associated to $X$.

\begin{Lemma}\label{lem:albwelldef}For $i\in\{1,2\}$, let $(X,\dist_i,\m_i)$ be an $\RCD(0,N)$-structure on $X$, and let $\phi_i$ be a splitting of $(\tilde{X},\tilde{\dist}_i,\tilde{\m}_i)$ with soul $(\overline{X}_i,\overline{\dist}_i,\overline{\m}_i)$. If $(X,\dist_1,\m_1)$ and $(X,\dist_2,\m_2)$ are isomorphic, then $(\overline{X}_1,\overline{\dist}_1,\overline{\m}_1)$ is isomorphic to $(\overline{X}_2,\overline{\dist}_2,\overline{\m}_2)$, and $(\R^{k}/\Gamma(\phi_1),\dist_{\Gamma(\phi_1)})$ is isometric to $(\R^{k}/\Gamma(\phi_2),\dist_{\Gamma(\phi_2)})$.
\end{Lemma}

\begin{proof}Let us fix an isomorphism $\phi\colon(X,\dist_1,\m_1)\to(X,\dist_2,\m_2)$. We can lift $\phi$ to the universal covers to get an isomorphism $\tilde{\phi}\colon(\tilde{X},\tilde{\dist}_1,\tilde{\m}_1)\to(\tilde{X},\tilde{\dist}_2,\tilde{\m}_2)$ such that $p\circ\tilde{\phi}=\phi\circ p$. Now, let $\mu\coloneqq \phi_2\circ\tilde{\phi}\circ\phi_1^{-1}$. Since, $\overline{X}_1$ and $\overline{X}_2$ are compact, $\mu$ is of the form $\mu=(\mu_S,\mu_{\R})$, where $\mu_S\colon(\overline{X},\overline{\dist}_1,\overline{\m}_1)\to(\overline{X},\overline{\dist}_2,\overline{\m}_2)$ is an isomorphism, and $\mu_{\R}\in\Iso(\R^k)$ (where $k\coloneqq k(X)$). In particular $(\overline{X}_1,\overline{\dist}_1,\overline{\m}_1)$ is isomorphic to $(\overline{X}_2,\overline{\dist}_2,\overline{\m}_2)$.\\
We are going to show that $\Gamma(\phi_2)=\mu_{\R}\Gamma(\phi_1)\mu_{\R}^{-1}$. Let $\bar{x}_1\in\overline{X}_1$, let $t\in\R^k$, let $\alpha\in\overline{\pi}_1(X)$ and define $\tilde{z}\coloneqq\phi_1^{-1}(\bar{x}_1,t)$. By definition of the soul and Euclidian homomorphisms associated to $\phi_1$ and $\phi_2$, we have:
$$
\mu(\rho^{\phi_1}_S(\alpha)\cdot\bar{x}_1,\rho^{\phi_1}_{\R}(\alpha)\cdot t)=(\rho^{\phi_2}_S(\eta)\cdot\mu_S(\bar{x}_1),\rho^{\phi_2}_{\R}(\eta)\cdot\mu_{\R}(t))=(\mu_S(\rho^{\phi_1}_S(\alpha)\cdot\bar{x}_1),\mu_{\R}(\rho^{\phi_1}_{\R}(\alpha)\cdot t)),
$$
where $\eta\coloneqq\tilde{\phi}_*(\alpha)$, and $\tilde{\phi}_*$ is the automorphism of $\overline{\pi}_1(X)$ defined by $\tilde{\phi}_*(\alpha)\coloneqq\tilde{\phi}\circ\alpha\circ\tilde{\phi}^{-1}$. In particular, for every $t\in\R^{k}$ and $\alpha\in\overline{\pi}_1(X)$, we have $\mu_{\R}(\rho^{\phi_1}_{\R}(\alpha)\cdot t)=\rho^{\phi_2}_{\R}\circ\tilde{\phi}_*(\alpha)\cdot\mu_{\R}(t)$. Thus, for every $\alpha\in\overline{\pi}_1(X)$, we have $\mu_{\R}\circ\rho^{\phi_1}_{\R}(\alpha)\circ\mu_{\R}^{-1}=\rho^{\phi_2}_{\R}\circ\tilde{\phi}_*(\alpha)$. In particular, by definition of $\Gamma(\phi_1)$ and $\Gamma(\phi_2)$, and since $\tilde{\phi}_*(\overline{\pi}_1(X))=\overline{\pi}_1(X)$, we have $\Gamma(\phi_2)=\mu_{\R}\Gamma(\phi_1)\mu_{\R}^{-1}$. In conclusion, using Lemma 4.1 in \cite{Bettiol-Derdzinski-Piccione_18}, $(\R^k/\Gamma(\phi_1),\dist_{\Gamma(\phi_1)})$ is isometric to $(\R^k/\Gamma(\phi_2),\dist_{\Gamma(\phi_2)})$, which concludes the proof.
\end{proof}

Thanks to Lemma \ref{lem:albwelldef}, we can define the Albanese and soul maps.

\begin{Definition}[Albanese and soul maps]Given an $\RCD(0,N)$-structure $(X,\dist,\m)$ on $X$, and given a splitting $\phi$ of $(\tilde{X},\tilde{\dist},\tilde{\m})$ with soul $(\overline{X},\overline{\dist},\overline{\m})$, we define:
$$
\mathcal{A}([X,\dist,\m])\coloneqq[\R^k/\Gamma(\phi),\dist_{\Gamma(\phi)}]\in\mathscr{M}_{\text{flat}}(A(X)),
$$
and:
$$
\mathcal{S}([X,\dist,\m])\coloneqq[\overline{X},\overline{\dist},\overline{\m}]\in\mathfrak{RCD}(0,N-k(X)).
$$
The map $\mathcal{A} \colon \M_{0,N}(X)\to\mathscr{M}_{\text{flat}}(A(X))$ is called the \emph{Albanese map associated to $X$}, and the map 
$$\mathcal{S}\colon\M_{0,N}(X) \to\mathfrak{RCD}(0,N-k(X))$$
 is called the \emph{soul map associated to $X$}.
\end{Definition}

We end this section with the following surjectivity result.

\begin{Proposition}The Albanese map associated to $X$ is surjective from $\M_{0,N}(X)$ onto $\mathscr{M}_{\text{flat}}(A(X))$.
\end{Proposition}

\begin{proof}First of all, let $(X,\dist_0,\m_0)$ be a reference $\RCD(0,N)$-structure on $X$, and let $\phi_0$ be a splitting of its lift $(\tilde{X},\tilde{\dist}_0,\tilde{\m}_0)$ with soul $(\overline{X}_0,\overline{\dist}_0,\overline{\m}_0)$. Now, let $\Gamma\in\Gamma(X)$ and let us show that there is some $(X,\dist,\m)\in\kR_{0,N}(X)$ such that $\mathcal{A}([X,\dist,\m])=[\R^k/\Gamma,\dist_{\Gamma}]$.\\
Since $\Gamma(\phi_0)\in\Gamma(X)$, there is $\alpha\in\Aff(\R^k)$ such that $\Gamma=\alpha\Gamma(\phi_0)\alpha^{-1}$. Now, let $\psi\coloneqq (\id_{\overline{X}_0},\alpha)\circ\phi_0$, and consider the metric measure structure $(\tilde{\dist},\tilde{\m})$ defined as the pull back by $\psi$ of $(\overline{\dist}_0\times\dist_{\mathrm{eucli}},\overline{\m}_0\otimes\mathcal{L}_k)$. Note that $\psi$ is a homeomorphism, and $(\overline{X}_0\times{\R^k},\overline{\dist}_0\times\dist_{\mathrm{eucli}},\overline{\m}_0\otimes\mathcal{L}_k)$ is an $\RCD(0,N)$ space; hence, $(\tilde{X},\tilde{\dist},\tilde{\m})$ is an $\RCD(0,N)$-structure on $\tilde{X}$.\\
Now, we are going to show that $(\tilde{X},\tilde{\dist},\tilde{\m})$ is the lift of some $(X,\dist,\m)$. Thanks to Remark \ref{rem:pushlift}, it is equivalent to show that $\overline{\pi}_1(X)\subset\Iso_{\mathrm{m.m.s.}}(\tilde{X},\tilde{\dist},\tilde{\m})$, which is itself equivalent to $\psi_{*}(\overline{\pi}_1(X))\subset \Iso_{\mathrm{m.m.s.}}(\overline{X}_0\times{\R^k},\overline{\dist}_0\times\dist_{\mathrm{eucli}},\overline{\m}_0\otimes\mathcal{L}_k)$. Let $\eta\in\overline{\pi}_1(X)$, then $\psi_*(\eta)=\psi\circ\eta\circ\psi^{-1}=(\id_{\overline{X}_0},\alpha)\circ\phi_{0_*}(\eta)\circ(\id_{\overline{X}_0},\alpha^{-1})=(\rho_S^{\phi_0}(\eta),\alpha\circ\rho_{\R}^{\phi_0}(\eta)\circ\alpha^{-1})$.
Note that $\rho_S^{\phi_0}(\eta)\in\Iso_{\mathrm{m.m.s.}}(\overline{X}_0,\overline{\dist}_0,\overline{\m}_0)$ and $\alpha\circ\rho_{\R}^{\phi_0}(\eta)\circ\alpha^{-1}\in\alpha\Gamma(\phi_0)\alpha^{-1}=\Gamma\subset\Iso(\R^k)$; hence, $\psi_*(\eta)\in\Iso_{\mathrm{m.m.s.}}(\overline{X}_0\times{\R^k},\overline{\dist}_0\times\dist_{\mathrm{eucli}},\overline{\m}_0\otimes\mathcal{L}_k)$. In conclusion, there is an $\RCD(0,N)$-structure $(X,\dist,\m)\in\kR_{0,N}(X)$ whose lift is $(\tilde{X},\tilde{\dist},\tilde{\m})$. By construction, $\psi$ is a splitting of $(\tilde{X},\tilde{\dist},\tilde{\m})$ with soul $(\overline{X}_0,\overline{\dist}_0,\overline{\m}_0)$. Moreover, we have seen above that, for every $\eta\in\overline{\pi}_1(X)$, we have $\rho_{\R}^{\psi}(\eta)={p_{\R}^{\psi}}_*\circ\psi_*(\eta)=\alpha\circ\rho_{\R}^{\phi_0}(\eta)\circ\alpha^{-1}$. Hence, $\Gamma(\psi)=\alpha\Gamma(\phi_0)\alpha^{-1}=\Gamma$, and we get $\mathcal{A}([X,\dist,\m])=[\R^k/\Gamma,\dist_{\Gamma}]$.
\end{proof}

\section{Proof of the main results}\label{proof of main}

\subsection{Proof of Theorem \ref{thmA}}\label{proofA}

First of all, let us introduce the systole associated to an $\RCD(0,N)$-structure on $X$. Finding a uniform lower bound on the systoles associated to a sequence will be the key to prove Theorem \ref{thmA}.

\begin{Definition}[Systole of an $\RCD(0,N)$-structure]The systole associated to an $\RCD(0,N)$-structure $(X,\dist,\m)$ on $X$ is the quantity $\mathrm{sys}(X,\dist)\coloneqq\inf\{\tilde{\dist}(\eta\cdot\tilde{x},\tilde{x})\}$, where the infimum is taken over all point $\tilde{x}\in\tilde{X}$ and $\eta\in\overline{\pi}_1(X)\backslash\{\id\}$. Whenever $\overline{\pi}_1(X)$ is trivial, we define $\mathrm{sys}(X,\dist)\coloneqq\infty$.
\end{Definition}

The following proposition relates the systole of an $\RCD(0,N)$-structure $(X,\dist,\m)$ on $X$ and the quantity $\delta(X,\dist)$ introduced in Theorem \ref{3.3}.

\begin{Proposition}\label{prop6.1}Let $(X,\dist,\m)$ be an $\RCD(0,N)$-structure on $X$. Then, $\mathrm{sys}(X,\dist)=2\delta(X,\dist)$, where $\delta(X,\dist)$ is defined in Theorem \ref{3.3}.
\end{Proposition}

\begin{proof}Let $\delta<\delta(X,\dist)$, let $\eta\in\overline{\pi}_1(X)\backslash\{\id\}$, and let $\tilde{x}\in\tilde{X}$. Then, by Proposition \ref{prop3.4} and Theorem \ref{3.3}, $p$ induces a homeomorphism from $B_{\tilde{\dist}}(\tilde{x},\delta)$ (resp. $B_{\tilde{\dist}}(\eta\cdot\tilde{x},\delta)$) onto $B_{\dist}(x,\delta)$, where $x\coloneqq p(\tilde{x})$. Seeking for a contradiction, assume that there exists $\tilde{y}\in B_{\tilde{\dist}}(\tilde{x},\delta)\cap B_{\tilde{\dist}}(\eta\cdot\tilde{x},\delta)$. Then, $\dist(\eta\cdot\tilde{y},\eta\cdot\tilde{x})=\dist(\tilde{y},\tilde{x})<\delta$. In particular, $\tilde{y}$ and $\eta\cdot\tilde{y}$ are two distinct elements of $B_{\tilde{\dist}}(\eta\cdot\tilde{x},\delta)$ which have the same image under $p$, which is the contradiction we were looking for. Hence, $B_{\tilde{\dist}}(\tilde{x},\delta)\cap B_{\tilde{\dist}}(\eta\cdot\tilde{x},\delta)=\varnothing$. In particular, $\tilde{\dist}(\eta\cdot\tilde{x},\tilde{x})\geq2\delta$; thus, $2\delta\leq\mathrm{sys}(X,\dist)$. Since that holds for every $\delta<\delta(X,\dist)$, we have $2\delta(X,\dist)\leq\mathrm{sys}(X,\dist)$.\\
Now assume that $\delta(X,\dist)<\delta$. Then, there is some $x\in X$ such that $B_{\dist}(x,\delta)$ is not evenly covered by $p$. Therefore, given any $\tilde{x}\in p^{-1}(x)$, there exists $\tilde{y}_i\in B_{\tilde{\dist}}(\tilde{x},\delta)$ ($i\in\{1,2\}$) such that $p\tilde{y}_1=p\tilde{y}_2$, $\tilde{y}_1\neq \tilde{y}_2$. Hence, there exists $\gamma\in\overline{\pi}_1(X)\backslash\{\id\}$ such that $\tilde{y}_2=\gamma\tilde{y}_1$; thus, $\tilde{\dist}(\tilde{y}_1, \tilde{y}_2)=\tilde{\dist}(\gamma \tilde{y}_1, \tilde{y}_1)\leq 2\delta$. Therefore, we have $\mathrm{sys}(X,\dist)\leq 2\delta$. Thus, letting $\delta$ go to $\delta(X,\dist)$, we finally obtain $\mathrm{sys}(X,\dist)\leq 2\delta(X,\dist)$.
\end{proof}

The next result shows that we can find a positive uniform lower bound on the systoles associated to a converging sequence of $\RCD(0,N)$-structures on $X$.

\begin{Proposition}\label{prop6.2}Assume that $\{(X,\dist_n,\m_n)\}$ converges to $(X,\dist_{\infty},\m_{\infty})$ in the mGH-topology, where, for every $n\in\N\cup\{\infty\}$, $(X,\dist_n,\m_n)$ is an $\RCD(0,N)$-structure on $X$. Then $
0<\inf_{n\in\N}\{\delta(X,\dist_n)\}$.
\end{Proposition}

\begin{proof}First of all, observe that by Theorem \ref{3.3}, $\delta(X,\dist_n)>0$ for every $n\in\N$. In particular, it is sufficient to prove that there exists a constant $\delta>0$ such that $\delta(X,\dist_n)\geq\delta$ whenever $n$ is large enough.\\
We define $\epsilon_n\coloneqq\GH((X,\dist_n),(X,\dist_{\infty}))\to0$, $\delta_2\coloneqq\delta(X,\dist_{\infty})/2$, and $\delta_1\coloneqq\delta(X,\dist_{\infty})/3$. Whenever $n$ is large enough, we have $\delta_1>20\epsilon_n$ and $\delta_2>\delta_1+10\epsilon_n$; hence, by Theorem 3.4 of \cite{Wei-Sormani_01}, there is a surjective group homomorphism $\psi_n\colon G(\delta_1,\dist_n)\to G(\delta_2,\dist_{\infty})$. Moreover, since $\delta_2<\delta(X,\dist_{\infty})$, then $G(\delta_2,\dist_{\infty})$ is isomorphic to $\overline{\pi}_1(X)$ by Proposition \ref{3.3}. Now, fixing $\tilde{x}\in\tilde{X}$, and $x_1\in X_{\dist_n}^{\delta_1}$, such that $x\coloneqq p(\tilde{x})=p_{\dist_n}^{\delta_1}(x_1)$, we have a surjective homomorphism:
$$
q\colon \pi_1(X,x)/p_*\pi_1(\tilde{X},\tilde{x})\to(\pi_1(X,x)/p_*\pi_1(\tilde{X},\tilde{x}))/({{p}_{\dist_n}^{\delta_1}}_*\pi_1(X_{\dist_n}^{\delta_1},x_1)/{p}_*\pi_1(\tilde{X},\tilde{x})).
$$
However, the domain of $q$ is isomorphic to $\overline{\pi}_1(X)$, whereas its codomain is isomorphic to $G(\delta_1,\dist_n)$. Therefore, $q$ gives rise to a surjective homomorphism $\nu_n$ from $\overline{\pi}_1(X)$ onto $G(\delta_1,\dist_n)$. Hence, we have a surjective group homomorphism:
$$
\begin{tikzcd}
\overline{\pi}_1(X)\arrow[r,"\nu_n"]& G(\delta_1,\dist_n)\arrow[r,"\psi_n"]&G(\delta_2,\dist_{\infty})\arrow[r,"\sim"] &\overline{\pi}_1(X).
\end{tikzcd}
$$
However, $\overline{\pi}_1(X)$ is a Hopfian group by Proposition \ref{prop4.2}; thus the homomorphism above has to be an isomorphism. In particular, $\nu_n\colon\overline{\pi}_1(X)\rightarrow G(\delta_1,\dist_n)$ has to be injective; hence, it is an isomorphism, and it implies that $q$ is also an isomorphism. In particular, we necessarily have ${{p}_{\dist_n}^{\delta_1}}_*\pi_1(X_{\dist_n}^{\delta_1},x_1)={p}_*\pi_1(\tilde{X},\tilde{x})$; hence, by the classification Theorem (see Theorem 2, Chapter 2, Section 5 in \cite{Spanier_81}), $(X_{\dist_n}^{\delta_1},X,{p}_{\dist_n}^{\delta_1})$ is equivalent to $(\tilde{X},X,p)$. In particular, every ball of radius $\delta_1$ in $(X,\dist_n)$ is evenly covered by $p$; thus, $\delta(X,\dist_n)\geq\delta_1$, which concludes the proof.
\end{proof}

The following proposition is a converse to Proposition \ref{prop6.2}; it will be essential to prove the converse implication of Theorem \ref{thmA}.

\begin{Proposition}\label{prop:equidelta}Assume that $\{(\tilde{X},\tilde{\dist}_n,\tilde{\m}_n,\tilde{*}_n)\}$ converges to $(\tilde{X},\tilde{\dist}_{\infty},\tilde{\m}_{\infty},\tilde{*}_{\infty})$ in the equivariant pmGH-topology, where, for every $n\in\N\cup\{\infty\}$, $(\tilde{X},\tilde{\dist}_n,\tilde{\m}_n,\tilde{*}_n)$ is an equivariant pointed $\RCD(0,N)$-structure on $\tilde{X}$. Then $0<\inf_{n\in\N}\{\delta(X,\dist_n)\}$, where $(X,\dist_n,\m_n)$ is the push-forward of $(\tilde{X},\tilde{\dist}_n,\tilde{\m}_n)$.
\end{Proposition}

\begin{proof}We fix a sequence $\{(\tilde{f}_n,\tilde{g}_n,\phi_n,\epsilon_n)\}$ realizing the equivariant pointed convergence. Looking for a contradiction, assume that $\inf_{n\in\N}\{\mathrm{sys}(X,\dist_n)\}=0$. Without loss of generality, we can assume (passing to a subsequence if necessary) that there exist sequences $\{\tilde{x}_n\}$ in $\tilde{X}$ and $\{\gamma_n\}$ in $\overline{\pi}_1(X)\backslash\{\id\}$ such that $\tilde{\dist}_n(\gamma_n\tilde{x}_n,\tilde{x}_n)\to0$. However, when $n$ is large enough so that $\tilde{\dist}_n(\gamma_n\tilde{x}_n,\tilde{x}_n)\leq\epsilon_n^{-1}$, we have $\mathrm{sys}(X,\dist_{\infty})\leq \tilde{\dist}_{\infty}(\tilde{f}_n(\tilde{x}_n),\phi_n(\gamma_n)\tilde{f}_n(\tilde{x}_n))=\tilde{\dist}_{\infty}(\tilde{f}_n(\tilde{x}_n),\tilde{f}_n(\gamma_n\tilde{x}_n))\leq \tilde{\dist}_n(\gamma_n\tilde{x}_n,\tilde{x}_n)+\epsilon_n\to0$. Therefore, $\mathrm{sys}(X,\dist_{\infty})=0=\delta(X,\dist_{\infty})$ (using Proposition \ref{prop6.1}), which is the contradiction we were looking for. Hence $0<\inf_{n\in\N}\{\mathrm{sys}(X,\dist_n)\}$; therefore, thanks to Proposition \ref{prop6.1}, we have $0<\inf_{n\in\N}\{\delta(X,\dist_n)\}$.
\end{proof}

We can now prove Theorem \ref{thmA}.

\begin{proof}[Proof of Theorem \ref{thmA}, direct implication]

Assume that $\{\mathcal{X}_n=(X,\dist_n,\m_n,*_n)\}$ converges in the pmGH-topology to $\mathcal{X}_{\infty}=(X,\dist_{\infty},\m_{\infty},*_{\infty})$. Let us prove that $\{\tilde{\mathcal{X}}_n\}$ converges in the equivariant pmGH-topology to $\tilde{\mathcal{X}}_{\infty}$.\\

\proofpart{I}{Construction of the realizing sequence $\{\tilde{f}_n,\tilde{g}_n,\psi_n,\epsilon_n\}$}

First of all, we fix a sequence $\{f_n,g_n,\epsilon_n'\}$ realizing the convergence of $\{\mathcal{X}_n\}$ to $\mathcal{X}_{\infty}$ in the pmGH-topology. Then, we define $\delta\coloneqq\inf_{n\in\N\cup\{\infty\}}\{\delta(X,\dist_n)\}$, which satisfies $\delta>0$ thanks to Proposition \ref{prop6.2}. By Proposition \ref{prop6.1}, we have $\mu_0\coloneqq\inf_{n\in\N\cup\{\infty\}}\{\mathrm{sys}(X,\dist_n)\}= 2\delta>0$. We define $\alpha\coloneqq \delta/2$, and we assume that $n$ is large enough so that:
\begin{equation}\label{eqthA.1}
5\epsilon_n'<\alpha<\mu_0/2-3\epsilon_n'/2.
\end{equation}
Theorem \ref{3.3} implies that $(\tilde{X},\tilde{\dist}_{\infty},\tilde{\m}_{\infty})$ is isomorphic to $(X_{\dist_{\infty}}^{\alpha},\dist_{\infty,\alpha},\m_{\infty,\alpha})$, since $\alpha<\delta(X,\dist_{\infty})$. Now, thanks to Theorem 16 of \cite{Reviron_08} (and the construction in its proof), there exists a triple $(\tilde{f}_n,\tilde{g}_n,\psi_n)$ such that: 
\begin{itemize}
\item$\tilde{f}_n\colon(\tilde{X},\tilde{*}_n)\to(\tilde{X},\tilde{*}_{\infty})$ (resp. $\tilde{g}_n\colon(\tilde{X},\tilde{*}_{\infty})\to(\tilde{X},\tilde{*}_{n})$) satisfy $p\circ\tilde{f}_n=f_n\circ p$ (resp. $p\circ\tilde{g}_n=g_n\circ p$),
\item for every $\tilde{x}\in\tilde{X}$, we have $\tilde{\dist}_n(\tilde{g}_n\circ\tilde{f}_n(\tilde{x}),\tilde{x})\leq\epsilon_n'$ and $\tilde{\dist}_{\infty}(\tilde{f}_n\circ\tilde{g}_n(\tilde{x}),\tilde{x})\leq\epsilon_n'$,
\item for every $\tilde{x}\in\tilde{X}$ and $\eta\in\overline{\pi}_1(X)$, we have $\tilde{f}_n(\eta\cdot\tilde{x})=\psi_n(\eta)\cdot\tilde{f}_n(\tilde{x})$ and $\tilde{g}_n(\eta\cdot\tilde{x})=\psi_n^{-1}(\eta)\cdot\tilde{g}_n(\tilde{x})$.
\end{itemize}
Moreover, using inequality \ref{eqthA.1}, Theorem 16 of \cite{Reviron_08} assures that, for every $\tilde{x},\tilde{y}\in\tilde{X}$, we have:
$$
\lvert\tilde{\dist}_{\infty}(\tilde{f}_n(\tilde{x}),\tilde{f}_n(\tilde{y}))-\tilde{\dist}_n(\tilde{x},\tilde{y})\rvert \leq 3\epsilon_n'(\tilde{\dist}_n(\tilde{x},\tilde{y})/\alpha+1),
$$
and:
$$
\lvert\tilde{\dist}_{n}(\tilde{g}_n(\tilde{x}),\tilde{g}_n(\tilde{y}))-\tilde{\dist}_{\infty}(\tilde{x},\tilde{y})\rvert \leq 3\epsilon_n'(\tilde{\dist}_{\infty}(\tilde{x},\tilde{y})/\alpha+1).
$$
We fix $C>0$ such that $C+3/\alpha\leq C^2$, and we define $\epsilon_n\coloneqq C\sqrt{\epsilon_n'}$. When $n$ is large enough so that $\epsilon_n'\leq\epsilon_n$, we have:
\begin{equation}\label{proofthmAeq1}
(\tilde{f}_n,\tilde{g}_n,\psi_n,\epsilon_n)\ \text{satisfy point (i) to (iv) of Definition \ref{def:equivtopo} w.r.t.}\ \tilde{\mathcal{X}}_n\ \text{and}\ \tilde{\mathcal{X}}_{\infty}.
\end{equation}
\\
Let us prove that, when $n$ is large enough, $\tilde{f}_n$ and $\tilde{g}_n$ are Borel maps. Let $\tilde{x}\in\tilde{X}$, and let $r<\delta/3$. Thanks to Proposition \ref{prop3.4} and property \ref{proofthmAeq1}, we easily get:
\begin{equation*}
\tilde{f}_n^{-1}(\tilde{B}_{\infty}(\tilde{x},r))= (f_n\circ p)^{-1}(B_{\infty}(x,r))\cap\tilde{B}_n(\tilde{g}_n(\tilde{x}),r+2\epsilon_n),
\end{equation*}
when $n$ is large enough so that $\delta/3+4\epsilon_n<\delta/2<\epsilon_n^{-1}$, and where $x\coloneqq p(\tilde{x})$. However, $f_n$ is a Borel map, and $p$ is continuous; therefore $\tilde{f}_n^{-1}(\tilde{B}_{\infty}(\tilde{x},r))$ is a Borel subset of $\tilde{X}$. We have shown that when $n$ is large enough, the pre-image by $\tilde{f}_n$ of balls of radius $r<\delta/3$ are Borel subsets of $\tilde{X}$. Therefore, for $n$ large enough, $\tilde{f}_n$ is a Borel map, and the same is true for $\tilde{g}_n$ with the same procedure.\\

\proofpart{II}{Measured convergence}

Our goal here is to prove that (making $\epsilon_n$ larger if necessary but keeping $\epsilon_n\to0$), we have:
$$
\max\{\dist_{\mathcal{P}}^{\{{\epsilon_n}^{-1}\}}(\tilde{f}_{n_*}\tilde{\m}_n,\tilde{\m}_{\infty}),\dist_{\mathcal{P}}^{\{{\epsilon_n}^{-1}\}}(\tilde{g}_{n_*}\tilde{\m}_{\infty},\tilde{\m}_{n})\}\leq\epsilon_n.
$$
This is implied by the fact that, $\{\dist_{\mathcal{P}}^{\{R\}}(\tilde{f}_{n_*}\tilde{\m}_n,\tilde{\m}_{\infty})\}$ and $\{\dist_{\mathcal{P}}^{\{R\}}(\tilde{g}_{n_*}\tilde{\m}_{\infty},\tilde{\m}_{n})\}$ converge to $0$ as $n$ goes to infinity, for every $R>0$.\\

First of all, observe that $\lim_{n\to\infty}\dist_{\mathcal{P}}^{\{R\}}(\tilde{f}_{n_*}\tilde{\m}_n,\tilde{\m}_{\infty})=0$ for every $R>0$ if and only if $\{\tilde{f}_{n_*}\tilde{\m}_n\}$ converge to $\tilde{\m}_{\infty}$ in the weak-$*$ topology. Then, note that the space $\mathcal{M}_{\mathrm{loc}}(\tilde{X},\tilde{\dist}_{\infty})$ of Radon measures on $(\tilde{X},\tilde{\dist}_{\infty})$ endowed with the weak-$*$ topology is metrizable (see Theorem A2.6.III in \cite{Daley-Vere-Jones_03}). Hence, it is sufficient to show that any subsequence of $\{\tilde{f}_{n_*}\tilde{\m}_n\}$ admits a subsequence converging to $\tilde{\m}_{\infty}$. Without loss of generality (reindexing the sequence if necessary), let us just show that $\{\tilde{f}_{n_*}\tilde{\m}_n\}$ admits a subsequence converging to $\tilde{\m}_{\infty}$.\\

First, let us show that $\{\tilde{f}_{n_*}\tilde{\m}_n\}$ is precompact, which is implied by the uniform boundedness of $$\{\tilde{f}_{n_*}\tilde{\m}_n({B}_{\tilde{\dist}_{\infty}}(R))\},$$ for every $R>0$ (see Theorem A2.6.IV and Theorem A2.4.I in \cite{Daley-Vere-Jones_03}). We define 
\[r_0\coloneqq\inf_{n\in\N\cup\{\infty\}}\{\delta(X,\dist_n)\}/2 \quad \text{and} \quad M\coloneqq\sup_{n\in\N\cup\{\infty\}}\{\m_n(X)\}.\]
 Observe that $r_0$ is positive thanks to Proposition \ref{prop6.2}, and $M$ is finite since $\{\m_n(X)\}$ converges to $\m_{\infty}(X)$, which is finite. Thanks to point (v) of Proposition \ref{prop3.4}, we have $\tilde{\m}_{n}(B_{\tilde{\dist}_{n}}(r_0))=\m_{n}(B_{\dist_{n}}(r_0))\leq M$, for every $n\in\N$. Then, thanks to property \ref{proofthmAeq1}, we have $\tilde{f}_n^{-1}({B}_{\tilde{\dist}_{\infty}}(R))\subset B_{\tilde{\dist}_n}(2R)$, for every $R>0$, and $n$ sufficiently large. Now, consider the following two cases:
\begin{itemize}
\item if $R\leq r_0/2$, we get $\tilde{f}_{{n_*}}\tilde{\m}_{n}({B}_{\tilde{\dist}_{\infty}}(R))\leq \tilde{\m}_{n}(B_{\tilde{\dist}_{n}}(2R))\leq \tilde{\m}_{n}(B_{\tilde{\dist}_{n}}(r_0))\leq M$, when $n$ is sufficiently large, 
\item if $R>r_0/2$, thanks to Bishop--Gromov inequality for $\RCD(0,N)$ spaces (see Theorem 6.2 in \cite{Bacher-Sturm_10}), we get $\tilde{f}_{{n_*}}\tilde{\m}_{n}({B}_{\tilde{\dist}_{\infty}}(R))\leq\tilde{\m}_{n}(B_{\tilde{\dist}_{n}}(2R))\leq M(2R/r_0)^N$, when $n$ is sufficiently large.
\end{itemize}
In particular, for every $R>0$, the sequence $\{\tilde{f}_{{n_*}}\tilde{\m}_{n}({B}_{\tilde{\dist}_{\infty}}(R))\}$ is uniformly bounded; hence $\{\tilde{f}_{{n_*}}\tilde{\m}_{n}\}$ is precompact.\\
Now, passing to a subsequence if necessary, we can assume that $\{\tilde{f}_{{n_*}}\tilde{\m}_{n}\}$ is converging to some $\m\in\mathcal{M}_{\mathrm{loc}}(\tilde{X},\tilde{\dist}_{\infty})$. Let us show that $\m=\tilde{\m}_{\infty}$. Note that it is sufficient to prove that, for every $\tilde{x}\in\tilde{X}$ and $0<r< r_0$, we have $\m(B_{\tilde{\dist}_{\infty}}(\tilde{x},r))=\tilde{\m}_{\infty}(B_{\tilde{\dist}_{\infty}}(\tilde{x},r))$; since small balls generate the Borel $\sigma$-algebra of $\tilde{X}$.\\
First, observe that, for every $n\in\N$, we have $\tilde{\m}_{n}(B_{\tilde{\dist}_{n}}(r_0))=\m_{n}(B_{\dist_{n}}(r_0))\geq m$, where
 \[ m\coloneqq\inf_{n\in\N\cup\{\infty\}}\{\m_n(B_{\dist_n}(r_0))\} .\]
 In addition, $m$ is positive since $\{\m_{n}(B_{\dist_n}(r_0))\}$ is a sequence of positive numbers converging to $\m_{\infty}(B_{\dist_{\infty}}(r_0))$, which is positive. Therefore, $\{\tilde{\mathcal{X}}_n\}$ is a sequence of pointed $\RCD(0,N)$ spaces with measures uniformly bounded from below; hence (thanks to Theorem 7.2 in \cite{Gigli-Mondino-Savare_15}), any limit point in the pmGH-topology is a full support $\RCD(0,N)$ space. However, the sequence converges in the pmGH-topology to $(\tilde{X},\tilde{\dist}_{\infty},\m,\tilde{*}_{\infty})$. Thus, $(\tilde{X},\tilde{\dist}_{\infty},\m,\tilde{*}_{\infty})$ is a full support $\RCD(0,N)$ space. In particular, thanks to Theorem 30.11 in \cite{Villani_09}, we have $\m(\partial B_{\tilde{\dist}_{\infty}}(\tilde{x},R))=0$ for every $R>0$ and $\tilde{x}\in\tilde{X}$. Hence, thanks to Proposition A2.6.II in \cite{Daley-Vere-Jones_03}, for every $R>0$ and $\tilde{x}\in\tilde{X}$ we have:
\begin{equation}\label{eqn}
\m(B_{\tilde{\dist}_{\infty}}(\tilde{x},R))=\lim_{n\to\infty}\tilde{f}_{{n_*}}\tilde{\m}_{n}(B_{\tilde{\dist}_{\infty}}(\tilde{x},R)),
\end{equation}
Now, let $\tilde{x}\in\tilde{X}$ and $0<r< r_0$, and let us show that we have $\m(B_{\tilde{\dist}_{\infty}}(\tilde{x},r))=\tilde{\m}_{\infty}(B_{\tilde{\dist}_{\infty}}(\tilde{x},r))$. First, when $n$ is large enough so that $r\leq \epsilon_n^{-1}$, we can use property \ref{proofthmAeq1} to get:
$$
B_{\tilde{\dist}_{n}}(\tilde{g}_{n}(\tilde{x}),r-2\epsilon_n)\subset \tilde{f}_{n}^{-1}(B_{\tilde{\dist}_{\infty}}(\tilde{x},r))\subset B_{\tilde{\dist}_{n}}(\tilde{g}_{n}(\tilde{x}),r+2\epsilon_n)).
$$
In particular, defining $A\coloneqq\m(B_{\tilde{\dist}_{\infty}}(\tilde{x},r))$ and using equation \ref{eqn}, we have:
$$
\limsup_{n\to\infty}\tilde{\m}_{n}(B_{\tilde{\dist}_{n}}(\tilde{g}_{n}(\tilde{x}),r-2\epsilon_n))\leq A\leq\liminf_{n\to\infty}\tilde{\m}_{n}(B_{\tilde{\dist}_{n}}(\tilde{g}_{n}(\tilde{x}),r+2\epsilon_n)).
$$
Moreover, when $n$ is large enough, we have $r+2\epsilon_n<r_0<\delta/2$; hence, point (v) of Proposition \ref{prop3.4} implies:
$$
\limsup_{n\to\infty}\m_{n}({B}_{\dist_{n}}(g_{n}(x),r-2\epsilon_n))\leq A\leq\liminf_{n\to\infty}\m_{n}({B}_{\dist_{n}}(g_{n}(x),r+2\epsilon_n)),
$$
where $x\coloneqq p(\tilde{x})$. Now, observe that when $n$ is large enough so that $r+4\epsilon_n\leq\epsilon_n^{-1}$, we can use property \ref{proofthmAeq1} to get:
\begin{align*}
{B}_{\dist_{n}}(g_{n}(x),r+2\epsilon_n)&\subset f_{n}^{-1}({B}_{\dist_{\infty}}(x,r+4\epsilon_n)),\\
f_{n}^{-1}({B}_{\dist_{\infty}}(x,r-4\epsilon_n))&\subset {B}_{\dist_{n}}(g_{n}(x),r-2\epsilon_n).
\end{align*}
In particular, for every $\eta>0$, we have:
$$
\limsup_{n\to\infty}f_{{n_*}}\m_{n}({B}_{\dist_{\infty}}(x,r-\eta))\leq A\leq\liminf_{n\to\infty}f_{{n_*}}\m_{n}({B}_{\dist_{\infty}}(x,r+\eta)).
$$
However, since $\{f_{n_*}\m_n\}$ converges to $\m_{\infty}$, and since $\mathcal{X}_{\infty}$ is a full support $\RCD(0,N)$ space, we can apply Theorem 30.11 in \cite{Villani_09} and Proposition A2.3.II in \cite{Daley-Vere-Jones_03} to get:
\begin{align*}
\limsup_{n\to\infty}f_{{n_*}}\m_{n}({B}_{\dist_{\infty}}(x,r-\eta))&=\m_{\infty}(B_{\dist_{\infty}}(x,r-\eta)),\\
\liminf_{n\to\infty}f_{{n_*}}\m_{n}({B}_{\dist_{\infty}}(x,r+\eta))&=\m_{\infty}(B_{\dist_{\infty}}(x,r+\eta)).
\end{align*}
Hence, for every $\eta>0$, we have:
$$
\m_{\infty}(B_{\dist_{\infty}}(x,r-\eta))\leq\m(B_{\tilde{\dist}_{\infty}}(\tilde{x},r))\leq\m_{\infty}(B_{\dist_{\infty}}(x,r+\eta));
$$
and, letting $\eta$ goes to $0$, we have $\m(B_{\tilde{\dist}_{\infty}}(\tilde{x},r))=\m_{\infty}(B_{\dist_{\infty}}(x,r))=\tilde{\m}_{\infty}(B_{\tilde{\dist}_{\infty}}(\tilde{x},r))$ (using $r<r_0<\delta/2$ for the last equality). Therefore $\{\tilde{f}_{n_*}\tilde{\m}_n\}$ converges to $\tilde{\m}_{\infty}$.\\

For every $R'>0$, we define $\epsilon(n,R')\coloneqq\dist_{\mathcal{P}}^{\{R'\}}(\tilde{f}_{n_*}\tilde{\m}_n,\tilde{\m}_{\infty})$. Thanks to the discussion above, we have $\lim_{n\to\infty}\epsilon(n,R')\to0$, for every $R'>0$. Let $R>0$, and let us show that $\lim_{n\to\infty}\dist_{\mathcal{P}}^{\{R\}}(\tilde{g}_{n_*}\tilde{\m}_{\infty},\tilde{\m}_{n})=0$.\\
Let $A\subset \overline{B}_{\tilde{\dist}_n}(R)$, and observe that we have $\tilde{\m}_n(A)\leq \tilde{f}_{n_*}\tilde{\m}_n(\tilde{f}_n(A))$. Also, when $n$ is large enough, we can use property \ref{proofthmAeq1} to get $\tilde{f}_n(A)\subset \overline{B}_{\tilde{\dist}_{\infty}}(2R)$. Therefore, we have $\tilde{\m}_n(A)\leq\tilde{\m}_{\infty}(\{\tilde{f}_n(A)\}^{\epsilon(n,2R)})+\epsilon(n,2R)$. Then, when $n$ is large enough, we can use property \ref{proofthmAeq1} to obtain $\{\tilde{f}_n(A)\}^{\epsilon(n,2R)}\subset\tilde{g}_n^{-1}(\{A\}^{2\epsilon_n+\epsilon(n,2R)})$. Thus, we have $\tilde{\m}_n(A)\leq\tilde{g}_{n_*}\tilde{\m}_{\infty}(\{A\}^{2\epsilon_n+\epsilon(n,2R)})+\epsilon(n,2R)$. Applying the same arguments, we also have $\tilde{g}_{n_*}\tilde{\m}_{\infty}(A)\leq \tilde{\m}_n(\{A\}^{2\epsilon_n+\epsilon(n,2R)})+\epsilon(n,2R)$. Therefore, $\dist_{\mathcal{P}}^{\{R\}}(\tilde{g}_{n_*}\tilde{\m}_{\infty},\tilde{\m}_{n})\leq\epsilon(n,2R)+2\epsilon_n$; in particular, $\lim_{n\to\infty}\dist_{\mathcal{P}}^{\{R\}}(\tilde{g}_{n_*}\tilde{\m}_{\infty},\tilde{\m}_{n})=0$. This concludes the proof.
\end{proof}

\begin{proof}[Proof of Theorem \ref{thmA}, converse implication]Assume that $\{\tilde{\mathcal{X}}_n=(\tilde{X},\tilde{\dist}_n,\tilde{\m}_n,\tilde{*}_n)\}$ converges in the equivariant pmGH-topology to $\tilde{\mathcal{X}}_{\infty}=(\tilde{X},\tilde{\dist}_{\infty},\tilde{\m}_{\infty},\tilde{*}_{\infty})$. Let us prove that $\{{\mathcal{X}}_n=({X},{\dist}_n,{\m}_n,{*}_n)\}$ converges in the pmGH-topology to ${\mathcal{X}}_{\infty}=({X},{\dist}_{\infty},{\m}_{\infty},{*}_{\infty})$.\\

Let $\{\tilde{f}_n,\tilde{g}_n,\phi_n,\epsilon_n\}$ be a sequence realizing the convergence of $\{\tilde{\mathcal{X}}_n\}$ to $\tilde{\mathcal{X}}_{\infty}$ in the equivariant pmGH-topology. Thanks to the equivariant requirement, there exists pointed Borel maps $f_n\colon (X,*_n)\to(X,*_{\infty})$ and $g_n\colon (X,*_{\infty})\to(X,*_n)$ such that $p\circ\tilde{f}_n=f_n\circ p$ and $p\circ \tilde{g}_n=g_n\circ p$.\\

Let us fix $x\in X$ and $\tilde{x}\in p^{-1}(x)$. Observe that $\dist_n(g_n(f_n(x)),x)=\inf\{\tilde{\dist}_n(\tilde{y},\tilde{x})\}$, where the infimum is taken over all $\tilde{y}\in\tilde{X}$ such that $p(\tilde{y})=g_n(f_n(x))$. However, we have $p(\tilde{g}_n(\tilde{f}_n(\tilde{x})))=g_n(f_n(x))$. Therefore, we have $\dist_n(g_n(f_n(x)),x)\leq \tilde{\dist}_n(\tilde{g}_n(\tilde{f}_n(\tilde{x})),\tilde{x})\leq \epsilon_n$. The same argument shows that $\dist_{\infty}(f_n(g_n(x)),x)\leq\epsilon_n$.\\

Let $y_i\in X$ ($i\in\{1,2\}$) and let $\tilde{y}_i$ such that $p(\tilde{y}_i)=y_i$ and $\tilde{\dist}_{\infty}(\tilde{y}_1,\tilde{y}_2)=\dist_{\infty}(y_1,y_2)$. Assume that $D_{\infty}\coloneqq\Diam(X,\dist_{\infty})\leq\epsilon_n^{-1}$ and observe that, since $p(\tilde{g}_n(\tilde{y}_i))=g_n(y_i)$, we have:
\begin{align}\label{equiv4} 
\dist_n(g_n(y_1),g_n(y_2))-\dist_{\infty}(y_1,y_2)&\leq \tilde{\dist}_n(\tilde{g}_n(\tilde{y}_1),\tilde{g}_n(\tilde{y}_2))-\tilde{\dist}_{\infty}(\tilde{y}_1,\tilde{y}_2)\\
&\leq\epsilon_n.
\end{align}
Then, let $x_i\in X$ ($i\in\{1,2\}$) such that $\dist_n(x_1,x_2)\leq\epsilon_n^{-1}$, and fix $\tilde{x}_i\in\tilde{X}$ such that $p(\tilde{x}_i)=x_i$ and $\tilde{\dist}_{n}(\tilde{x}_1,\tilde{x}_2)=\dist_n(x_1,x_2)$. Observe that we have $p(\tilde{f}_n(\tilde{x}_i))=f_n(x_i)$, therefore:
\begin{align}\label{equiv5} 
\dist_{\infty}(f_n(x_1),f_n(x_2))-\dist_n(x_1,x_2)&\leq \tilde{\dist}_{\infty}(\tilde{f}_n(\tilde{x}_1),\tilde{f}_n(\tilde{x}_2))-\tilde{\dist}_n(\tilde{x}_1,\tilde{x}_2)\\
&\leq\epsilon_n.
\end{align}
Let us show that $\{D_n\coloneqq\Diam(X,\dist_n)\}$ is a bounded sequence. Let $x_i\in X$ ($i\in\{1,2\}$) and observe that thanks to inequality \ref{equiv4}, we have $\dist_n(g_n(f_n(x_1)),g_n(f_n(x_2)))\leq\epsilon_n+D_{\infty}$ (when $D_{\infty}\leq\epsilon_n^{-1}$). However, we have $\lvert \dist_n(x_1,x_2)-\dist_n(g_n(f_n(x_1)),g_n(f_n(x_2)))\rvert\leq \dist_n(g_n(f_n(x_1)),x_1)+\dist_n(g_n(f_n(x_2)),x_2)\leq 2\epsilon_n$. Therefore, $\dist_n(x_1,x_2)\leq D_{\infty}+3\epsilon_n$. We can conclude that $\{D_n\}$ is bounded.\\
Since $\{D_n\}$ is bounded, we have (thanks to inequality \ref{equiv5}):
$$
\forall x_1,x_2\in X,\dist_{\infty}(f_n(x_1),f_n(x_2))-\dist_n(x_1,x_2)\leq \epsilon_n,
$$
when $n$ is large enough. Also, we have $\dist_{n}(x_1,x_2)\leq\dist_{n}(g_n(f_n(x_1)),g_n(f_n(x_2)))+2\epsilon_n$. Therefor, using inequality \ref{equiv4} we obtain, $\dist_{n}(x_1,x_2)-\dist_{\infty}(f_n(x_1),f_n(x_2))\leq 3\epsilon_n$. Hence, we can conclude that $\Dis(f_n)\leq 3\epsilon_n$. The same argument also gives $\Dis(g_n)\leq3\epsilon_n$, which concludes the proof of the second metric requirement.\\

Finally, using Lemma \ref{prop:equidelta}, and applying exactly the same procedure as in Part II of the direct implication, we can prove that (making $\epsilon_n$ smaller if necessary but keeping $\epsilon_n\to 0$) we have: 
\[ \max\{\dist_{\mathcal{P}}(f_{n_*}\m_n,\m_{\infty}),\dist_{\mathcal{P}}(g_{n_*}\m_{\infty},\m_{n})\}\leq\epsilon_n\, .\]
 Hence, $\{f_n,g_n,\epsilon_n\}$ is a sequence realizing the convergence of $\{\mathcal{X}_n\}$ to $\mathcal{X}_{\infty}$ in the pmGH-topology, which concludes the proof.

\end{proof}

\subsection{Proof of Theorem \ref{thmB}}\label{proofB}

In this section, we give a proof of Theorem \ref{thmB}. Let $\{(X,\dist_n,\m_n)\}$ be a sequence converging in the mGH-topology to $(X,\dist_{\infty},\m_{\infty})$, where for every $n\in\N\cup\{\infty\}$, $(X,\dist_n,\m_n)$ is an $\RCD(0,N)$-structure on $X$. For every $n\in\N\cup\{\infty\}$, we fix $\phi_n$ a splitting of $(\tilde{{X}},\tilde{\dist}_n,\tilde{\m}_n)$ with soul $(\overline{{X}}_n,\overline{\dist}_n,\overline{\m}_n)$, and we denote $k\coloneqq k(X)$ the splitting degree of $X$ (see Corollary \ref{cor5.1} for the definition of the splitting degree). To conclude, we are going to prove that:
\begin{equation}\label{property1}
\{\overline{{X}}_n,\overline{\dist}_n,\overline{\m}_n\}\text{ converges to }(\overline{{X}}_{\infty},\overline{\dist}_{\infty},\overline{\m}_{\infty})\text{ in the mGH-topology},
\end{equation}
and:
\begin{equation}\label{property2}
\{\R^k/\Gamma(\phi_n),\dist_{\Gamma(\phi_n)}\}\text{ converges to }(\R^k/\Gamma(\phi_{\infty}),\dist_{\Gamma(\phi_{\infty})})\text{ in the GH-topology}.\end{equation}
Observe that, since $X$ is compact, we can find a family $\{*_n\}_{n\in\N\cup\{\infty\}}$ of points in $X$ such that $\{(X,\dist_n,\m_n,*_n)\}$ converges to $(X,\dist_{\infty},\m_{\infty},*_{\infty})$ in the pmGH-topology. Then, for every $n\in\N\cup\{\infty\}$, let us fix $\tilde{*}_n\in p^{-1}(*_n)$. Observe that without loss of generality, we can assume that, for every $n\in\N\cup\{\infty\}$, we have $p_{\R^k}(\phi_n(\tilde{*}_n))=0$. For every $n\in\N\cup\{\infty\}$, we denote:
\begin{itemize}
\item $\mathcal{X}_n\coloneqq(X,\dist_n,\m_n,*_n)$
\item $\tilde{\mathcal{X}}_n\coloneqq(\tilde{{X}},\tilde{\dist}_n,\tilde{\m}_n,\tilde{*}_n)=p_{\tilde{*}_n}^*(\mathcal{X}_n)$, \item $\overline{\mathcal{X}}_n\coloneqq(\overline{{X}}_n,\overline{\dist}_n,\overline{\m}_n,\overline{*}_n)$, where $\overline{*}_n\coloneqq p_{\overline{X}_n}(\phi_n(\tilde{*}_n))$. 
\end{itemize}
Thanks to Theorem \ref{thmA}, $\{\tilde{\mathcal{X}}_n\}$ converges to $\tilde{\mathcal{X}}_{\infty}$ in the equivariant pmGH-topology. Thanks to the proof of Theorem \ref{thmA}, there exists a sequence $\{f_n,g_n,\epsilon_n\}$ (resp. $\{\tilde{f}_n,\tilde{g}_n,\phi_n,\epsilon_n\}$) realizing the convergence of $\{\mathcal{X}_n\}$ to $\mathcal{X}_{\infty}$ (resp. of $\{\tilde{\mathcal{X}}_n\}$ to $\tilde{\mathcal{X}}_{\infty}$) in the pmGH-topology (resp. in the equivariant pmGH-topology), such that $p\circ\tilde{f}_n=f_n\circ p$ and $p\circ\tilde{g}_n=g_n\circ p$. Finally, for every $n\in\N$, we define:
\begin{itemize}
\item $k_n\coloneqq \phi_{\infty}\circ\tilde{f}_n\circ\phi_n^{-1}$, $k_n^{\R}\coloneqq p_{\R^k}\circ k_n(\overline{*}_n,\cdot)$, and $k_n^{S}\coloneqq p_{\overline{X}_{\infty}}\circ k_n(\cdot, 0)$,
\item $l_n\coloneqq \phi_{n}\circ\tilde{g}_n\circ\phi_{\infty}^{-1}$, $l_n^{\R}\coloneqq p_{\R^k}\circ l_n(\overline{*}_{\infty},\cdot)$, and $l_n^{S}\coloneqq p_{\overline{X}_{n}}\circ l_n(\cdot, 0)$.
\end{itemize}
The main difficulty of the proof will be to prove that $k_n$ and $l_n$ almost split. More precisely, we will show that for we have $k_n\simeq (k_n^{S},k_n^{\R})$ and $l_n\simeq (l_n^{S},l_n^{\R})$ (where we will give a precise meaning to $\simeq$). Then, we will deduce property \ref{property1} and property \ref{property2} from that.\\

First of all, we prove that $\{\Diam(\overline{X},\overline{\dist}_n)\}$ is bounded.

\begin{Proposition}\label{prop7.1}The sequence $\{\Diam(\overline{X},\overline{\dist}_n)\}$ is bounded.
\end{Proposition}

\begin{proof}Looking for a contradiction, let us suppose that $\limsup_{n\to\infty}\Diam(\overline{X}_n,\overline{\dist}_n)=\infty$. Passing to a subsequence if necessary, we can assume that $\Diam(\overline{X}_n,\overline{\dist}_n)> 2^{n+1}$, for every $n\in\N$. Hence, there are sequences $\{\overline{x}_n\}$ and $\{\overline{z}_n\}$ such that, for every $n\in\N$, we have $\overline{x}_n,\overline{z}_n\in\overline{X}_n$, and $\overline{d}_n(\overline{x}_n,\overline{z}_n)=2^{n+1}$.\\
For every $n\in\N$, let $\overline{\gamma}_n\colon[-2^n,2^n]\to\overline{X}_n$ be a minimizing geodesic parametrized by arc length from $\overline{x}_n$ to $\overline{z}_n$, and let us denote $\tilde{\gamma}_n\coloneqq(\overline{\gamma}_n,0)$. Thanks to Proposition \ref{prop3.6}, there exists $\eta\in\overline{\pi}_1(X)$ such that $\eta\tilde{\gamma}_n(0)\in B_{\overline{X}_n\times\R^k}((\overline{*}_n,0),D)$, where $D\coloneqq\sup_{n\in\N\cup\{\infty\}}\{\Diam(X,\dist_n)\}<\infty$ ($D$ being finite because $(X,\dist_n)$ converges to $(X,\dist_{\infty})$ in the GH-topology).\\
Then, let us define $\tilde{\beta}_n\coloneqq \eta\tilde{\gamma}_n$, and denote $\overline{\beta}_n\coloneqq p_{\overline{X}_n}(\tilde{\beta}_n)$, and $v_n\coloneqq p_{\R^k}(\tilde{\beta}_n)=\rho_{\R}^{\phi_n}(\eta)(0)$. Observe that the sequence $\{\tilde{\beta}_n\}$ consists of isometric embeddings such that $\tilde{\beta}_n(0)\in B_{\overline{X}_n\times\R^k}((\overline{*}_n,0),D)$. Moreover, $\{k_n,l_n\}$ realizes the convergence of $\{\overline{\mathcal{X}}_n\times(\R^k,0)\}$ to $(\overline{\mathcal{X}}_{\infty}\times(\R^k,0))$ in the pmGH-topology. Therefore, thanks to Arzel\`{a}--Ascoli Theorem (see Proposition 27.20 in \cite{Villani_09}), we can assume (passing to a subsequence if necessary) that $\{k_n\circ\tilde{\beta}_n\}$ converges locally uniformly to an isometric embedding $\tilde{\beta}\colon\R\to\overline{X}_{\infty}\times\R^k$. However, $(\overline{X}_{\infty},\overline{\dist}_{\infty})$ is compact; thus, applying Lemma 1 of \cite{Shen-Wei_91}, there exist $a,b\in\R^k$ and $\overline{y}_{\infty}\in\overline{X}_{\infty}$ such that, for every $t\in\R$, $\tilde{\beta}(t)=(\overline{y}_{\infty},at+b)$, and $\lVert a\rVert =1$.\\
Now, we define $\overline{y}_n\coloneqq\overline{\beta}_n(0)$ and, for $u\in\R^k$, $\Phi_n(u)\coloneqq(\overline{y}_n,u)\in\overline{X}_{n}\times\R^k$. Observe that $\{\Phi_n\}$ is a sequence of isometric embeddings such that, for every $n\in\N$, we have $\Phi_n(0)\in B_{\overline{X}_n\times\R^k}((\overline{*}_n,0),D)$. Therefore, thanks to Arzel\`{a}--Ascoli Theorem (see Proposition 27.20 in \cite{Villani_09}), we can assume (passing to a subsequence if necessary) that $\{k_n\circ\Phi_n\}$ converges locally uniformly to an isometric embedding $\Phi\colon\R^k\to\overline{X}_{\infty}\times\R^k$. Moreover, since $\overline{X}_{\infty}$ is compact, we can easily deduce from Lemma 1 of \cite{Shen-Wei_91} that there exist $\phi\in\Iso(\R^k)$ and $\overline{z}_{\infty}\in\overline{X}_{\infty}$ such that $\Phi(t)=(\overline{z}_{\infty},\phi(t))$, for every $t\in\R^k$.\\
Notice that $\tilde{\beta}_n(0)=(\overline{y}_n,v_n)=\Phi_n(v_n)$. Moreover, observe that $\lvert v_n\rvert \leq D$; hence, passing to a subsequence if necessary, we can assume that $v_n\to v\in\R^k$. Thus, we have $\dist_{\overline{X}_{\infty}\times\R^k}(\Phi(v),\tilde{\beta}(0))=\lim_{n\to\infty}\dist_{\overline{X}_{\infty}\times\R^k}(k_n\circ \Phi_n(v_n),k_n\circ\tilde{\beta}_n(0))=0$. In particular, $\overline{y}_{\infty}=\overline{z}_{\infty}$, and $\phi(v)=b$. Now, let $c\in\R^{k}$ such that $[\phi-\phi(0)](c)=a$; thus, we have $\Phi(ct+v)=at+b=\tilde{\beta}(t)$, for every $t\in\R$.\\
Now, observe that we have $0=\dist_{\overline{X}_{\infty}\times\R^k}(\Phi(c+v),\tilde{\beta}(1))$. Therefore:
$$
\dist_{\overline{X}_{n}\times\R^k}(\Phi_n(c+v_n), \tilde{\beta}_n(1))=(1+\lVert c\rVert^2)^{1/2}\to0.
$$
Therefore, we have $0=(1+\lVert c\rVert^2)^{1/2}>0$, which is the contradiction we were looking for.
\end{proof}

Thanks to Proposition \ref{prop7.1}, we can introduce the following notations.

\begin{Notation}\label{not9.1}We denote $D\coloneqq\sup_{n\in\N\cup{\infty}}\{\Diam(X,\dist_n)\}<\infty$ (finiteness being granted by the convergence of $\{X,\dist_n\}$ to $(X,\dist_{\infty})$ in the GH-topology). We also denote $\overline{D}\coloneqq\sup_{n\in\N\cup\{\infty\}}\{\Diam(\overline{X}_n,\overline{\dist}_n)\}<\infty$.
\end{Notation}

Our first goal will be to obtain a convergence result on the following ``splitting quantities''.

\begin{Notation}[Splitting quantities]\label{not10}Given $n\in\N$ and $R>0$, we define:
\begin{itemize}
\item[(i)] $\alpha(n,R)\coloneqq\sup\{\dist_{\overline{X}_{\infty}\times\R^k}(k_n(\overline{y}_n,t),(k_n^{S}(\overline{y}_n),k_n^{\R}(t)))\}$, the supremum being taken over $\overline{y}_n\in\overline{X}_n$ and $\lvert t\rvert\leq R$,
\item[(ii)] $\beta(n,R)\coloneqq\sup\{\dist_{\overline{X}_{n}\times\R^k}(l_n(\overline{y}_{\infty},t),(l_n^{S}(\overline{y}_{\infty}),l_n^{\R}(t)))\}$, the supremum being taken over $\overline{y}_{\infty}\in \overline{X}_{\infty}$ and $\lvert t\rvert\leq R$.
\end{itemize}
\end{Notation}

The following next two technical lemmas will be our main ingredients in the proof of the convergence result of the splitting quantities.

\begin{Lemma}\label{lem7.1}Let $\{\overline{y}_n\}$ be a sequence such that, for every $n\in\N$, $\overline{y}_n\in\overline{X}_n$. For every $n\in\N$, and $t\in\R^k$, we define $\Phi_n\colon t\in\R^k\rightarrow (\overline{y}_n,t)\in\overline{X}_n\times\R^k$. Then, the sequence of maps $\{k_n\circ \Phi_n\colon\R^k\to\overline{X}_{\infty}\times\R^k\}$ admits a subsequence converging locally uniformly to a map $\Phi\colon\R^k\to\overline{X}_{\infty}\times\R^k$. Moreover, for any such limit $\Phi$, there exists $\overline{y}_{\infty}$, and $\phi\in\Or_k(\R)$ such that $\forall t\in\R^k,\Phi(t)=(\overline{y}_{\infty},\phi(t))$.
\end{Lemma}

\begin{proof}Observe that, for every $n\in\N$, $\Phi_n$ is an isometric embedding that satisfies 
\[\Phi_n(0)\in B_{\overline{X}_n\times\R^k}((\overline{*}_n,0),\overline{D})\, .\]
Therefore, applying Arzel\`{a}--Ascoli Theorem (see Proposition 27.20 in \cite{Villani_09}) as in the proof of Proposition \ref{prop7.1}, we can assume without loss of generality that $\{k_n\circ\Phi_n\}$ converges locally uniformly to an isometric embedding $\Phi\colon\R^k\to\overline{X}_{\infty}\times\R^k$. Moreover, using Lemma 1 of \cite{Shen-Wei_91}, there exist $\phi\in\Iso(\R^k)$ and $\overline{y}_{\infty}\in\overline{X}_{\infty}$ such that $\Phi(t)=(\overline{y}_{\infty},\phi(t))$, for every $t\in\R^k$. To conclude, we need to show that $\phi(0)=0$. First, observe that:
$$\overline{\dist}_{\infty}(\overline{*}_{\infty},\overline{y}_{\infty})\leq\dist_{\overline{X}_{\infty}\times\R^k}((\overline{*}_{\infty},0),\Phi(0))\leq \overline{\dist}_{n}(\overline{*}_{n},\overline{y}_n)+u_n,
$$
whenever $n$ is large enough (so that $\overline{D}\leq\epsilon_n^{-1}$), and where $u_n=\epsilon_n+\dist_{\overline{X}_{\infty}\times\R^k}(\Phi(0),k_n\circ\Phi_n(0))\to0$. Now, let $t\in\R^k$ such that $\phi(t)=0$, and observe that:
$$
\overline{\dist}_{n}(\overline{*}_{n},\overline{y}_n)\leq \dist_{\overline{X}_{n}\times\R^k}((\overline{*}_{n},0),(\overline{y}_n,t))\leq \dist_{\overline{X}_{\infty}\times\R^k}((\overline{*}_{\infty},0),\Phi(t))+v_{n}=\overline{\dist}_{\infty}(\overline{*}_{\infty},\overline{y}_{\infty})+v_{n},
$$
when $n$ is large enough (so that $(\overline{D}^2+\lvert t\rvert^{2})^{1/2}\leq\epsilon_n^{-1}$), and where $v_n\coloneqq \epsilon_n+\dist_{\overline{X}_{\infty}\times\R^k}(k_n\circ\Phi_n(t),\Phi(t))\to0$. Hence, combining the two inequalities above, we obtain:
$$
\lvert \dist_{\overline{X}_{\infty}\times\R^k}((\overline{*}_{\infty},0),\Phi(0))-\overline{\dist}_{\infty}(\overline{*}_{\infty},\overline{y}_{\infty})\rvert\leq u_n+v_n\to0.
$$
In particular, $\overline{\dist}^2_{\infty}(\overline{*}_{\infty},\overline{y}_{\infty})=\overline{\dist}^2_{\infty}(\overline{*}_{\infty},\overline{y}_{\infty})+\lvert \phi(0)\rvert^2$. In conclusion, $\phi(0)=0$.
\end{proof}

\begin{Lemma}\label{lem7.2}Let $\{\overline{y}_n\}$ and $\{\overline{z}_n\}$ be sequences such that, for every $n\in\N$, $\overline{y}_n,\overline{z}_n\in\overline{X}_n$. For every $n\in\N$ and $t\in\R^k$, we define $\Phi_n(t)\coloneqq(\overline{y}_n,t)$ and $\Psi_n(t)\coloneqq(\overline{z}_n,t)$. Assume that (passing to a subsequence if necessary), the sequences of maps $\{k_n\circ\Phi_n\}$ and $\{k_n\circ\Psi_n\}$ converge locally uniformly, respectively to $\Phi=(\overline{y}_{\infty},\phi)$ and $\Psi=(\overline{z}_{\infty},\psi)$, where $\overline{y}_{\infty},\overline{z}_{\infty}\in\overline{X}_{\infty}$ and $\phi,\psi\in\Or_k(\R)$. Then, we necessarily have $\phi=\psi$.
\end{Lemma}

\begin{proof}Looking for a contradiction, let us suppose that $\phi\neq\psi$. In that case, there exists $t\in\R^k\backslash\{0\}$ such that $\phi(t)\neq\psi(t)$, which implies $\lim_{s\to\infty}\dist_{\overline{X}_{\infty}\times\R^k}(\Phi(st),\Psi(st))=\infty$. In particular, there exists $s\in\R$ such that $\dist_{\overline{X}_{\infty}\times\R^k}(\Phi(st),\Psi(st))\geq\overline{D}+1$.
However, using $\overline{\dist}_n(\overline{y}_n,\overline{z}_n)=\dist_{\overline{X}_{n}\times\R^k}(\Phi_n(st),\Psi_n(st))$, we have:
\begin{equation}\label{abc}
\dist_{\overline{X}_{\infty}\times\R^k}(\Phi(st),\Psi(st))\leq\overline{\dist}_n(\overline{y}_n,\overline{z}_n)+u_n\leq\overline{D}+u_n,
\end{equation}
where $u_n=\epsilon_n+\dist_{\overline{X}_{\infty}\times\R^k}(\Phi(st),k_n\circ\Phi_n(st))+\dist_{\overline{X}_{\infty}\times\R^k}(\Psi(st),k_n\circ\Psi_n(st))$, and when $n$ is large enough (so that $\overline{D}\leq\epsilon_n^{-1}$). Now, observe that $\lim_{n\to\infty}u_n=0$; therefore, passing to the limit in inequality \ref{abc}, we have $\dist_{\overline{X}_{\infty}\times\R^k}(\Phi(st),\Psi(st))\leq \overline{D}$, which contradicts $\dist_{\overline{X}_{\infty}\times\R^k}(\Phi(st),\Psi(st))\geq\overline{D}+1$.
\end{proof}

We can now state the convergence result on the splitting quantities.

\begin{Lemma}\label{prop7.2}For every $R>0$, we have $\lim_{n\to\infty}\alpha(n,R)=\lim_{n\to\infty}\beta(n,R)=0$.
\end{Lemma}

\begin{proof}\proofpart{I}{$\lim_{n\to\infty}\alpha(n,R)=0$}

Looking for a contradiction, we assume that $\lim_{n\to\infty}\alpha(n,R)\neq0$. Passing to a subsequence if necessary, there exist $\epsilon>0$, and sequences $\{\overline{y}_n\}$ and $\{t_n\}$ such that:
\begin{equation}\label{eqpart1}
\epsilon\leq\dist_{\overline{X}_{\infty}\times\R^k}(k_n(\overline{y}_n,t_n),(k_n^{S}(\overline{y}_n),k_n^{\R}(t_n))),
\end{equation}
$\overline{y}_n\in\overline{X}_n$, and $\lvert t_n\rvert \leq R$. Moreover, since $\{t_n\}$ is bounded, we can assume that $t_n\to t$. Now, applying Lemma \ref{lem7.1}, and passing to a subsequence if necessary, we can assume that $\{k_n\circ\Phi_n\}$ converges locally uniformly to $\Phi=(\overline{y}_{\infty},\phi)$, where $\Phi_n(s)\coloneqq(\overline{y}_n,s)$, $\phi\in\Or_k(\R)$ and $\overline{y}_{\infty}\in\overline{X}_{\infty}$. In particular, we have:
\begin{equation}\label{eqc}
\lim_{n\to\infty}\overline{\dist}_{\infty}(p_{\overline{X}_{\infty}}\circ k_n(\overline{y}_n,t_n),k_n^{S}(\overline{y}_n))=\overline{\dist}_{\infty}(p_{\overline{X}_{\infty}}\circ \Phi(t),p_{\overline{X}_{\infty}}\circ\Phi(0))=0.
\end{equation}
Now, applying Lemma \ref{lem7.1} and Lemma \ref{lem7.2}, we can also assume that $\{k_n\circ\Psi_n\}$ converges locally uniformly to $\Psi=(\overline{z}_{\infty},\phi)$, where $\Psi_n(s)\coloneqq(\overline{*}_n,s)$, and $\overline{z}_{\infty}\in\overline{X}_{\infty}$. Thus, we have:
\begin{equation}\label{eqb}
\lim_{n\to\infty}\dist_{\mathrm{eucli}}(p_{\R^k}\circ k_n(\overline{y}_n,t_n),k_n^{\R}(t_n))=\dist_{\mathrm{eucli}}(p_{\R^k}\circ\Phi(t),p_{\R^k}\circ\Psi(t))=0.
\end{equation}
Hence, using equations \ref{eqc} and \ref{eqb}, we have $\lim_{n\to\infty}\dist_{\overline{X}_{\infty}\times\R^k}(k_n(\overline{y}_n,t_n),(k_n^{S}(\overline{y}_n),k_n^{\R}(t_n)))=0$, which contradicts inequality \ref{eqpart1}.

\proofpart{II}{$\lim_{n\to\infty}\beta(n,R)=0$}

Looking for a contradiction, we assume that $\lim_{n\to\infty}\beta(n,R)\neq0$. Passing to a subsequence if necessary, there exist $\epsilon>0$, and sequences $\{\overline{y}_{\infty}^{(n)}\}\in\overline{X}_{\infty}^{\N}$ and $\lvert t_n\rvert \leq R$ such that:
\begin{equation}\label{eqpart2}
\epsilon\leq\dist_{\overline{X}_{n}\times\R^k}(l_n(\overline{y}_{\infty}^{(n)},t_n),(l_n^{S}(\overline{y}_{\infty}^{(n)}),l_n^{\R}(t_n))).
\end{equation}
Observe that $\{t_n\}$ is bounded and $\overline{X}_{\infty}$ is compact; therefore, we can assume that $t_n\to t$ and $\overline{y}_{\infty}^{(n)}\to\overline{y}_{\infty}$.\\
Now, let us define $(\overline{y}_n,s_n)\coloneqq l_n(\overline{y}_{\infty}^{(n)},t_n)$, and $\Phi_n(u)\coloneqq(\overline{y}_n,u)$, $u\in\R^k$. Applying Lemma \ref{lem7.1}, we can assume (passing to a subsequence if necessary) that $\{k_n\circ \Phi_n\}$ converges locally uniformly to $(\overline{z}_{\infty},\phi)$, where $\overline{z}_{\infty}\in\overline{X}_{\infty}$ and $\phi\in\Or_k(\R)$. Observe that $\{s_n\}$ is bounded since $\{(\overline{y}_{\infty}^{(n)},t_n)\}$ converges. Therefore, we can assume that $s_n\to s$. However, we have $\lim_{n\to\infty}k_n\circ\Phi_n(s_n)=\Phi(s)=\lim_{n\to\infty}k_n\circ l_n(\overline{y}_{\infty}^{(n)},t_n)=(\overline{y}_{\infty},t)$. Thus, $\overline{z}_{\infty}=\overline{y}_{\infty}$, and $\phi(s)=t$. Now, observe that $\overline{\dist}_{n}(p_{\overline{X}_{n}}\circ l_n(\overline{y}_{\infty}^{(n)},t_n),l_n^{S}(\overline{y}_{\infty}^{(n)}))\leq{\dist}_{\overline{X}_n\times\R^k}(\Phi_n(0),l_n(\Phi(0)))+{\dist}_{\overline{X}_n\times\R^k}(l_n(\overline{y}_{\infty},0),l_n(\overline{y}_{\infty}^{(n)},0))$.
Therefore:
\begin{equation}\label{eqf}
\lim_{n\to\infty}\overline{\dist}_{n}(p_{\overline{X}_{n}}\circ l_n(\overline{y}_{\infty}^{(n)},t_n),l_n^{S}(\overline{y}_{\infty}^{(n)}))=0.
\end{equation}
Then, applying Lemma \ref{lem7.1} and Lemma \ref{lem7.2}, we can also assume that $\{k_n\circ\Psi_n\}$ converges locally uniformly to $\Psi=(\overline{z}'_{\infty},\phi)$, where $\Psi_n(u)\coloneqq(\overline{*}_n,u)$, and $\overline{z}'_{\infty}\in\overline{X}_{\infty}$. Moreover, $\Psi(0)=\lim_{n\to\infty}k_n\circ\Psi_n(0)=(\overline{*}_{\infty},0)$; therefore, $\overline{z}'_{\infty}=\overline{*}_{\infty}$. Now, using $p_{\R^k}\circ l_n(\overline{y}_{\infty}^{(n)},t_n)=s_n$, and $l_n^{\R}(t_n)=p_{\R^k}\circ l_n(\Psi(\phi^{-1}(t_n)))$, observe that $\dist_{\mathrm{eucli}}(p_{\R^k}\circ l_n(\overline{y}_{\infty}^{(n)},t_n),l_n^{\R}(t_n))\leq\dist_{\mathrm{eucli}}(s_n,\phi^{-1}(t_n))+\dist_{\overline{X}_n\times\R^k}(\Psi_n(\phi^{-1}(t_n)),l_n(\Psi(\phi^{-1}(t_n))))$. Therefore, using $\lim_{n\to\infty}s_n=s=\phi^{-1}(t)=\lim_{n\to\infty}\phi^{-1}(t_n)$, we obtain:
\begin{equation}\label{eqbeta1}
\lim_{n\to\infty}\dist_{\mathrm{eucli}}(p_{\R^k}\circ l_n(\overline{y}_{\infty}^{(n)},t_n),l_n^{\R}(t_n))=0.
\end{equation}
Finally, observe that equations \ref{eqf} and \ref{eqbeta1} contradict inequality \ref{eqpart2}, which concludes the proof.
\end{proof}

The continuity of the soul map is a consequence of the following proposition, which gives us property \ref{property1} as a corollary.

\begin{Proposition}\label{propsplitconv}The sequence $\{k_n^{S},l_n^S\}$ (resp. $\{k_n^{\R},l_n^{\R}\}$) realizes the convergence of $\{\overline{X}_n,\overline{\dist}_n,\overline{\m}_n\}$ (resp. $\{\R^k,\lvert\cdot\rvert,\mathcal{L}_k,0\}$) to $(\overline{X}_{\infty},\overline{\dist}_{\infty},\overline{\m}_{\infty})$ (resp. $(\R^k,\lvert\cdot\rvert,\mathcal{L}_k,0)$) in the mGH-topology (resp. pmGH-topology).
\end{Proposition}

\begin{proof}\proofpart{I}{$\{k_n^{\R},l_n^{\R}\}$ realizes the convergence of $\{\R^k,\lvert\cdot\rvert,\mathcal{L}_k,0\}$ to $(\R^k,\lvert\cdot\rvert,\mathcal{L}_k,0)$}

We are going to show that there exists a map $\epsilon^{\R}\colon\N\times\R_{\geq0}\to\R_{\geq0}$ such that for every $R>0$:
\begin{itemize}
\item[(i)] $\lim_{n\to\infty}\epsilon^{\R}(n,R)=0$,
\item[(ii)] for every $\lvert t\rvert \leq R$, we have $\dist_{\mathrm{eucli}}(k_n^{\R}\circ {l_n^{\R}}(t),t)\leq \epsilon^{\R}(n,R)$, and $\dist_{\mathrm{eucli}}(l_n^{\R}\circ {k_n^{\R}}(t),t)\leq \epsilon^{\R}(n,R)$ (when $n$ is large enough),
\item[(iii)] $\max\{\Dis({k_n^{\R}}_{\lvert B_{\R^k}(0,R)}),\Dis({l_n^{\R}}_{\lvert B_{\R^k}(0,R)})\}\leq\epsilon^{\R}(n,R)$ (when $n$ is large enough).
\end{itemize}
Then we will prove that $\{k_{n_*}^{\R}\mathcal{L}_k\}$ converges to $\mathcal{L}_k$ for the weak-$*$ topology.\\

Let $R>0$, and let $t\in\R^k$ such that $\lvert t\rvert \leq R$. Observe that $\dist_{\overline{X}_{\infty}\times\R^k}(k_n(\overline{*}_n,t),(\overline{*}_{\infty},k_n^{\R}(t)))\leq\alpha(n,R)$. Hence, thanks to Lemma \ref{prop7.2}, we get $\dist_{\overline{X}_{n}\times\R^k}(l_n\circ k_n(\overline{*}_n,t),l_n(\overline{*}_{\infty},k_n^{\R}(t)))\leq \alpha(n,R)+\epsilon_n\leq R$ (when $n$ is large enough). In addition, we have $\dist_{\overline{X}_{\infty}\times\R^k}((\overline{*}_{\infty},0),k_n(\overline{*}_n,t))\leq R+\epsilon_n$ (when $n$ is large enough). In particular, this implies $\lvert k_n^{\R}(t)\rvert\leq R+\alpha(n,R)+\epsilon_n\leq 2R$. Therefore, we get $\dist_{\overline{X}_{n}\times\R^k}(l_n(\overline{*}_{\infty},k_n^{\R}(t)),(\overline{*}_n,l_n^{\R}\circ k_n^{\R}(t)))\leq \beta(n,2R)$. In conclusion, when $n$ is large enough, we obtain $\dist_{\mathrm{eucli}}(t,l_n^{\R}\circ k_n^{\R}(t))\leq 2\epsilon_n+\alpha(n,R)+\beta(n,2R)$. The same strategy leads to $\dist_{\mathrm{eucli}}(t,k_n^{\R}\circ l_n^{\R}(t))\leq 2\epsilon_n+\beta(n,R)+\alpha(n,2R)$, for $n$ large enough. Therefore, we get point (ii) if we set $\epsilon^{\R}(n,R)\coloneqq2(\epsilon_n+\alpha(n,2R)+\beta(n,2R))$. Moreover, thanks to Lemma \ref{prop7.2}, we have $\lim_{n\to\infty}\epsilon^{\R}(n,R)=0$, for every $R>0$.\\
Now, given $t_1,t_2\in \R^k$ such that $\lvert t_1\lvert \leq R$ and $\lvert t_2\lvert \leq R$, we define:
$$
A\coloneqq\lvert \lvert \dist_{\mathrm{eucli}}(k_n^{\R}(t_1),k_n^{\R}(t_2))-\dist_{\mathrm{eucli}}(t_1,t_2) \rvert - \lvert \dist_{\overline{X}_{\infty}\times\R^k}(k_n(\overline{*}_n,t_1),k_n(\overline{*}_n,t_2))- \dist_{\mathrm{eucli}}(t_1,t_2)\rvert \rvert.
$$
Using $k_n^S(\overline{*}_{n})=\overline{*}_{\infty}$, we get:
\begin{align*}
A&\leq \lvert \dist_{\mathrm{eucli}}(k_n^{\R}(t_1),k_n^{\R}(t_2))- \dist_{\overline{X}_{\infty}\times\R^k}(k_n(\overline{*}_n,t_1),k_n(\overline{*}_n,t_2)) \rvert\\
&\leq \dist_{\overline{X}_{\infty}\times\R^k}(k_n(\overline{*}_n,t_1),(\overline{*}_{\infty},k_n^{\R}(t_1)))+ \dist_{\overline{X}_{\infty}\times\R^k}(k_n(\overline{*}_n,t_2),(\overline{*}_{\infty},k_n^{\R}(t_2)))\\
&\leq 2\alpha(n,R).
\end{align*}
Hence, for $n$ large enough, we have:
\begin{align*}
\lvert \dist_{\mathrm{eucli}}(k_n^{\R}(t_1),k_n^{\R}(t_2))-\dist_{\mathrm{eucli}}(t_1,t_2) \rvert&\leq A+\lvert \dist_{\overline{X}_{\infty}\times\R^k}(k_n(\overline{*}_n,t_1),k_n(\overline{*}_n,t_2))- \dist_{\mathrm{eucli}}(t_1,t_2)\rvert\\
&\leq 2\alpha(n,R)+\epsilon_n.
\end{align*}
In particular, this implies $\Dis({k_n^{\R}}_{\lvert B_{\R^k}(0,R)})\leq2\alpha(n,R)+\epsilon_n$. Moreover, $k_n$ and $l_n$ playing symmetric roles, we also have $\Dis({l_n^{\R}}_{\lvert B_{\R^k}(0,R)})\leq2\beta(n,R)+\epsilon_n$. We replace $\epsilon^{\R}(n,R)$ by $\epsilon^{\R}(n,R)+2\alpha(n,R)+2\beta(n,R)+\epsilon_n$. This concludes the proof of (i), (ii), and (iii) (thanks to Lemma \ref{prop7.2}).\\

Let us prove that $\{k_{n_*}^{\R}\mathcal{L}_k\}$ converges to $\mathcal{L}_k$ for the weak-$*$ topology. Here, the strategy will be the same as in the proof of Theorem \ref{thmA}. More precisely, the weak-$*$ topology on $\mathcal{M}_{\mathrm{loc}}(\R^k)$ is metrizable; therefore, it is equivalent to prove that every subsequence of $\{k_{n_*}^{\R}\mathcal{L}_k\}$ admits a subsequence converging to $\mathcal{L}_k$ in the weak-$*$ topology. Let us just prove that $\{k_{n_*}^{\R}\mathcal{L}_k\}$ admits a subsequence converging to $\mathcal{L}_k$ in the weak-$*$ topology (the proof for a subsequence being exactly the same). First of all, observe that thanks to Lemma \ref{lem7.1}, we can assume (passing to a subsequence if necessary) that $\{t\in\R^k\to k_n(\overline{*}_n,t)\in\overline{X}_{\infty}\times\R^k\}$ converges locally uniformly to $(\overline{y}_{\infty},\phi)$ for some $\overline{y}_{\infty}\in\overline{X}_{\infty}$ and $\phi\in\Or_k(\R)$. In particular, $\{k_n^{\R}\}$ converges locally uniformly to $\phi$. Then, notice that $\overline{y}_{\infty}=\lim_{n\to\infty}p_{\overline{X}_{\infty}}\circ k_n(\overline{*}_n,0)=\overline{*}_{\infty}$. Now let $R>0$, and let $f\in\Class_c(\R^k)$ be a continuous function such that $\Spt(f)\subset B_{\R^k}(0,R)$, and let us show that $\lim_{n\to\infty}\int_{\R^k}f\dist k_{n_*}^{\R}\mathcal{L}_k=\int_{\R^k}f\dist \mathcal{L}_k$.\\
First, observe that if $f\circ k_n^{\R}(t)\neq 0$, then $\lvert k_n^{\R}(t)\lvert\leq R$. Therefore, $\dist_{\mathrm{eucli}}(k_n^{\R}(t),k_n^{\R}(0))\leq R+\dist_{\mathrm{eucli}}(k_n^{\R}(0),0)$. In particular, we have: 
\[\dist_{\overline{X}_{\infty}\times\R^k}(k_n(\overline{*}_n,t),k_n(\overline{*}_n,0))\leq\tilde{R}_n\,,\]
where $\tilde{R}_n\coloneqq(\overline{D}^2+(R+\dist_{\mathrm{eucli}}(k_n^{\R}(0),0))^2)^{1/2}$, and $\overline{D}$ is defined in Notation \ref{not9.1}. Hence, if $f\circ k_n^{\R}(t)\neq 0$, then we have $\lvert t\rvert =\dist_{\overline{X}_{n}\times\R^k}((\overline{*}_n,t),(\overline{*}_n,0)) \leq 3\epsilon_n+\tilde{R}_n\leq 2\tilde{R}\coloneqq 2(\overline{D}^2+R^2)^{1/2}
$ (then $n$ is sufficiently large). Hence, whenever $n$ is large enough, we have:
\begin{align*}
\lvert\int_{\R^k}f\dist k_{n_*}^{\R}\mathcal{L}_k-\int_{\R^k}f\dist \phi_*\mathcal{L}_k\rvert&\leq \int_{B(2\tilde{R})}\lvert f\circ k_n^{\R}(t)-f\circ \phi(t)\rvert \dist\mathcal{L}_k(t)\\
&\leq \mathcal{L}_k(B(2\tilde{R}))\omega_n,
\end{align*}
where $\omega_n\coloneqq\sup_{\lvert x-y\rvert\leq\nu_n}\{\lvert f(x)-f(y)\rvert\}$ and $\nu_n\coloneqq\sup_{t\in B(2\tilde{R})}\{\lvert k_n^{\R}(t)-\phi(t)\rvert\}$. Observe that $k_n^{\R}$ converges locally uniformly to $\phi$; thus $\nu_n\to0$. In particular, since $f$ has compact support, $\omega_n\to0$; hence $\lim_{n\to\infty}\int_{\R^k}f\dist k_{n_*}^{\R}\mathcal{L}_k=\int_{\R^k}f\dist \phi_*\mathcal{L}_k$. In conclusion, passing to a subsequence if necessary, $\{k_{n_*}^{\R}\mathcal{L}_k\}$ converges in the weak-$*$ topology to $\phi_*\mathcal{L}_k=\mathcal{L}_k$ (using $\phi\in\Or_k(\R)$).

\proofpart{II}{$\{k_n^{S},l_n^{S}\}$ realizes the convergence of $\{\overline{X}_n,\overline{\dist}_n,\overline{\m}_n\}$ to $(\overline{X}_{\infty},\overline{\dist}_{\infty},\overline{\m}_{\infty})$}
Let $\overline{y}_n\in\overline{X}_n$, and observe that $\dist_{\overline{X}_{\infty}\times\R^k}(k_n(\overline{y}_n,0),(k_n^S(\overline{y}_n),0))\leq\alpha(n,0)$. Thanks to Lemma \ref{prop7.2}, we have $\dist_{\overline{X}_{n}\times\R^k}(l_n\circ k_n(\overline{y}_n,0),l_n(k_n^S(\overline{y}_n),0))\leq\alpha(n,0)+\epsilon_n$ (when $n$ is large enough). However, we have $\dist_{\overline{X}_{n}\times\R^k}(l_n(k_n^S(\overline{y}_n),0),(l_n^S\circ k_n^S(\overline{y}_n),0))\leq\beta(n,0)$. In conclusion, we have: 
\begin{equation}\label{thbpart21}
\overline{\dist}_n(\overline{y}_n,l_n^S\circ k_n^S(\overline{y}_n))=\dist_{\overline{X}_{n}\times\R^k}((\overline{y}_n,0),(l_n^S\circ k_n^S(\overline{y}_n),0))\leq \alpha(n,0)+\beta(n,0)+2\epsilon_n\eqqcolon\epsilon_n^S.
\end{equation}
Since $k_n$ and $l_n$ play symmetric roles, we also have:
\begin{equation}\label{thbpart22}
\overline{\dist}_{\infty}(\overline{y}_{\infty},k_n^S\circ l_n^S(\overline{y}_{\infty}))\leq \epsilon_n^S,
\end{equation}
for every $\overline{y}_{\infty}\in\overline{X}_{\infty}$. Observe that, thanks to Lemma \ref{prop7.2}, we have $\lim_{n\to\infty}\epsilon_n^S=0$.\\
Now, let $y_1,y_2\in\overline{X}_n$, and define:
$$
A\coloneqq\lvert \lvert \overline{\dist}_{\infty}(k_n^S(y_1),k_n^S(y_2))-\overline{\dist}_n(y_1,y_2) \rvert - \lvert \dist_{\overline{X}_{\infty}\times\R^k}(k_n(y_1,0),k_n(y_2,0))- \overline{\dist}_n(y_1,y_2) \rvert \rvert.
$$
Then, we have:
\begin{align*}
A&\leq \lvert \overline{\dist}_{\infty}(k_n^S(y_1),k_n^S(y_2))-\dist_{\overline{X}_{\infty}\times\R^k}(k_n(y_1,0),k_n(y_2,0)) \rvert\\
&\leq \dist_{\overline{X}_{\infty}\times\R^k}(k_n(y_1,0),(k_n^S(y_1),0))+\dist_{\overline{X}_{\infty}\times\R^k}(k_n(y_2,0),(k_n^S(y_2),0))\\
&\leq 2\alpha(n,0).
\end{align*}
Hence, we finally get (for $n$ large enough):
\begin{align*}
\lvert \overline{\dist}_{\infty}(k_n^S(y_1),k_n^S(y_2))-\overline{\dist}_n(y_1,y_2) \rvert&\leq A+\lvert \dist_{\overline{X}_{\infty}\times\R^k}(k_n(y_1,0),k_n(y_2,0))- \overline{\dist}_n(y_1,y_2) \rvert \rvert\\
&\leq 2\alpha(n,0)+\epsilon_n.
\end{align*}
In particular, we have $\Dis(k_n^S)\leq2\alpha(n,0)+\epsilon_n$. Then, $k_n$ and $l_n$ playing symmetric roles, we also have $\Dis(l_n^S)\leq2\beta(n,0)+\epsilon_n$. Finally, replacing $\epsilon_n^S$ by $\epsilon_n^S+\epsilon_n+2\max\{\alpha(n,0),\beta(n,0)\}$, and applying Lemma \ref{prop7.2}, we have $\epsilon_n^S\to0$; therefore, using inequalities \ref{thbpart21} and \ref{thbpart22}, we can conclude that $\{\overline{X}_n,\overline{\dist}_n\}$ converges to $(\overline{X}_{\infty},\overline{\dist}_{\infty})$ in the GH-topology.\\

Now let us prove that $\{k_{n_*}^S\overline{\m}_n\}$ converges to $\overline{\m}_{\infty}$ in the weak-$*$ topology. As we've seen in Part I, it is sufficient to show that $\{k_{n_*}^S\overline{\m}_n\}$ admits a subsequence converging to $\overline{\m}_{\infty}$ in the weak-$*$ topology (the weak-$*$ topology being metrizable).\\
First, let us show that $\{k_{n_*}^S\overline{\m}_n\}$ is precompact in the space $\mathcal{M}(\overline{X}_{\infty})$ of Radon measure on $\overline{X}_{\infty}$, which is implied by the uniform boundedness of the sequence $\{\overline{\m}_n(\overline{X}_n)\}$. Let us fix $r_0\in(0,\delta/2)$, where $\delta\coloneqq\inf_{n\in\N\cup\{\infty\}}\{\delta(X,\dist_n)\}$ ($\delta$ being positive thanks to Proposition \ref{prop6.2}). Then, using point (v) of Proposition \ref{prop3.4} and Theorem \ref{3.3}, observe that $\overline{\m}_n\otimes\mathcal{L}_k(B_{\overline{X}_n\times\R^k}((\overline{*}_n,0),r_0))=\m_n(B_{\dist_n}(*_n,r_0))\leq M$, where $M\coloneqq\sup_{n\in\N\cup\{\infty\}}\{\m_n(X)\}$ is finite. Moreover, notice that $B_{\overline{\dist}_n}(\overline{*}_n,r_0/\sqrt{2})\times B_{\R^k}(0,r_0/\sqrt{2})\subset B_{\overline{X}_n\times\R^k}((\overline{*}_n,0),r_0)$. Hence, $\overline{\m}_n(B_{\overline{\dist}_n}(\overline{*}_n,r_0/\sqrt{2}))\leq (\sqrt{2}/r_0)^kM/\omega_k$, where $\omega_k\coloneqq\mathcal{L}_k(B_{\R^k}(0,1))$. In particular, for every $n\in\N$, we can apply Bishop--Gromov inequality (see Theorem 6.2 in \cite{Bacher-Sturm_10}), and get $\overline{\m}_n(\overline{X}_n)\leq (\sqrt{2}/r_0)^N\overline{D}^{N-k}M/\omega_k\eqqcolon\overline{M}$, where $\overline{D}$ is defined in Noation \ref{not9.1}. In conclusion, $\{\overline{\m}_n(\overline{X}_n)\}$ is uniformly bounded; thus $\{k_{n_*}^S\overline{\m}_n\}$ is precompact in the weak-$*$ topology.\\
Now, passing to a subsequence if necessary, we can assume that $\{k_{n_*}^S\overline{\m}_n\}$ converges to some Radon measure $\overline{\m}$ on $\overline{X}_{\infty}$. We need to prove that $\overline{\m}=\overline{\m}_{\infty}$. Observe that it is equivalent to prove that $\overline{\m}\otimes\mathcal{L}_k=\overline{\m}_{\infty}\otimes\mathcal{L}_k$. First, thanks to the first part of the proof, $\{k_{n_*}^{\R}\mathcal{L}_k\}$ converges to $\mathcal{L}_k$ in the weak-$*$ topology; hence $\{(k_n^S,k_n^{\R})_*[\overline{\m}_n\otimes\mathcal{L}_k]\}$ converges to $\overline{\m}\otimes\mathcal{L}_k$ in the weak-$*$ topology. In addition, thanks to Theorem \ref{thmA}, $\{k_{n_*}[\overline{\m}_n\otimes\mathcal{L}_k]\}$ converges to $\overline{\m}_{\infty}\otimes\mathcal{L}_k$. Now, let $\phi\in\Class_c(\overline{X}_{\infty}\times\R^k)$ and $R>0$ such that $\Spt(\phi)\subset B_{\overline{X}_{\infty}\times\R^k}(R)$. Then, proceeding as in Part I of the proof, we obtain:
$$
\Spt(\phi\circ(k_n^S,k_n^{\R}) )\cup\Spt(\phi\circ k_n)\subset \overline{X}_{\infty}\times B_{\R^k}(0,2\tilde{R}),
$$
when $n$ is sufficiently large and where $\tilde{R}\coloneqq(\overline{D}^2+R^2)^{1/2}$. In particular, we have:
\begin{align*}
\lvert \int_{\overline{X}_{n}\times\R^k}\phi(k_n^S(\overline{x}),k_n^{\R}(t))- \phi(k_n(\overline{x},t))\dist\overline{\m}_n\otimes\mathcal{L}_k(\overline{x},t) \rvert&\leq (2\tilde{R})^k\omega_k\overline{M}\omega_{\phi}(\alpha(n,2\tilde{R})),
\end{align*}
where $\omega_{\phi}$ is the modulus of uniform continuity associated to $\phi$. Then, thanks to Lemma \ref{prop7.2}, we have $\lim_{n\to\infty}\omega_{\phi}(\alpha(n,2\tilde{R}))=0$. Thus, for every $\phi\in\Class_c(\overline{X}_{\infty}\times\R^k)$, we have:
$$
\lim_{n\to\infty}\int_{\overline{X}_{n}\times\R^k}\phi\dist(k_n^S,k_n^{\R})_*[\overline{\m}_n\otimes\mathcal{L}_k]=\lim_{n\to\infty}\int_{\overline{X}_{n}\times\R^k}\phi\dist k_{n_*}[\overline{\m}_n\otimes\mathcal{L}_k].
$$
In particular, $\{(k_n^S,k_n^{\R})_*[\overline{\m}_n\otimes\mathcal{L}_k]\}$ and $\{k_{n_*}[\overline{\m}_n\otimes\mathcal{L}_k]\}$ have the same limit, i.e. $\overline{\m}\otimes\mathcal{L}_k=\overline{\m}_{\infty}\otimes\mathcal{L}_k$. This concludes the proof.
\end{proof}

Inspired by the proof of Theorem 5.4 in \cite{Tuschmann-Wiemeler_17}, we introduce the following "shrunk" metrics.

\begin{Definition}\label{7.1}Given $n\in\N\cup\{\infty\}$, and $m\in\N$, let $\tilde{\dist}_{n,m}\coloneqq\phi_n^{*}(2^{-m}\overline{\dist}_n\times\dist_{\mathrm{eucli}})$, and let $(X,\dist_{n,m},\m_n)$ be the push-forward of $(\tilde{X},\tilde{\dist}_{m,n},\tilde{m}_n)$ (see Proposition \ref{propX}).
\end{Definition}

\begin{Remark}\label{rem7.1}Given $n\in\N\cup\{\infty\}$, $m\in\N$, and $\tilde{x},\tilde{y}\in \tilde{X}$, we have $\tilde{\dist}_{n,m}(\tilde{x},\tilde{y})\leq \tilde{\dist}_n(\tilde{x},\tilde{y})$. Indeed, defining $(\overline{x},t_x)\coloneqq\phi_n(x)$ and $(\overline{y},t_y)\coloneqq\phi_n(y)$, we have:
$$
\tilde{\dist}_{n,m}(\tilde{x},\tilde{y})^2=2^{-2m}\overline{\dist}_n^2(\overline{x},\overline{y})+\dist_{\mathrm{eucli}}^2(t_x,t_y)\leq \overline{\dist}_n^2(\overline{x},\overline{y})+\dist_{\mathrm{eucli}}^2(t_x,t_y)=\tilde{\dist}_{n}(\tilde{x},\tilde{y})^2.
$$
In particular, this implies $\dist_{n,m}\leq\dist_n$.
\end{Remark}

The following lemma shows that the "shrunk" metrics associated to the sequence $\{X,\dist_n,\m_n\}$ are close to the corresponding Albanese varieties.

\begin{Lemma}\label{prop7.5}We have $\GH([X,\dist_{n,m}],\mathcal{A}([X,\dist_n,\m_n]))\leq2^{-m+1}\overline{D}$, for every $m\in\N$, and $n\in\N\cup\{\infty\}$ (where $\overline{D}$ is defined in Notation \ref{not9.1}).
\end{Lemma}

\begin{proof}First of all, observe that there exists a continuous map $a\colon X\to (\R^k/\Gamma(\phi_n),\dist_{\Gamma(\phi_n)})$ such that $a\circ p=q\circ p_{\R^k}\circ\phi_n$, where $q\colon\R^k\to\R^k/\Gamma(\phi_n)$ is the quotient map. Notice that, $q\circ p_{\R^k}\circ\phi_n$ is surjective; hence, $a$ is also surjective. Now, let $y,z\in X$ and let us show that $\lvert \dist_{n,m}(y,z)-\dist_{\Gamma(\phi_n)}(a(y),a(z))\rvert \leq 2^{-m}\overline{D}$.\\
First, thanks to Proposition \ref{propX}, there exists $\tilde{y}\in p^{-1}(y)$ and $\tilde{z}\in p^{-1}(z)$ such that $\dist_{n,m}(y,z)=\tilde{\dist}_{m,n}(\tilde{y},\tilde{z})$. Let $(\overline{y},t_y)\coloneqq\phi_n(\tilde{y})$ and $(\overline{z},t_z)\coloneqq\phi_n(\tilde{z})$. Observe that $a(y)=q(t_y)$ and $a(z)=q(t_z)$, and $\dist^2_{n,m}(y,z)=2^{-2m}\overline{\dist}^2_n(\overline{y},\overline{z})+\dist_{\mathrm{eucli}}^2(t_y,t_z).$ Now, note that, by definition of $\dist_{\Gamma(\phi_n)}$, we have $\dist_{\Gamma(\phi_n)}(a(y),a(z))\leq \dist_{\mathrm{eucli}}(t_y,t_z)$; in particular:
\begin{equation}\label{eqm}
0\leq \dist_{n,m}(y,z)-\dist_{\Gamma(\phi_n)}(a(y),a(z)).
\end{equation}
Notice that, by definition of $\dist_{\Gamma(\phi_n)}$, there exists $\eta\in\overline{\pi}_1(X)$ such that $\dist_{\Gamma(\phi_n)}(a(y),a(z))=\dist_{\mathrm{eucli}}(t_y, t_2)$, where $t_2\coloneqq\rho_{\R}^{\phi_n}(\eta)\cdot t_z$. Then, observe that $\tilde{\dist}_{n,m}(\tilde{y},\tilde{z})=\dist_{n,m}(y,z)\leq \tilde{\dist}_{n,m}(\tilde{y},\eta\cdot\tilde{z})$; hence:
\begin{align*}
 \dist_{n,m}(y,z)-\dist_{\Gamma(\phi_n)}(a(y),a(z))&\leq (2^{-2m}\overline{\dist}^2_n(\overline{y},\rho_S^{\phi_n}(\eta)\cdot \overline{z})+\dist_{\mathrm{eucli}}(t_y,t_2))^{1/2}-\dist_{\mathrm{eucli}}(t_y,t_2)\\
 &\leq 2^{-m}\overline{\dist}_n(\overline{y},\rho_S^{\phi_n}(\eta)\cdot \overline{z})+\dist_{\mathrm{eucli}}(t_y,t_2)-\dist_{\mathrm{eucli}}(t_y,t_2)\\
 &\leq 2^{-m}\overline{D}.
\end{align*}
In particular, thanks to inequality \ref{eqm}, we have $\lvert \dist_{n,m}(y,z)-\dist_{\Gamma(\phi_n)}(a(y),a(z))\rvert \leq 2^{-m}\overline{D}$. Therefore, recalling that $a$ is surjective, and using Corollary 7.3.28 of \cite{Burago-Ivanov_01}, we get $\GH([X,\dist_{n,m}],\mathcal{A}([X,\dist_n,\m_n]))\leq2^{-m+1}\overline{D}$.
\end{proof}

To prove property \ref{property2}, we will need to obtain a convergence result on the following quantities.

\begin{Notation}\label{not9.3}Given $R>0$, $n,m\in\N$, we denote:
\begin{itemize}
\item[(i)]$\epsilon(n,m,R)\coloneqq\sup\{\lvert \tilde{\dist}_{\infty,m}(\tilde{f}_n(\tilde{y}_1),\tilde{f}_n(\tilde{y}_2))-\tilde{\dist}_{n,m}(\tilde{y}_1,\tilde{y}_2)\rvert\}$, the supremum being taken over $\tilde{y}_i\in \tilde{B}_n(R)$,
\item[(ii)]$\epsilon'(n,m,R)\coloneqq\sup\{\lvert \tilde{\dist}_{n,m}(\tilde{g}_n(\tilde{y}_1),\tilde{g}_n(\tilde{y}_2))-\tilde{\dist}_{\infty,m}(\tilde{y}_1,\tilde{y}_2)\rvert\}$, the supremum being taken over ${\tilde{y}_i\in \tilde{B}_{\infty}(R)}$.
\end{itemize}
\end{Notation}

\begin{Lemma}\label{lem9.5}For every $R>0$, $\lim_{n,m\to\infty}\epsilon(n,m,R)=\lim_{n,m\to\infty}\epsilon'(n,m,R)=0$.
\end{Lemma}

\begin{proof}Let us only prove that $\lim_{n,m\to\infty}\epsilon(n,m,R)=0$, the proof for $\epsilon'(n,m,R)$ being exactly the same. Let $\tilde{y}_i\in\tilde{B}_n(R)$, $i\in\{1,2\}$. For $i\in\{1,2\}$, we denote $(\overline{y}_i,t_i)\coloneqq \phi_n(\tilde{y}_i)$, $\tilde{y}_i^{\infty}\coloneqq \tilde{f}_n(\tilde{y}_i)$, $(\overline{y}^{\infty}_i,t^{\infty}_i)\coloneqq \phi_{\infty}(\tilde{y}^{\infty}_i)$, and $A\coloneqq\lvert \tilde{\dist}_{\infty,m}(\tilde{f}_n(\tilde{y}_1),\tilde{f}_n(\tilde{y}_2))-\tilde{\dist}_{n,m}(\tilde{y}_1,\tilde{y}_2)\rvert$. Using the fact that, for every $x,y\in\R_{\geq0}$, we have $\lvert \sqrt{x}-\sqrt{y}\rvert\leq\sqrt{\lvert x-y\rvert}$, we get:
\begin{align*}
A&\leq (\lvert 2^{-2m}(\overline{\dist}_{\infty}^2(\overline{y}^{\infty}_1,
,\overline{y}^{\infty}_2) - \overline{\dist}_{n}^2(\overline{y}_1,
,\overline{y}_2)) +(\dist_{\mathrm{eucli}}^2(t^{\infty}_1,t^{\infty}_2)-\dist_{\mathrm{eucli}}^2(t_1,t_2))\rvert)^{1/2}\\
&\leq 2^{-m}(\lvert (\overline{\dist}_{\infty}^2(\overline{y}^{\infty}_1,
,\overline{y}^{\infty}_2) - \overline{\dist}_{n}^2(\overline{y}_1,
,\overline{y}_2)) \rvert)^{1/2}+ (\lvert \dist_{\mathrm{eucli}}^2(t^{\infty}_1,t^{\infty}_2)-\dist_{\mathrm{eucli}}^2(t_1,t_2)\rvert)^{1/2}.
\end{align*}
However, note that $\lvert (\overline{\dist}_{\infty}^2(\overline{y}^{\infty}_1,\overline{y}^{\infty}_2) -\overline{\dist}_{n}^2(\overline{y}_1,
\overline{y}_2)) \rvert \leq 2\overline{D}^2$, 
where $\overline{D}$ is introduced in Notation \ref{not9.1}. Then, observe that $\dist_{\mathrm{eucli}}(t^{\infty}_1,t^{\infty}_2)\leq \tilde{\dist}_{\infty}(\tilde{y}_1^{\infty},\tilde{y}_2^{\infty})\leq \tilde{\dist}_n(\tilde{y}_1,\tilde{y}_2)+\epsilon_n\leq 2R+\epsilon_n$, when $n$ is large enough. Therefore, $\dist_{\mathrm{eucli}}(t^{\infty}_1,t^{\infty}_2)+\dist_{\mathrm{eucli}}(t_1,t_2)\leq4R+\epsilon_n$. Then, using $\lvert t_i\rvert \leq R$, and denoting $B\coloneqq\lvert \dist_{\mathrm{eucli}}(t^{\infty}_1,t^{\infty}_2)-\dist_{\mathrm{eucli}}(t_1,t_2)\rvert$, we have:
\begin{align*}
B&\leq \lvert \dist_{\mathrm{eucli}}(t^{\infty}_1,t^{\infty}_2)-\dist_{\mathrm{eucli}}(k_n^{\R}(t_1),k_n^{\R}(t_2))\rvert+\lvert \dist_{\mathrm{eucli}}(k_n^{\R}(t_1),k_n^{\R}(t_2))-\dist_{\mathrm{eucli}}(t_1,t_2)\rvert\\
&\leq \Dis(k^{\R}_{n_{\lvert B_{\R^k}(0,R)}})+\dist_{\mathrm{eucli}}(k_n^{\R}(t_1),t^{\infty}_1)+\dist_{\mathrm{eucli}}(k_n^{\R}(t_2),t^{\infty}_2)\\
&\leq \Dis(k^{\R}_{n_{\lvert B_{\R^k}(0,R)}})+2\alpha(n,R).
\end{align*}
In conclusion, we obtain:
$$
A\leq 2^{-m+1/2}\overline{D}+(4R+\epsilon_n)^{1/2}(\Dis(k^{\R}_{n_{\lvert B_{\R^k}(0,R)}})+2\alpha(n,R))^{1/2}\eqqcolon \tilde{\epsilon}(n,m,R).
$$
Therefore, passing to the supremum as $\tilde{y}_i\in\tilde{B}_n(R)$ ($i\in\{1,2\}$), we obtain $\epsilon(n,m,R)\leq\tilde{\epsilon}(n,m,R)$. Thanks to Lemma \ref{prop7.2} and Proposition \ref{propsplitconv}, we have $\lim_{n,m\to\infty}\tilde{\epsilon}(n,m,R)=0$.
\\Therefore, $\lim_{n,m\to\infty}{\epsilon}(n,m,R)=0$, which concludes the proof.
\end{proof}

We conclude this section with the following proposition, which states the continuity of the Albanese map by proving property \ref{property2}.

\begin{Proposition}The sequence $\{[\R^k/\Gamma(\phi_n),\dist_{\Gamma(\phi_n)}]=\mathcal{A}(X,\dist_n,\m_n)\}$ converges in the GH-topology to $[\R^k/\Gamma(\phi_{\infty}),\dist_{\Gamma(\phi_{\infty})}]=\mathcal{A}(X,\dist_{\infty},\m_{\infty})$.
\end{Proposition}

\begin{proof}Observe that for every $n\in\N\cup\{\infty\}$, we have $\GH([X,\dist_{n,n}],\mathcal{A}(X,\dist_n,\m_n))\leq2^{-n+1}\overline{D}$, thanks to Lemma \ref{prop7.5}
(where $\dist_{n,n}$ is defined in Definition \ref{7.1} and $\overline{D}$ is introduced in Notation \ref{not9.1}). In particular, using the triangle inequality for the Gromov--Hausdorff distance $\GH$, we obtain:
$$
\GH(\mathcal{A}(X,\dist_n,\m_n),\mathcal{A}(X,\dist_{\infty},\m_{\infty}))\leq2^{-n+2}\overline{D}+\GH([X,\dist_{n,n}],[X,\dist_{\infty,n}]).
$$
Therefore, to conclude, it is sufficient to prove that:
$$
\lim_{n,m\to\infty}\GH([X,\dist_{n,m}],[X,\dist_{\infty,m}])=0,
$$
which is what we are going to prove.\\

Let $n,m\in\N$, and $y_1,y_2\in X$. There exists $\tilde{y}_1\in p^{-1}(y_1)$, and $\tilde{y}_2\in p^{-1}(y_2)$, such that $\tilde{\dist}_n(\tilde{*}_n,\tilde{y}_1)=\dist_n(*_n,y_1)$, and $\tilde{\dist}_{n,m}(\tilde{y}_1,\tilde{y}_2)=\dist_{n,m}(y_1,y_2)$. Then, for $i\in\{1,2\}$, we denote $(\overline{y}_i,t_i)\coloneqq\phi_n(\tilde{y}_i)$. Observe that $\tilde{\dist}_n(\tilde{y}_1,\tilde{y}_2)=(\overline{\dist}_n^2(\overline{y}_1,\overline{y}_2)+\dist_{\mathrm{eucli}}^2(t_1,t_2))^{1/2}$, where $\overline{\dist}_n(\overline{y}_1,\overline{y}_2)\leq \overline{D}$, and $\dist_{\mathrm{eucli}}(t_1,t_2)\leq\dist_{n,m}(y_1,y_2)\leq \dist_n(y_1,y_2)\leq D$ (using Remark \ref{rem7.1}). Therefore, we get:
\begin{equation}\label{eqo}
\tilde{\dist}_n(\tilde{*}_n,\tilde{y}_2)\leq D+(\overline{D}^2+D^2)^{1/2}\eqqcolon \tilde{D}.
\end{equation}
Now, using $\dist_{n,m}(y_1,y_2)=\tilde{\dist}_{n,m}(\tilde{y}_1,\tilde{y}_2)$, and $\tilde{f}_n(\tilde{y}_i)\in p^{-1}(f_n(y_i))$, we have:
\begin{align*}
\dist_{\infty,m}(f_n(y_1),f_n(y_2))-\dist_{n,m}(y_1,y_2)&\leq \tilde{\dist}_{\infty,m}(\tilde{f}_n(\tilde{y}_1),\tilde{f}_n(\tilde{y}_2))-\tilde{\dist}_{n,m}(\tilde{y}_1,\tilde{y}_2)\\
&\leq \epsilon(n,m,\tilde{D}),
\end{align*}
where $\epsilon(n,m,R)$ is introduced in Notation \ref{not9.3}. Since $\tilde{f}_n$ and $\tilde{g}_n$ play symmetric roles, we also have:
$$
\forall y_1',y_2'\in X,\dist_{n,m}(g_n(y_1'),g_n(y_2'))-\dist_{\infty,m}(y_1,y_2)\leq \epsilon'(n,m,\tilde{D}),
$$
where $\epsilon'(n,m,\tilde{D})$ is also introduced in Notation \ref{not9.3}. In particular, this implies:
$$ 
\forall y_1,y_2\in X,\dist_{n,m}(g_n\circ f_n(y_1),g_n\circ f_n(y_2))-\dist_{\infty,m}(f_n(y_1),f_n(y_2))\leq \epsilon'(n,m,\tilde{D}).
$$
Hence, defining $A\coloneqq\dist_{n,m}(y_1,y_2)-\dist_{\infty,m}(f_n(y_1),f_n(y_2))$, we have:
\begin{align*}
A&\leq \dist_{n,m}(y_1,y_2)-\dist_{n,m}(g_n\circ f_n(y_1),g_n\circ f_n(y_2))+\epsilon'(n,m,\tilde{D})\\
&\leq\dist_{n,m}(g_n\circ f_n(y_1),y_1)+\dist_{n,m}(g_n\circ f_n(y_2),y_2)+\epsilon'(n,m,\tilde{D})\\
&\leq \dist_{n}(g_n\circ f_n(y_1),y_1)+\dist_{n}(g_n\circ f_n(y_2),y_2)+\epsilon'(n,m,\tilde{D})\\
&\leq 2\epsilon_n+\epsilon'(n,m,\tilde{D}).
\end{align*}
In conclusion, we have:
\begin{equation}\label{eqp}
\lvert \dist_{\infty,m}(f_n(y_1),f_n(y_2))-\dist_{n,m}(y_1,y_2)\lvert \leq2\epsilon_n+\epsilon(n,m,\tilde{D})+\epsilon'(n,m,\tilde{D}).
\end{equation}
Moreover, since $\dist_{\infty,m}\leq\dist_{\infty}$, and since $f_n$ is an $\epsilon_n$-isometry from $(X,\dist_n)$ onto $(X,\dist_{\infty})$, we have:
\begin{equation}\label{eqq}
\forall x\in X, \exists y\in X, \dist_{\infty,m}(x,f_n(y))\leq\epsilon_n.
\end{equation}
Hence, thanks to inequalities \ref{eqp} and \ref{eqq}, $f_n$ is a $2\epsilon_n+\epsilon(n,m,\tilde{D})+\epsilon'(n,m,\tilde{D})$-isometry from $(X,\dist_{n,m})$ to $(X,\dist_{\infty,m})$. Therefore, using Corollary 7.3.28 of \cite{Burago-Ivanov_01}, we have:
$$
\GH([X,\dist_{n,m}],[X,\dist_{\infty,m}])\leq2(2\epsilon_n+\epsilon(n,m,\tilde{D})+\epsilon'(n,m,\tilde{D})).
$$
However, thanks to Lemma \ref{lem9.5}, we have $\lim_{n,m\to\infty}2\epsilon_n+\epsilon(n,m,\tilde{D})+\epsilon'(n,m,\tilde{D})=0$, which concludes the proof.

\end{proof}

\subsection{Proof of Theorem \ref{thmC}}\label{proofC}

The proof of Theorem \ref{thmC} is inspired by the proof of Theorem 1.1 in \cite{Tuschmann-Wiemeler_17} and uses some of the computations realized in \cite{Perez_20}.\\

First of all, using Theorem \ref{thmB}, we are going to prove the following result.

\begin{Proposition}\label{prop: retraction}Let $N\in[1,\infty)$, let $X$ be a compact topological space that admits an $\RCD(0,N)$-structure such that $\overline{\pi}_1(X)=0$ (see Theorem \ref{th3.2} for the definition of $\overline{\pi}_1(X)$), and let $\Gamma$ be a Bieberbach subgroup of $\R^k$ ($k\geq2$). Then, the moduli space $\M_{0,N+k}(X\times \R^k/\Gamma)$ retracts onto $\mathscr{M}_{\mathrm{flat}}(\R^k/\Gamma)$.
\end{Proposition}

\begin{proof}Let us describe the crystallographic class $\Gamma(X\times \R^k/\Gamma)$ (introduced in Proposition \ref{prop5.1}). First, observe that, since $\overline{\pi}_1(X)=0$, then the universal cover of $X\times \R^k/\Gamma$ is $X\times \R^k$, and the covering map is just $\id_{X}\times q$, where $q\colon \R^k\to\R^k/\Gamma$ is the usual quotient map. Now, let $g$ be the flat Riemannian metric on $\R^k/\Gamma$ such that $q$ is a local isometry, and fix an $\RCD(0,N)$-structure $(X,\dist_0,\m_0)$ on $X$. Observe that $(X,\dist_0,\m_0)\times(\R^k/\Gamma,\dist_g,\m_g)$ is an $\RCD(0,N+k)$-structure on $X\times\R^k/\Gamma$, where $d_g$ and $\m_{g}$ are respectively the Riemannian distance and measure associated to $g$. Moreover, the lifted $\RCD(0,N+k)$-structure on $X\times \R^k$ is equal to $(X,\dist_0,\m_0)\times(\R^k,\dist_{\mathrm{eucli}},\mathcal{L}_k)$. In particular, the identity map $\id_{X\times\R^k}$ is a splitting of $(X,\dist_0,\m_0)\times\R^k$. Moreover, since $\overline{\pi}_1(X\times\R^k/\Gamma)$ acts trivially on $X$, we have $\Gamma(\id_{X\times\R^k})=\Gamma$. Hence, the crystallographic class $\Gamma(X\times\R^k/\Gamma)$ is equal to the set of crystallographic subgroup of $\Iso(\R^k)$ that are isomorphic to $\Gamma$. This implies, thanks to Remark \ref{rem5.1}, that $\mathscr{M}_{\text{flat}}(A(X\times\R^k/\Gamma))$ is isometric to $\mathscr{M}_{\mathrm{flat}}(\R^k/\Gamma)$.\\

Now, thanks to Theorem \ref{thmB}, the Albanese map associated to $X\times\R^k/\Gamma$ is continuous from $\M_{0,N+k}(X\times\R^k/\Gamma)$ onto $\mathscr{M}_{\text{flat}}(A(X\times\R^k/\Gamma))$. Hence, it gives rise to a continuous surjective map $\phi$ from $\M_{0,N+k}(X\times\R^k/\Gamma)$ onto $\mathscr{M}_{\mathrm{flat}}(\R^k/\Gamma)$. Given $[\R^k/\Gamma,\dist]\in\mathscr{M}_{\mathrm{flat}}(\R^k/\Gamma)$, we define:
$$
s([\R^k/\Gamma,\dist])\coloneqq[(X,\dist_0,\m_0)\times(\R^k/\Gamma,\dist,\mathcal{H}_{\dist})]\in\M_{0,N+k}(X\times\R^k/\Gamma),
$$
where $\mathcal{H}_{\dist}$ is the Hausdorff measure associated to $(\R^k/\Gamma,\dist)$. Observe that $s$ is a section of $\phi$; therefore, we only have to show that $s$ is continuous in order to conclude the proof.\\

Let us show that $s\colon\mathscr{M}_{\mathrm{flat}}(\R^k/\Gamma)\to\M_{0,N+k}(X\times\R^k/\Gamma)$ is continuous. To do so, let $\{(\R^k/\Gamma,\dist_n)\}$ converge in the Gromov--Hausdorff sense to $(\R^k/\Gamma,\dist_{\infty})$, where, for every $n\in\N\cup\{\infty\}$, $\dist_{n}$ is a flat metric on $\R^k/\Gamma$. Let us prove that $\{(X,\dist_0,\m_0)\times(\R^k/\Gamma,\dist_n,\mathcal{H}_{\dist_n})\}$ converges in the measured Gromov--Hausdorff sense to $(X,\dist_0,\m_0)\times(\R^k/\Gamma,\dist_{\infty},\mathcal{H}_{\dist_{\infty}})$. Observe that it is sufficient to prove that $\{(\R^k/\Gamma,\dist_{n},\mathcal{H}_{\dist_{n}})\}$ converges in the mGH sense to $(\R^k/\Gamma,\dist_{\infty},\mathcal{H}_{\dist_{\infty}})$. However, since $\dist_{\infty}$ is a flat metric on $\R^k/\Gamma$, the Hausdorff dimension of $(\R^k/\Gamma,\dist_{\infty})$ is equal to $k$. In particular, by Theorem 1.2 of \cite{DePhilippis-Gigli_18}, $\{(\R^k/\Gamma,\dist_{n},\mathcal{H}_{\dist_{n}})\}$ converges in the mGH sense to $(\R^k/\Gamma,\dist_{\infty},\mathcal{H}_{\dist_{\infty}})$. In conclusion $s$ is continuous.
\end{proof}

Proposition \ref{prop: retraction} implies that the homotopy groups of $\mathscr{M}_{\mathrm{flat}}(\R^k/\Gamma)$ inject in those of $\M_{0,N+k}(X\times\R^k/\Gamma)$. Therefore, the topology of $\M_{0,N+k}(X\times\R^k/\Gamma)$ is, in a way, at least as complicated as the topology of $\mathscr{M}_{\mathrm{flat}}(\R^k/\Gamma)$. Thankfully, informations on moduli spaces of flat metrics have been derived in \cite{Tuschmann-Wiemeler_17} (in the case of the torus $T^k$, with $k\geq4$ and $k\neq8,9,10$) and in \cite{Perez_20} (in the case of $3$ and $4$-dimensional closed flat Riemannian manifolds). We are now able to prove Theorem \ref{thmC}.

\begin{proof}[Proof of Theorem \ref{thmC}]Observe that thanks to Theorem 3.4.3 of \cite{Perez_20} and Proposition 5.5 of \cite{Tuschmann-Wiemeler_17} the moduli space $\mathscr{M}_{\mathrm{flat}}(N)$ has non-trivial higher rational homotopy groups. Therefore, Proposition \ref{prop: retraction} concludes the proof.
\end{proof}

Let us now prove Corollary \ref{corB}.

\begin{proof}[Proof of Corollary \ref{corB}]First of all, observe that thanks to Theorem 3.4.3 of \cite{Perez_20}, the moduli space of flat metrics on $X_3\coloneqq\mathbb{S}^1\times\mathbb{K}^2$ is homotopy equivalent to a circle (where $\mathbb{K}^2$ is the Klein bottle). Then, let us define $X_4\coloneqq [0,1]\times X_3$, and $X_N\coloneqq \mathbb{S}^{N-3}\times X_3$ ($N\geq5$). Thanks to Proposition \ref{prop: retraction}, for every $N\geq3$, $\M_{0,N}(X_N)$ retracts onto $\mathscr{M}_{\mathrm{flat}}(X_3)$. In particular, for every $N\geq3$, $\M_{0,N}(X_N)$ has non trivial fundamental group.\\
To conclude the proof, we apply the same idea, using the fact that $\pi_3(\mathbb{T}^4)\otimes\Q\simeq\Q$, and $\pi_5(\mathbb{T}^5)\otimes\Q\simeq\Q$ (see Proposition 5.5 of \cite{Tuschmann-Wiemeler_17}).
\end{proof}

\section*{Appendix}\label{Appendix}

Before proving Proposition \ref{prop:equivtopo}, let us point out that the following technicality on the Prokhorov distance.

\begin{remAppendix}\label{rem:prokhorov}There are various notions of distances between restricted measures. Indeed, given a pointed complete separable metric space $(Y,\dist,\ast)$ endowed with two boundedly finite measures $\m_1$ and $\m_2$, and $R>0$, we can define:
$$
\dist_{\mathcal{P}}^R(\m_1,\m_2)\coloneqq\inf\Bigg\{\epsilon>0,\forall A\subset \overline{B}_R, A\ \text{closed implies}\ \Bigg\lvert\begin{tabular}[c]{l}$\m_1(A)\leq\m_2(A^{\epsilon})+\epsilon$\\ $\m_2(A)\leq\m_1(A^{\epsilon})+\epsilon$\end{tabular}\Bigg\},
$$
or:
$$
\dist_{\mathcal{P}}(\m_1^{R},\m_2^{R})\coloneqq\inf\Bigg\{\epsilon>0,\forall A\subset X, A\ \text{closed implies}\ \Bigg\lvert\begin{tabular}[c]{l}$\m_1^R(A)\leq\m_2^R(A^{\epsilon})+\epsilon$\\ $\m_2^R(A)\leq\m_1^R(A^{\epsilon})+\epsilon$\end{tabular}\Bigg\},
$$
where $\m_1^R$ (resp. $\m_2^R$) is the restriction of $\m_1$ (resp. $\m_2$) to $\overline{B}_R$. Let us see how these two notions differ.\\
First of all, we easily obtain:
\begin{align}\label{equiv1}\dist_{\mathcal{P}}^R(\m_1,\m_2)\leq\dist_{\mathcal{P}}(\m_1^{R},\m_2^{R}).
\end{align}\\
Then, if $\dist_{\mathcal{P}}^R(\m_1,\m_2)\leq\epsilon$, we have:
\begin{align}\label{equiv2}
\dist_{\mathcal{P}}(\m_1^{R},\m_2^{R})\leq \epsilon+\m_1(\overline{B}_{R+\epsilon}\backslash\overline{B}_{R})+\m_2(\overline{B}_{R+\epsilon}\backslash\overline{B}_{R}).\end{align}
Therefore, both notions lead to the same convergence. Indeed, given a boundedly finite measure $\m_{\infty}$ and a sequence $\{\m_k\}$ of boundedly finite measures, we have the following equivalence thanks to inequalities \ref{equiv1} and \ref{equiv2}:
$$
\m_k \overset{\ast}{\rightharpoonup}\m_{\infty}\iff \forall R>0,\dist_{\mathcal{P}}(\m_k^{R},\m_{\infty}^{R})\to0\iff \forall R>0,\dist_{\mathcal{P}}^R(\m_k,\m_{\infty})\to0.
$$
In Definition \ref{def:equivtopo}, we are using $\dist_{\mathcal{P}}^R(\m_1,\m_2)$ instead of $\dist_{\mathcal{P}}(\m_1^{R},\m_2^{R})$ because the first one is a non-decreasing function of $R$ whereas the other one is not. Doing so, we are losing the triangle inequality. In fact, given three boundedly finite measure $\m_i$ ($i\in\{1,2,3\}$) such that $\dist_{\mathcal{P}}^R(\m_1,\m_2)\leq\delta$ and $\dist_{\mathcal{P}}^R(\m_2,\m_3)\leq\eta$, we only have:
\begin{align}\label{equiv3}
\dist_{\mathcal{P}}^{R-(\delta+\eta)}(\m_1,\m_3)\leq\delta+\eta.
\end{align}
This will be sufficient to investigate the properties of $\D_{\mathrm{p}}^{\mathrm{eq}}$.
\end{remAppendix}

\begin{proof}[Proof of Proposition \ref{prop:equivtopo}]\proofpart{I}{Properties of $\Dpeq$}
Observe that $\Dpeq$ is symmetric, nonnegative and invariant under equivariant isomorphisms. However, $\Dpeq$ does not satisfy the triangle inequality a priori. Nevertheless, it will be sufficient for our purposes to show that $\Dpeq$ satisfy the following modified triangle inequality:
\begin{align}\label{trig}
\Dpeq(\tilde{\mathcal{X}}_1,\tilde{\mathcal{X}}_3) &\leq 4(\Dpeq(\tilde{\mathcal{X}}_1,\tilde{\mathcal{X}}_2)+\Dpeq(\tilde{\mathcal{X}}_2,\tilde{\mathcal{X}}_3))
\end{align}
for any three equivariant pointed $\RCD(0,N)$-structures $\tilde{\mathcal{X}}_i=(\tilde{X},\tilde{\dist}_i,\tilde{\m}_i,\tilde{\ast}_i)\in\kR_{0,N}^{\mathrm{p,eq}}(\tilde{X})$ ($i\in\{1,2,3\}$).\\
Observe that the inequality is trivially true whenever $\Dpeq(\tilde{\mathcal{X}}_1,\tilde{\mathcal{X}}_2)$ or $\Dpeq(\tilde{\mathcal{X}}_2,\tilde{\mathcal{X}}_3)$ is equal to $1/24$.\\
Now, assume that for $i\in\{1,2\}$, we have $\Dpeq(\tilde{\mathcal{X}}_{i},\tilde{\mathcal{X}}_{i+1})\leq\epsilon_{i,i+1}$, where $\epsilon_{i,i+1}\in(0,1/24)$, and let $(f_{i,i+1},g_{i,i+1},\phi_{i,i+1})$ be an associated equivariant pointed $\epsilon_{i,i+1}$-isometry. We define $f\coloneqq f_{23}\circ f_{12}$, $g\coloneqq g_{12}\circ g_{23}$, and $\phi\coloneqq\phi_{23}\circ\phi_{12}$. We want to show that $(f,g,\phi)$ is an equivariant pointed $4(\epsilon_{12}+\epsilon_{23})$-isometry between $\tilde{\mathcal{X}}_1$ and $\tilde{\mathcal{X}}_3$.\\
First, point (i) and point (ii) of Definition \ref{def:equivtopo} are trivially satisfied.\\
Then, let $x,y\in\tilde{X}$ such that $\dist_{1}(x,y)\leq[4(\epsilon_{12}+\epsilon_{23})]^{-1}$. Notice that $[4(\epsilon_{12}+\epsilon_{23})]^{-1}\leq\epsilon_{12}^{-1}$. Therefore $\tilde{\dist}_2(f_{12}(x),f_{12}(y))\leq\epsilon_{12}+[4(\epsilon_{12}+\epsilon_{23})]^{-1}\leq \epsilon_{23}^{-1}$ (where we used $\epsilon_{i,i+1}<1/24$). Hence, we have $\lvert\tilde{\dist}_3(f(x),f(y))-\tilde{\dist}_2(f_{12}(x),f_{12}(y))\rvert\leq \epsilon_{23}$ and $\lvert \tilde{\dist}_2(f_{12}(x),f_{12}(y))-\tilde{\dist}_1(x,y)\rvert\leq\epsilon_{12}$. Thus, we get:
\begin{align*}
\lvert \tilde{\dist}_3(f(x),f(y))-\tilde{\dist}_1(x,y)\rvert &\leq \epsilon_{12}+\epsilon_{23}\\
&\leq 4(\epsilon_{12}+\epsilon_{23}).
\end{align*}
Using the same argument for $g$, we can conclude that point (iii) of Definition \ref{def:equivtopo} is satisfied.\\
Now, let $x\in\tilde{X}$ and observe that $\tilde{\dist}_2(g_{23}f_{23}f_{12}(x),f_{12}(x))\leq\epsilon_{23}\leq\epsilon_{12}^{-1}$. Hence, we have:
\begin{align*}
\tilde{\dist}_1(gf(x),x)&\leq\tilde{\dist}_{1}(x,g_{12}f_{12}(x))+\tilde{\dist}_1(g_{12}g_{23}f_{23}f_{12}(x),g_{12}f_{12}(x))\\
&\leq \epsilon_{12}+\epsilon_{23}+\epsilon_{12}\\
&\leq 4(\epsilon_{12}+\epsilon_{23}).
\end{align*}
Arguing the same way for $f\circ g$, we can conclude that point (iv) of Definition \ref{def:equivtopo} is satisfied.\\
Now, let us show that point (v) of Definition \ref{def:equivtopo} is also satisfied. First of all, let us found an upper bound on $\dist_{P}^{[2(\epsilon_{12}+\epsilon_{23})]^{-1}}(f_*\tilde{\m}_1,f_{{23}_*}\tilde{\m}_2)$. Let $A$ be a closed subset of $\overline{B}_3([2(\epsilon_{12}+\epsilon_{23})]^{-1})$. We easily get that $f_{23}^{-1}(A)\subset\overline{B}_2(\epsilon_{12}^{-1})$, which implies:
$$
f_{*}\tilde{\m}_1(A)=f_{{23}_*}(f_{{12}_*}\tilde{\m}_1)(A)\leq\tilde{\m}_2((f_{23}^{-1}(A))^{\epsilon_{12}})+\epsilon_{12}.
$$
Then, note that we have $(f_{23}^{-1}(A))^{\epsilon_{12}}\subset f_{23}^{-1}((A)^{\epsilon_{12}+\epsilon_{23}})$. Therefore, we have $f_{*}\tilde{\m}_1(A)\leq f_{{23}_*}\tilde{\m}_2(A^{\epsilon_{12}+\epsilon_{23}})+\epsilon_{12}$. Doing the same in the opposite direction gives us:
$$
\dist_{P}^{[2(\epsilon_{12}+\epsilon_{23})]^{-1}}(f_*\tilde{\m}_1,f_{{23}_*}\tilde{\m}_2)\leq\epsilon_{12}+\epsilon_{23}.
$$
Finally, observe that $[4(\epsilon_{12}+\epsilon_{23})]^{-1}\leq [2(\epsilon_{12}+\epsilon_{23})]^{-1}-({\epsilon_{12}+2\epsilon_{23}})$. Thus, applying inequality \ref{equiv3} of Remark \ref{rem:prokhorov}, we have:
\begin{align*}
\dist_{P}^{[4(\epsilon_{12}+\epsilon_{23})]^{-1}}(f_*\tilde{\m}_1,\tilde{\m}_3)&\leq\dist_{P}^{[2(\epsilon_{12}+\epsilon_{23})]^{-1}-({\epsilon_{12}+2\epsilon_{23}})}(f_*\tilde{\m}_1,\tilde{\m}_3)\\
&\leq \dist_{P}^{[2(\epsilon_{12}+\epsilon_{23})]^{-1}}(f_{{23}_*}\tilde{\m}_2,\tilde{\m}_3)+\dist_{P}^{[2(\epsilon_{12}+\epsilon_{23})]^{-1}}(f_{{23}_*}\tilde{\m}_2,f_*\tilde{\m}_1)\\
&\leq \epsilon_{12}+2\epsilon_{23}\\
&\leq 4(\epsilon_{12}+\epsilon_{23}).
\end{align*}
Applying the same argument, we also get $\dist_{P}^{[4(\epsilon_{12}+\epsilon_{23})]^{-1}}(g_*\tilde{\m}_3,\tilde{\m}_1)\leq4(\epsilon_{12}+\epsilon_{23})$. Therefore, point (v) of Definition \ref{def:equivtopo} is satisfied.\\
This concludes the proof of the modified triangle inequality \ref{trig}.\\

\proofpart{II}{Hausdorff uniform structure}

Let $\mathcal{B}\coloneqq\{\{\Dpeq\leq2^{-n}\},n\in\N\}\subset\mathcal{P}\Big(\M_{0,N}^{\mathrm{p,eq}}(\tilde{X})\times\M_{0,N}^{\mathrm{p,eq}}(\tilde{X})\Big)$. Thanks to the fact that $\Dpeq$ is well defined on $\M_{0,N}^{\mathrm{p,eq}}(\tilde{X})$, symmetric, nonnegative, and satisfies the modified triangle inequality \ref{trig}, we can easily check the axioms introduced p.141 of \cite{Bourbaki_Topo-Generale_07}. Therefore, $\mathcal{B}$ is a fundamental system of neighborhood of a uniform structure on $\M_{0,N}^{\mathrm{p,eq}}(\tilde{X})$.\\

Let us now show that the uniform structure is Hausdorff. Let $\tilde{\mathcal{X}}_i=(\tilde{X},\tilde{\dist}_i,\tilde{\m}_i,\tilde{\ast}_i)\in\kR_{0,N}^{\mathrm{p,eq}}(\tilde{X})$ ($i\in\{1,2\}$) such that $\Dpeq(\tilde{\mathcal{X}}_1,\tilde{\mathcal{X}}_2)=0$. We need to prove that $\tilde{\mathcal{X}}_1$ and $\tilde{\mathcal{X}}_2$ are equivariantly isomorphic.\\
First, since $\Dpeq(\tilde{\mathcal{X}}_1,\tilde{\mathcal{X}}_2)=0$, there is a sequence of equivariant pointed $\epsilon_n$-isometries $(f_n,g_n,\phi_n)$ such that $\epsilon_n\to0$. Let us fix a countable dense subset $\mathcal{D}_1$ in $(\tilde{X},\tilde{\dist}_1)$ and observe that for every $x\in\mathcal{D}_1$, the sequence $\{f_n(x)\}$ is bounded. Therefore, applying Cantor's diagonal argument, we can assume that there is a map $f\colon \mathcal{D}_1\to (\tilde{X},\tilde{\dist}_2)$ such that, for every $x\in\mathcal{D}_1$, we have $f_n(x)\to f(x)$. Then, observe that given $x,y\in\mathcal{D}_1$, and $n$ large enough, we have $\lvert\tilde{\dist}_2(f_n(x),f_n(y))-\tilde{\dist}_1(x,y)\rvert\leq \epsilon_n\to 0$. Therefore, $f$ is an isometric embedding, hence, can be extended in a unique way into an isometric embedding $f\colon (\tilde{X},\tilde{\dist}_1)\to(\tilde{X},\tilde{\dist}_2)$. Then, it is not hard to prove that, given any sequence $\{x_n\}$ in $\tilde{X}$ converging to $x\in\tilde{X}$, we have $f_n(x_n)\to f(x)$. Now, we can apply the same procedure to the sequence $\{g_n\}$, and get an isometric embedding $g\colon(\tilde{X},\tilde{\dist}_2)\to(\tilde{X},\tilde{\dist}_1)$ such that, given any sequence $\{x_n\}$ in $\tilde{X}$ converging to $x\in\tilde{X}$, we have $g_n(x_n)\to g(x)$. In particular, given $x\in\tilde{X}$, we have $\lim_{n\to\infty}\tilde{\dist}_1(g_n(f_n(x)),x)=0=\tilde{\dist}_1(g(f(x)),x)$. Therefore, $f$ and $g$ are respectively inverse to each other. Moreover, we have $f(\tilde{*}_1)=\lim_{n\to\infty}f_n(\tilde{*}_1)=\tilde{*}_2$. The same argument gives us $g(\tilde{*}_2)=\tilde{*}_1$. Hence, $f\colon (\tilde{X},\tilde{\dist}_1,\tilde{*}_1)\to(\tilde{X},\tilde{\dist}_2,\tilde{*}_2)$ is an isomorphism of pointed metric space. Now, observe that $\dist_{\mathcal{P}}^{\epsilon_n^{-1}}(f_{n_*}\tilde{\m}_1,\tilde{\m}_2)\leq\epsilon_n\to 0$. Hence, thanks to Remark \ref{rem:prokhorov}, $\{f_{n_*}\tilde{\m}_1\}$ converges to $\tilde{\m}_2$ in the weak-$*$ topology. Let us show that $f_*\tilde{\m}_1=\tilde{\m}_2$. We fix $R>0$ and $h\in\Class^0(\tilde{X})$ such that $\Spt(h)\subset\overline{B}_2(R)$. Observe that we have $\int_{\tilde{X}}h\dist\tilde{\m}_2=\lim_{n\to\infty}\int_{\tilde{X}}h\dist f_{n_*}\tilde{\m}_1$. However, $h\circ f_n$ is point-wise converging to $h\circ f$. Also, whenever $n$ is large enough, we have $\Spt(h\circ f_n)\subset \overline{B}_1(2R)$. Hence, applying the dominated convergence theorem, we obtain:
$$
\int_{\tilde{X}}h\dist\tilde{\m}_2=\lim_{n\to\infty}\int_{\tilde{X}}h\dist f_{n_*}\tilde{\m}_1=\int_{\tilde{X}}h\dist f{_*}\tilde{\m}_1.
$$
Therefore, since $\tilde{\m}_2$ and $f_{*}\tilde{\m}_1$ are Radon measures, we necessarily have $\tilde{\m}_2=f_{*}\tilde{\m}_1$. Finally, given $\gamma\in\overline{\pi}_1(X)$ and $x\in\tilde{X}$, we have:
$$
p(f(x))=\lim_{n\to\infty}p(f_n(x))=\lim_{n\to\infty}p(f_n(\gamma x))=p(f(\gamma x)).
$$
Thus, $p(f\gamma f^{-1})=p$. We define $\phi\in\Iso(\overline{\pi}_1(X))$ by $\phi(\gamma)\coloneqq f\gamma f^{-1}$, which satisfies $f(\gamma x)=\phi(\gamma)f(x)$.\\
Hence, we can conclude that $\tilde{\mathcal{X}}_1$ is equivariantly isomorphic to $\tilde{\mathcal{X}}_2$.\\

\proofpart{III}{Metrizable uniform structure}

We have seen that $\Dpeq$ induces a Hausdorff uniform structure on $\M_{0,N}^{\mathrm{p,eq}}(\tilde{X})$. Moreover, $\mathcal{B}$ is a countable system of fundamental neighborhoods for this uniform structure. Therefore, thanks to Proposition 2 p.126 of \cite{Bourbaki_Topo-Generale_II_07}, there exists a distance $\dist\colon\M_{0,N}^{\mathrm{p,eq}}(\tilde{X})\times\M_{0,N}^{\mathrm{p,eq}}(\tilde{X})\to[0,+\infty]$ such that $\dist$ induces the same uniform structure as $\Dpeq$. Observe that we can assume, without loss of generality, that $\dist$ is finite (replacing $\dist$ by $\min\{1,\dist\}$ if necessary), which concludes the proof.
\end{proof}

\section*{Declaration about conflict of interests}
 On behalf of all authors, the corresponding author states that there is no conflict of interest.

\bibliographystyle{abbrv}
\bibliography{biblio.bib}
\end{document}